\documentclass{article}

\usepackage{microtype}
\usepackage{graphicx}
\usepackage{booktabs} % for professional tables

\usepackage{hyperref}

\usepackage[accepted]{icml2024}

\usepackage{amsmath}
\usepackage{amssymb}
\usepackage{mathtools}
\usepackage{amsthm}

\usepackage[capitalize,noabbrev]{cleveref}

\usepackage{lipsum}
\usepackage{amsfonts}
\usepackage{graphicx}
\usepackage{epstopdf}
\usepackage{algorithmic}
\ifpdf
  \DeclareGraphicsExtensions{.eps,.pdf,.png,.jpg}
\else
  \DeclareGraphicsExtensions{.eps}
\fi

\usepackage{bm}

\usepackage{mathtools}
\usepackage{cancel} 
\usepackage{xparse}
\usepackage{xfrac} % sfrac

\newcommand{\cH}{\mathcal{H}}

\newcommand{\cS}{\mathcal{S}}

\newcommand{\cZ}{\mathcal{Z}}

\newcommand{\bbR}{\mathbb{R}}

\newcommand{\bA}{\bm{A}}

\newcommand{\bC}{\bm{C}}
\newcommand{\bD}{\bm{D}}

\newcommand{\bI}{\bm{I}}

\newcommand{\bM}{\bm{M}}
\newcommand{\bN}{\bm{N}}

\newcommand{\bP}{\bm{P}}
\newcommand{\bQ}{\bm{Q}}
\newcommand{\bR}{\bm{R}}

\newcommand{\bT}{\bm{T}}

\newcommand{\ba}{\bm{a}}

\newcommand{\bc}{\bm{c}}
\newcommand{\bd}{\bm{d}}

\newcommand{\bs}{\bm{s}}

\newcommand{\bv}{\bm{v}}

\newcommand{\bx}{\bm{x}}
\newcommand{\by}{\bm{y}}
\newcommand{\bz}{\bm{z}}

\NewDocumentCommand{\norm}{mG{2}}{\big\|#1\big\|_{#2}}

\DeclareMathOperator{\rank}{rank}
\DeclareMathOperator{\Span}{span}
\newcommand{\argmin}{\mathop{\rm argmin}}
\newcommand{\argmax}{\mathop{\rm argmax}}

\NewDocumentCommand{\seqp}{mG{n}}{{#1}_1-\cdots+ {#1}_{#2}}
\NewDocumentCommand{\seqm}{mG{n}}{{#1}_1-\cdots- {#1}_{#2}}

\newcommand{\myparagraph}[1]{\textbf{#1.}}

\newtheorem{theorem}{Theorem}
\newtheorem{proposition}{Proposition}
\newtheorem{corollary}{Corollary}
\newtheorem{lemma}{Lemma}

\theoremstyle{remark}
\newtheorem{remark}{Remark}

\theoremstyle{definition}
\newtheorem{definition}{Definition}
\newtheorem{example}{Example}
\newtheorem{assumption}{Assumption}

\newcommand{\ra}[1]{\renewcommand{\arraystretch}{#1}}
\usepackage{tikz}
\usepackage{pgfplots}
\usepackage{enumitem}
\usepackage{booktabs}
\usepackage{multirow,rotating}
\usepackage{subcaption}
\crefformat{equation}{(#2#1#3)}

\usepackage[textsize=tiny]{todonotes}

\icmltitlerunning{On Optimal Stepsizes of Block Gradient Descent for Least-Squares}

\begin{document}

\twocolumn[
\icmltitle{Block Acceleration Without Momentum:\\ On Optimal Stepsizes of Block Gradient Descent for Least-Squares}



\icmlsetsymbol{equal}{*}

\begin{icmlauthorlist}
\icmlauthor{Liangzu Peng}{upenn,equal}
\icmlauthor{Wotao Yin}{ali}
\end{icmlauthorlist}

\icmlaffiliation{upenn}{Department of Electrical and Systems Engineering, University of Pennsylvania, Philadelphia, PA 19104, USA}
\icmlaffiliation{ali}{Decision Intelligence Lab, DAMO Academy, Alibaba Group US, Bellevue, WA 98004, USA}

\icmlcorrespondingauthor{Liangzu Peng}{lpenn@seas.upenn.edu}
\icmlcorrespondingauthor{Wotao Yin}{wotao.yin@alibaba-inc.com}

\icmlkeywords{Block Coordinate Descent, Alternating Projection, Least-Squares}

\vskip 0.3in
]



\printAffiliationsAndNotice{\icmlEqualContribution} % otherwise use the standard text.

\begin{abstract}
Block coordinate descent is a powerful algorithmic template suitable for big data optimization. This template admits a lot of variants including block gradient descent (BGD), which performs gradient descent on a selected block of variables, while keeping other variables fixed. For a very long time, the stepsize for each block has tacitly been set to one divided by the block-wise Lipschitz smoothness constant, imitating the vanilla stepsize rule for gradient descent (GD). However, such a choice for BGD has not yet been able to theoretically justify its empirical superiority over GD, as existing convergence rates for BGD have worse constants than GD in the deterministic cases. 

To discover such theoretical justification, we set up a simple environment where we consider BGD applied to least-squares with two blocks of variables. Assuming the data matrix corresponding to each block is orthogonal, we find optimal stepsizes of BGD in closed form, which provably lead to asymptotic convergence rates twice as fast as GD with Polyak's momentum; this means, under that orthogonality assumption, one can accelerate BGD by just tuning stepsizes and without adding any momentum. An application that satisfies this assumption is \textit{generalized alternating projection} between two subspaces, and applying our stepsizes to it improves the prior convergence rate that was once claimed, slightly inaccurately, to be optimal. The main proof idea is to minimize, in stepsize variables, the spectral radius of a matrix that controls convergence rates.

% Moreover, we apply our stepsizes to generalized alternating projection between two subspaces, improving the prior convergence rate that was once claimed, slightly inaccurately, to be optimal. The main proof idea is to minimize, in stepsize variables, the spectral radius of a matrix that controls convergence rates.

% Our technical devices include assuming block-wise orthogonality and minimizing the spectral radius of a matrix that controls convergence rates.

% We revisit three classic algorithms: alternating projection between two linear subspaces, the two-block Kaczmarz method for solving linear equations, and two-block Gauss-Seidel method for least-squares. At each iteration, all three algorithms solve some optimization problem to optimality, thereby inducing the common wisdom that none of them requires any stepsize parameters. On the contrary, we show that empowering these methods with two appropriately chosen stepsizes provably accelerates their vanilla versions. We call this two-block acceleration (TBA). For alternating projection, TBA alleviates the well-known negative result that vanilla alternating projection might be arbitrarily slow. In the case of the two-block Gauss-Seidel method for least-squares, TBA accelerates without momentum, and in some sense beats momentum: We prove TBA is at least twice as fast as the heavy ball method. Overall, the paper provides the first two-step analysis for these problems, as opposed to the prevalent single-step convergence proof in the literature.
\end{abstract}

\section{Introduction}
Block coordinate descent refers to a family of algorithms selecting and updating one block of variables at a time. In the span of more than six decades since its early appearance \citep{Hildreth-1957}, many variants of block coordinate descent have been proposed, analyzed, and recently tested on big data scenarios  \citep{Nesterov-SIAM-J-Opt2012,Xu-SIAM-J-IS2013,Wright-MP2015,Shi-arXiv2016v2,Lin-ICML2023b}. Despite the abundance of these exciting developments, one might ponder when one should prefer block coordinate descent to vanilla gradient descent (GD). Intuitively, block coordinate descent might be advantageous if the optimization variables admit a natural partition into blocks, which is the case in many applications \citep{Peng-arXiv2023}; or if the variable dimension is too high to fit the memory, in which case one is forced to optimize in a block-wise manner \citep{Nesterov-SIAM-J-Opt2012}; or if the given optimization problem is \textit{coordinate-friendly} \citep{Peng-AMSA2016}, meaning that minimizing over one block of coordinates while fixing others is computationally easy. 

Developing convergence theory to support the above intuition, however, is challenging. For example, consider \textit{block gradient descent} (BGD), a method that runs GD at every iteration on a selected block of variables. Nesterov was concerned that, if selecting blocks in a \textit{greedy} fashion (e.g., using the famous \textit{Gauss–Southwell} rule), then, following standard reasoning, one obtains a bound on convergence rates that might have a worse constant than that of GD \citep{Nesterov-SIAM-J-Opt2012}. This concern led him to a BGD variant that randomly selects blocks, whose convergence rate is proved to be better than GD. While theoretically appealing and having attracted a sequence of follow-up works, such a randomized variant ensures convergence only in expectation; this is perhaps why it is empirically slower\footnote{It is also shown that this randomized variant is in many cases empirically slower than the greedy version \citep{Nutini-ICML2015}.} for the least-squares problem than \textit{cyclic BGD}, a BGD variant that selects blocks in a cyclic fashion, see Tables 3.1 \& 3.2 of \citet{Beck-SIAM-J-2013}. We focus on the cyclic rule, and for short, we write BGD for cyclic BGD in the sequel. 

In the deterministic setting, the \textit{BGD versus GD} dilemma persists: BGD is empirically faster, but its current bound on (deterministic) convergence rates has, in general, a worse constant than GD; see Remark 3.3 of \citet{Beck-SIAM-J-2013}. Furthermore, this bittersweet dilemma manifests itself again if one extends, in a direct way, BGD and GD into their proximal versions (e.g., compare Theorems 10.15 and 11.18 of \citet{Beck-OptBook2017}), or into their accelerated versions (e.g., compare  Theorem 4.2 of \citet{Beck-SIAM-J-2013} and Theorem 2.2.2 of \citet{Nesterov-2004}).

How does one even reconcile these? It is now the case that there is some sub-optimality in the existing analysis of BGD methods, but it might also be the case that, after decades of development, making improvements is difficult.

%\subsection{Our Approach and Results: Overview}
Our approach to making progress features a return to the very basic setting, where we minimize the arguably simplest, very well-studied objective, \textit{least-squares}, 
\begin{equation}\label{eq:LS}
	\min_{\bx \in \bbR^n } F(\bx), \quad F(\bx):= \frac{1}{2}\cdot \| \bA \bx - \by \|_2^2.
\end{equation}
Here we assume the $m\times n$ matrix $A$ is of full column rank $n$ (necessarily, $m\geq n$). 
Instead of with multiple, we consider BGD with only two blocks. Specifically, we write  $\bx:=[\bx_1;\bx_2]$, $\bA:=[\bA_1\  \bA_2]$, and write $F$ and \cref{eq:LS} into
\begin{equation}\label{eq:LS2}
	\min_{\bx_1 \in \bbR^{n_1}, \bx_2 \in \bbR^{n_2} }  \frac{1}{2}\cdot \| \bA_1 \bx_1 + \bA_2 \bx_2  - \by \|_2^2.
\end{equation}
We run GD on each block in an alternating fashion with constant stepsizes $\gamma_1$ and $\gamma_2$ to minimize the least-squares objective $F$ and search for the global minimizer $\bx^*$:
\begin{equation}\label{eq:BGD}
	\begin{split}
		\bx^+_1 &= \bx_1 - \gamma_1 \cdot \nabla_{\bx_1} F(\bx) \\ 
		&= \bx_1 - \gamma_1 \cdot 
		\bA_1^\top ( \bA_1 \bx_1 + \bA_2 \bx_2  - \by ), \\ 
		\bx^+_2 &= \bx_2 - \gamma_2 \cdot \nabla_{\bx_2} F([\bx_1^+;\bx_2]) \\ 
		&= \bx_2 - \gamma_2 \cdot \bA_2^\top ( \bA_1 \bx_1^+ + \bA_2 \bx_2  - \by ).  
	\end{split} \tag{BGD}
\end{equation}
We denote by $\bx^+:=[\bx^+_1;\bx^+_2]$ and $\bx:=[\bx_1;\bx_2]$ the two consecutive iterates of \ref{eq:BGD}. It is clear that the convergence (rate) of \ref{eq:BGD} depends on stepsizes $\gamma_1,\gamma_2$ and that the fastest convergence of \ref{eq:BGD} is attained only if we set the stepsizes to be ``optimal''. Let the word ``optimal'' be vague for the moment, and we report the following result:

\begin{theorem}[Informal]\label{theorem:informal}
	Suppose \cref{assumption:BWO} below holds. Run  GD with Polyak's momentum (i.e., the heavy ball method) and \ref{eq:BGD}, respectively, with their ``optimal'' stepsizes. \ref{eq:BGD} is twice as fast as the heavy ball method (HB).
\end{theorem}
\begin{assumption}\label{assumption:BWO}
	$\bA_1^\top \bA_1$ and $\bA_2^\top \bA_2$ are identity matrices, and $\bA_2^\top \bA_1\neq 0$.
\end{assumption}
% We argue that \cref{assumption:BWO} is valid and reasonable, and we will justify this argument in \cref{section:assumption}.
The rest of the paper is organized as follows. In \cref{section:background} we establish notations and quantify optimal stepsizes. In \cref{section:assumption} we elaborate on \cref{assumption:BWO} and argue that \cref{assumption:BWO} is valid and reasonable. In \cref{section:corollary-main-result} we present optimal stepsizes for \ref{eq:BGD} in comparison to the heavy ball method, from which \cref{theorem:informal} follows. In \cref{section:related-work} we make connections to prior works, and in particular in \cref{subsection:review-GAP} we connect \ref{eq:BGD} to alternating projection. Then, in \cref{section:close-gap}, we apply our results to the problem of alternating projection, resulting in improvements over prior works. In \cref{section:experiments}, we perform basic experiments to verify our theory, and in \cref{section:conclusion} we conclude the paper.

% In \cref{section:spectrum,section:simplified-case,section:general-case}, we derive optimal stepsizes of \ref{eq:BGD} under \cref{assumption:BWO}, therefore proving \cref{theorem:informal}.

\section{Quantifying Optimality}\label{section:background} 
%Here we provide some background and necessary concepts. 
In \cref{subsection:GD-HB} we review gradient descent (GD) and the heavy ball method (HB), and clarify what ``optimal'' stepsizes mean for them. In \cref{subsection:BGD} we characterize optimal stepsizes for \ref{eq:BGD} and state \cref{theorem:informal} formally.

% is similar and we characterize optimal stepsizes for \ref{eq:BGD}. % A side product in \cref{subsection:BGD} is a formal statement of \cref{theorem:informal}.

% block gradient descent (BGD), and clarify what ``optimal'' stepsizes mean for them.  % In \cref{section:assumption}, we elaborate on \cref{assumption:BWO}. 

\subsection{Gradient Descent and Heavy Ball: A Review}\label{subsection:GD-HB} 
Recall that vanilla gradient descent applied to least-squares \cref{eq:LS} comes with a stepsize $\gamma>0$ and updates $\bx$ via
\begin{equation}\label{eq:GD}
    \bx^+ = \bx - \gamma \cdot \nabla F(\bx) = \bx - \gamma \cdot 
 \bA^\top ( \bA \bx  - \by ). \tag{GD}
\end{equation}
The two consecutive iterates $\bx^+$ and $\bx$ of \ref{eq:GD} satisfy
\begin{equation*}
     \bx^+ - \bx^* = (\bI - \gamma \bA^\top \bA) \cdot (\bx - \bx^*), 
\end{equation*}
so the rate of convergence to $\bx^*$ is dictated by the spectrum of $\bI - \gamma \bA^\top \bA$. Thus we wish to find a stepsize $\gamma$ to minimize the \textit{spectral radius} $\rho(\bI - \gamma \bA^\top \bA)$, where $\rho(\cdot)$ denotes the maximum magnitude of eigenvalues of a matrix. Explicitly:
\begin{equation*}
    \rho_{\textnormal{GD}}(\gamma):=\rho(\bI - \gamma \bA^\top \bA) = \max_{i=1,\dots,n} |\lambda_i(\bI - \gamma \bA^\top \bA)|.
\end{equation*} 
Here we used $\lambda_i(\cdot)$ to mean the $i$-th largest eigenvalue of a matrix. A folklore fact in optimization is that $\rho_{\textnormal{GD}}(\gamma)$ is minimized at $\gamma= \gamma^*:= \frac{2}{\lambda_1(\bA^\top \bA) + \lambda_n(\bA^\top \bA)  }$, that is
\begin{equation*}
    \begin{split}
        \rho_{\textnormal{GD}}^*:=\rho_{\textnormal{GD}}(\gamma^*) &= \min_{\gamma>0} \rho_{\textnormal{GD}}(\gamma) \\
        &= \frac{\lambda_1(\bA^\top \bA) - \lambda_n(\bA^\top \bA)}{\lambda_1(\bA^\top \bA) + \lambda_n(\bA^\top \bA)}.
    \end{split}
\end{equation*}
We say $\gamma^*$ is the optimal (constant) stepsize of \ref{eq:GD} applied to least-squares as it minimizes the spectral radius $\rho_{\textnormal{GD}}(\gamma)$.

In a similar style we review the heavy ball method \citep{Polyak-1964}. This method calculates the current iterate $\bx^{++}$ using two previous points $\bx^+$ and $\bx$, that is 
\begin{equation}\label{eq:HB}
    \bx^{++} = \bx^+ - \alpha\cdot \nabla F(\bx^+) + \beta \cdot (\bx^+ - \bx),  \tag{HB}
\end{equation}
where $\nabla F(\bx^+)$ is the gradient at $\bx^+$, i.e., $\nabla F(\bx^+)=\bA^\top ( \bA \bx^+  - \by )$. Compare this with \ref{eq:GD}. Note that \ref{eq:HB} has two parameters:  stepsize $\alpha>0$ and momentum coefficient $\beta\geq 0$. It is known that the iterates of \ref{eq:HB} satisfy
\begin{equation*}
    \begin{split}
        \begin{bmatrix}
        \bx^+ -\bx^* \\ 
        \bx^{++} - \bx^*
    \end{bmatrix} = 
    \bN(\alpha,\beta) \cdot \begin{bmatrix}
        \bx -\bx^* \\ 
        \bx^{+} - \bx^*
    \end{bmatrix}, \textnormal{where} \\
    \bN(\alpha,\beta):=\begin{bmatrix}
        0 & \bI \\ 
        -\beta \bI & (1+\beta)\bI - \alpha \bA^\top \bA
    \end{bmatrix}.
    \end{split}
\end{equation*}
Similarly, we call a stepsize $(\alpha,\beta)$ optimal if it minimizes  $\rho\big( \bN(\alpha,\beta) \big)=:\rho_{\textnormal{HB}}(\alpha,\beta)$. Let $\rho_{\textnormal{HB}}^*$ be the minimum of 
\begin{equation*}
    \min_{\alpha>0,\beta\geq 0} \rho_{\textnormal{HB}}(\alpha,\beta),
\end{equation*}
attained at $(\alpha^*,\beta^*)$. One can then prove \citep{Polyak-1964}:
\begin{equation*}
    \begin{split}
        &\ \rho_{\textnormal{HB}}^*= \frac{\sqrt{\lambda_1(\bA^\top \bA)} - \sqrt{\lambda_n(\bA^\top \bA)}}{\sqrt{\lambda_1(\bA^\top \bA)} + \sqrt{\lambda_n(\bA^\top \bA)}}, \\ 
        &\ \alpha^*= \left( \frac{2}{\sqrt{\lambda_1(\bA^\top \bA)} + \sqrt{\lambda_n(\bA^\top \bA)}} \right)^2,  \\ 
        &\ \beta^*= \left( \frac{\sqrt{\lambda_1(\bA^\top \bA)} - \sqrt{\lambda_n(\bA^\top \bA)}}{\sqrt{\lambda_1(\bA^\top \bA)} + \sqrt{\lambda_n(\bA^\top \bA)}} \right)^2.
    \end{split}
\end{equation*}
\subsection{Block Gradient Descent}\label{subsection:BGD}
The phrase ``optimal stepsizes''  of \ref{eq:BGD} bears a similar meaning to that of \ref{eq:HB}. To concretize this, we first rewrite the updates of \ref{eq:BGD} (let $\bI$ be the $n\times n$ identity matrix):
\begin{lemma}\label{lemma:error-decrease}
The iterates of \ref{eq:BGD} satisfy $\bx^+ - \bx^*= \bM(\gamma_1,\gamma_2) \cdot (\bx - \bx^*)$, with $\bM(\gamma_1,\gamma_2)$ defined as
    \begin{equation*}
        \tiny{
        \left(\bI - \begin{bmatrix}
            0 & 0 \\ 
            \gamma_2 \bA_2^\top \bA_1 & \gamma_2 \bA_2^\top \bA_2 
        \end{bmatrix} \Bigg) \Bigg(\bI - \begin{bmatrix}
            \gamma_1 \bA_1^\top \bA_1 & \gamma_1 \bA_1^\top \bA_2 \\ 
            0 & 0
        \end{bmatrix} \right).
        }
    \end{equation*}
\end{lemma}
% \begin{proof}%[Proof of \cref{lemma:error-decrease}]
%     The global minimizer $\bx^*$ of \cref{eq:LS2} satisfies the normal equation $\bA^\top \bA \bx^* = \bA^\top \by$, 
%     which means $\bA_1^\top \bA \bx^* = \bA_1^\top \by$ and $\bA_2^\top \bA \bx^* = \bA_2^\top \by$, so we can write the first equation $\bx^+_1 = \bx_1 - \gamma_1 \cdot \bA_1^\top  \bA (\bx - \bx^*)$ of \ref{eq:BGD} as
%     \begin{equation*}
%         \begin{bmatrix}
%             \bx^+_1 \\ 
%             \bx_2 
%         \end{bmatrix} - \bx^* = \Bigg(\bI - \begin{bmatrix}
%             \gamma_1 \bA_1^\top \bA_1 & \gamma_1 \bA_1^\top \bA_2 \\ 
%             0 & 0
%         \end{bmatrix} \Bigg)  (\bx - \bx^*). 
%     \end{equation*}
%     We finish by re-writing the second equation of \ref{eq:BGD} similarly and combining.
% \end{proof}
We call $(\gamma_1,\gamma_2)$ optimal if it minimizes  $\rho \big( \bM(\gamma_1,\gamma_2) \big)=:\rho_{\textnormal{BGD}}(\gamma_1,\gamma_2)$. Our contribution consists of, under \cref{assumption:BWO}, discovering stepsizes $\gamma_1^*,\gamma_2^*$ that satisfy
\begin{equation}\label{eq:minimize-sr}
    \rho_{\textnormal{BGD}}(\gamma_1^*,\gamma_2^*)= \min_{\gamma_1>0,\gamma_2>0} \rho_{\textnormal{BGD}}(\gamma_1,\gamma_2).% = \min_{\gamma_1>0,\gamma_2>0} \rho \big( \bM(\gamma_1,\gamma_2) \big)
\end{equation}
%This is not a trivial contribution as matrix $\bM(\gamma_1,\gamma_2)$, under  \cref{assumption:BWO}, is much more complicated than $\bN(\alpha,\beta)$, the latter appearing in \ref{eq:HB}. That said, 
With $\rho_{\textnormal{BGD}}^*:=\rho_{\textnormal{BGD}}(\gamma_1^*,\gamma_2^*)$, we  restate \cref{theorem:informal} below: 
% As we proceed, we will gradually gain insights into the structure of $\bM(\gamma_1,\gamma_2)$, and eventually we will see the precise expression of $\gamma_1^*,\gamma_2^*$, and $\rho_{\textnormal{BGD}}(\gamma_1^*,\gamma_2^*)=:\rho_{\textnormal{BGD}}^*$. At this stage, we can at least restate \cref{theorem:informal} as follows: 
\setcounter{theorem}{0}
\begin{theorem}\label{theorem:formal}
    \cref{assumption:BWO} implies $\rho_{\textnormal{BGD}}^* \leq (\rho_{\textnormal{HB}}^*)^2$. %  < (\rho_{\textnormal{GD}}^*)^2.$
\end{theorem}
% \begin{equation*}
%     \text{\cref{assumption:BWO} implies\ \ }  \rho_{\textnormal{BGD}}^* \leq (\rho_{\textnormal{HB}}^*)^2  < (\rho_{\textnormal{GD}}^*)^2.
% \end{equation*}
Vaguely put, a smaller spectral radius implies faster convergence. And the comparison is fair: All stepsizes are chosen to be optimal, minimizing their respective spectral radii; and all methods, namely \ref{eq:BGD}, \ref{eq:GD}, and \ref{eq:HB}, have comparable costs at each iteration. Moreover, $\rho_{\textnormal{BGD}}^* \leq (\rho_{\textnormal{HB}}^*)^2$ implies that, with optimal stepsizes, \ref{eq:BGD} is asymptotically at least twice as fast as \ref{eq:HB}.  Finally, we emphasize it is the presence of \cref{assumption:BWO} that remains to be interrogated.

\section{On \cref{assumption:BWO}}\label{section:assumption}
In \cref{section:assumption1-AP} We argue \cref{assumption:BWO} is reasonable, as it relates \ref{eq:BGD} to (generalized) \textit{alternating projection between two linear subspaces} in such a way that our results under \cref{assumption:BWO} directly applies and improves existing bounds of \citet{Falt-CDC2017} on convergence rates of \textit{generalized alternating projection}. Then, in \cref{subsection:role-of-assumption1}, we argue \cref{assumption:BWO} is also essential as it makes analysis possible.

\subsection{Justifying   \cref{assumption:BWO}}\label{section:assumption1-AP}
Here we discuss why \cref{assumption:BWO} is reasonable. Note that it requires the orthogonality of $\bA_j$ ($j=1,2$), that is $\bA_j^\top \bA_j=\bI_j$, where $\bI_j$ is the $n_j\times n_j$ identity matrix. First of all, one can always realize this assumption by orthogonalizing $\bA_1,\bA_2$ respectively (even though it entails some computational costs). Secondly, and more importantly, this block-wise orthogonality assumption bridges \ref{eq:BGD} and the classic method of alternating projection between two subspaces. To see this, let us consider the following lemma:
% Then observe \cref{eq:BWO->M} is obtained by applying the orthogonality assumption $\bA_1^\top \bA_1=\bI_1, \bA_2^\top\bA_2=\bI_2$, while one can always realize this assumption by orthogonalizing $\bA_1,\bA_2$ respectively. Finally, \cref{assumption:BWO} requires $\bC\neq 0$; this is to sidestep the trivial case where \ref{eq:BGD}, \ref{eq:GD}, and \ref{eq:HB} all converge to $\bx^*$ in just $1$ iteration with appropriate stepsizes.
% To see a connection between \ref{eq:BGD} and the classic method of alternating projection between two subspaces, let us first state the following basic lemma:
\begin{lemma}\label{lemma:BGD-AP}
	The iterates $\bx^+,\bx$ of \ref{eq:BGD} satisfy $\bA(\bx^+ - \bx^*) = \bT(\gamma_1,\gamma_2) \bA (\bx - \bx^*)$, with $\bT(\gamma_1,\gamma_2)$ defined as   
	\begin{equation*}
		\bT(\gamma_1,\gamma_2) := \left( \bI - \gamma_2 \bA_2 \bA_2^\top \right) \left( \bI - \gamma_1 \bA_1 \bA_1^\top \right).
	\end{equation*}
\end{lemma}
\cref{lemma:BGD-AP} implies \ref{eq:BGD} can be viewed as an algorithm that operates on the iterates $\bz^+:=\bA(\bx^+ - \bx^*)$ and $\bz:=\bA (\bx - \bx^*)$. Since $\bA$ is assumed to be full column rank, updating $\bz^+$ is equivalent to updating $\bx^+$. If \cref{assumption:BWO} holds, then $\bI - \bA_j \bA_j^\top$ ($j=1,2$) is an orthogonal projection, the update $\bz^+=\bT(1,1) \bz$ is precisely the vanilla alternating projection method between two subspaces \citep{Von-1951}, and the update $\bz^+=\bT(\gamma_1,\gamma_2) \bz$ can be viewed as \textit{generalized} alternating projection \citep{Falt-CDC2017}. 

With the above reasoning, we can intuitively conclude that studying the convergence of \ref{eq:BGD} under \cref{assumption:BWO} would also yield convergence guarantees for generalized alternating projection. A more detailed treatment from the perspective of alternating projection can be found in \cref{subsection:review-GAP}, where we rectify the slightly misleading claim of \citet{Falt-CDC2017} that their bound was optimal. 

Finally, \cref{assumption:BWO} requires $\bA_2^\top \bA_1\neq 0$; this is to sidestep the trivial case where \ref{eq:BGD}, \ref{eq:GD}, and \ref{eq:HB} all converge to $\bx^*$ in just $1$ iteration with appropriate stepsizes under the block-wise orthogonality assumption.% $\bA_j^\top \bA_j=\bI_j$ ($j=1,2$).

\subsection{The Technical Role of \cref{assumption:BWO}}\label{subsection:role-of-assumption1}
Note that $\bM(\gamma_1,\gamma_2)$ has a sophisticated expression, and so \cref{assumption:BWO} plays the role of simplifying, at least making analyzing \cref{eq:minimize-sr} possible. In particular, \cref{assumption:BWO} immediately simplifies the expression of $\bM(\gamma_1,\gamma_2)$:
\begin{lemma}\label{lemma:BWO->M}
    Define $\bC:=\bA_2^\top \bA_1$. Recall the definition of $\bM(\gamma_1,\gamma_2)$ in \cref{lemma:error-decrease}. \cref{assumption:BWO} implies
    \begin{equation*}%\label{eq:BWO->M}
        \small{
    	\bM(\gamma_1,\gamma_2) = \begin{bmatrix}
    		(1-\gamma_1) \bI_1 & -\gamma_1 \bC^\top  \\ 
    		-\gamma_2 (1- \gamma_1 ) \bC  & (1- \gamma_2) \bI_2  + \gamma_1 \gamma_2 \bC \bC^\top 
    	\end{bmatrix}.
    }
    \end{equation*}
\end{lemma}
We note that \cref{assumption:BWO} does not clean all obstacles, as $\bM(\gamma_1,\gamma_2)$ in \cref{lemma:BWO->M} is still complicated, e.g., it depends on $\gamma_1,\gamma_2$ quadratically; to compare, the corresponding matrix $\bN(\alpha,\beta)$ of \ref{eq:HB} has a linear dependency on $\alpha,\beta$. 

After a simplification through \cref{assumption:BWO}, $\bM(\gamma_1,\gamma_2)$ still maintains an interesting structure that will ultimately facilitate our understanding of \ref{eq:BGD}. For example, $\bM(\gamma_1,\gamma_2)$ depends on $\bC$, and $\bC$ is precisely the matrix that encodes the information about the relationship between the two blocks of variables $\bx_1$ and $\bx_2$. It is by leveraging the spectrum of $\bC$, and therefore of $\bM(\gamma_1,\gamma_2)$, that we will be able to show \ref{eq:BGD} enjoys faster convergence than its competitors \ref{eq:GD} and \ref{eq:HB}. To get prepared for the competition, we use \cref{assumption:BWO} and express the minimum spectral radii $\rho_{\textnormal{HB}}^*$ and $\rho_{\textnormal{GD}}^*$ in terms of the spectrum of $\bC$:
\begin{lemma}\label{lemma:rho-GD-eigs-C} 
Recall $\bA=[\bA_1, \bA_2]\in\bbR^{m\times (n_1+n_2)}$, $\bA$ is full rank, and $\bC=\bA_2^\top \bA_1\in\bbR^{n_2\times n_1}$. Suppose $\rank(\bC)=r$ and  \cref{assumption:BWO} holds. Then $\lambda_1(\bC\bC^\top)\neq 1$ and the maximum and minimum eigenvalues of $\bA^\top \bA$ are $1 + \sqrt{\lambda_1( \bC \bC^\top )}$ and $1 - \sqrt{\lambda_1( \bC \bC^\top )}$, respectively. Hence, the minimum spectral radii of \ref{eq:GD} and \ref{eq:HB} are given as
\begin{equation}\label{eq:min-GD-HB}
     \rho_{\textnormal{GD}}^* = \sqrt{\lambda_1(\bC\bC^\top)}, \   \rho_{\textnormal{HB}}^* = \frac{ 1- \sqrt{1-\lambda_1(\bC\bC^\top)} }{\sqrt{\lambda_1(\bC\bC^\top)}}.
\end{equation}
\end{lemma}

\begin{example}\label{example:BEM}
     With \cref{lemma:rho-GD-eigs-C} and \cref{assumption:BWO}, we have
    \begin{equation*}
        \rho \big( \bM(1, 1) \big) = \lambda_1(\bC\bC^\top) < \sqrt{\lambda_1(\bC\bC^\top)} =  \rho_{\textnormal{GD}}^*, 
    \end{equation*}
    meaning \ref{eq:BEM} converges asymptotically faster than \ref{eq:GD}.
\end{example}
Is it possible that $\gamma_1=\gamma_2=1$ actually minimizes $\rho \big( \bM(\gamma_1,\gamma_2) \big)$? After all, \ref{eq:BEM} already ensures the largest possible decrease of the objective for the present block, while all other stepsizes guarantee less! If that were true, then $\rho \big( \bM(1, 1) \big)$ would not necessarily be smaller than $\rho_{\textnormal{HB}}^*$, let alone our promise $\rho_{\textnormal{BGD}}^* \leq (\rho_{\textnormal{HB}}^*)^2$. The sole hope is that other stepsizes, though sub-optimal for the moment, might be more beneficial in the long run---and if so, we simply need to work harder to find them.

\section{Optimal Stepsizes and Spectral Radius}\label{section:corollary-main-result}
The inequality $\rho_{\textnormal{BGD}}^* \leq (\rho_{\textnormal{HB}}^*)^2$ and hence \cref{theorem:formal} are proved under \cref{assumption:BWO} by deriving a closed-form expression of the minimum value $\rho_{\textnormal{BGD}}^*$ of \cref{eq:minimize-sr} and the corresponding stepsizes. The derivation is summarized below:
% The derivation is summarized below with the values of  $\rho_{\textnormal{BGD}}^*$ and stepsizes given.
\begin{theorem}\label{corollary:main-result}
    Recall $\rho_{\textnormal{BGD}}^*=\rho_{\textnormal{BGD}}(\gamma_1^*,\gamma_2^*)$ in \cref{eq:minimize-sr} and $\bC:=\bA_2^\top \bA_1$. \cref{assumption:BWO} implies the following.

 %   \cref{assumption:BWO}, the minimum $\rho_{\textnormal{BGD}}^*=\rho_{\textnormal{BGD}}(\gamma_1,\gamma_2)$ of \cref{eq:minimize-sr} and the corresponding stepsizes $(\gamma_1,\gamma_2)$:
    \begin{itemize}[parsep=-1pt,topsep=0pt]
        \item If $\bC$ has full rank (i.e., $\rank(\bC)=\min\{n_1,n_2\}$), then
        \begin{equation}\label{eq:rho_BGD_fullrank}
            \rho_{\textnormal{BGD}}^*=\frac{\sqrt{1-\lambda_r(\bC\bC^\top)} - \sqrt{1-\lambda_1(\bC\bC^\top)}  }{\sqrt{1-\lambda_r(\bC\bC^\top)} + \sqrt{1-\lambda_1(\bC\bC^\top)}}.
        \end{equation}
        % and the corresponding stepsizes $(\gamma_1^*,\gamma_2^*)$ are given as
        % \begin{equation*}
        %     \begin{split}
        %         \gamma_1^* &= \frac{ \big(\sqrt{(1 + \nu_r )( 1 + \nu_1 )} + \sqrt{(1 - \nu_r  )( 1 - \nu_1  )} \big)^2 }{ (\nu_1+\nu_r)^2 } , \\ 
        %     \gamma_2^* &= \frac{ \big(\sqrt{(1 + \nu_r )( 1 + \nu_1 )} - \sqrt{(1 - \nu_r  )( 1 - \nu_1  )} \big)^2 }{ (\nu_1+\nu_r)^2 }.
        %     \end{split}
        % \end{equation*}
        % Here $\nu_1:=\sqrt{\zeta_1}$ and $\nu_r:=\sqrt{\zeta_r}$.
        \item If $\bC$ is rank-deficient, then
        \begin{equation}\label{eq:rho_BGD_rankdeficient}
            \rho_{\textnormal{BGD}}^*=\frac{1 - \sqrt{1-\lambda_1(\bC\bC^\top)}}{1+ \sqrt{1-\lambda_1(\bC\bC^\top)} }.
        \end{equation}
        and the corresponding stepsizes $(\gamma_1^*,\gamma_2^*)$ are given as
        \begin{equation}\label{eq:gamma1=gamma2-corollary}
             \gamma_1^*=\gamma_2^*=\frac{2}{1+ \sqrt{1-\lambda_1(\bC \bC^\top)}}.
        \end{equation}
    \end{itemize}
\end{theorem}
\begin{remark}
    If $r=1$, then we have $\lambda_1(\bC\bC^\top)= \lambda_r(\bC\bC^\top)$ and \cref{corollary:main-result} suggests the minimum spectral radius is $0$. When $\bC$ has full rank, the minimizers of $\rho\big( \bM(\gamma_1,\gamma_2) \big)$ have complicated forms, so we do not show them here; they can be found in the proof of \cref{subsection:proof-S11-fullrank}.    
\end{remark}
\begin{remark}\label{remark:rho-BGD<rho-HB}
    The full rank case exhibits a smaller minimum spectral radius than the rank-deficient case, and the minimum spectral radius in the rank-deficient case is precisely  $(\rho_{\textnormal{HB}}^*)^2$, where $\rho_{\textnormal{HB}}^*$ is defined in \cref{eq:min-GD-HB}. Explicitly, we have
\begin{equation*}
    \begin{split}
        \cref{eq:rho_BGD_fullrank}& =\frac{\sqrt{1-\lambda_r(\bC\bC^\top)} - \sqrt{1-\lambda_1(\bC\bC^\top)}  }{\sqrt{1-\lambda_r(\bC\bC^\top)} + \sqrt{1-\lambda_1(\bC\bC^\top)}} \\
    &\leq   \frac{1 - \sqrt{1-\lambda_1(\bC\bC^\top)}}{1+ \sqrt{1-\lambda_1(\bC\bC^\top)} } = \cref{eq:rho_BGD_rankdeficient} = (\rho_{\textnormal{HB}}^*)^2.
    \end{split}
\end{equation*}
    This consolidates \cref{theorem:formal}, and implies, under \cref{assumption:BWO}, \ref{eq:BGD} is eventually twice as fast as \ref{eq:HB}.
\end{remark}
\begin{remark}
    Interestingly, the optimal stepsizes $(\gamma_1^*,\gamma_2^*)$ derived in \cref{eq:gamma1=gamma2-corollary} are the same as the optimal stepsize $\alpha^*$ of \ref{eq:HB}. Indeed, with \cref{assumption:BWO} and \cref{lemma:rho-GD-eigs-C}, we have
    \begin{equation*}
        \begin{split}
            \gamma_1^*= \gamma_2^* &= \frac{2}{1+ \sqrt{1-\lambda_1(\bC \bC^\top)}} \\
            &= \Bigg( \frac{2}{\sqrt{\lambda_1(\bA^\top \bA)} + \sqrt{\lambda_n(\bA^\top \bA)}} \Bigg)^2=\alpha^*.
        \end{split}
    \end{equation*} 
    This offers some intuition as to why \ref{eq:BGD} is \textit{twice} as fast under \cref{assumption:BWO}: \ref{eq:BGD} and \ref{eq:HB} take the same stepsize, but \ref{eq:BGD} takes \textit{two} descent steps---with partial gradients and without momentum---in a single iteration, entailing a cost barely comparable to one descent step of \ref{eq:HB}.
\end{remark}

% \begin{itemize}[wide,parsep=-1pt,topsep=0pt]
%     \item Interestingly, the optimal stepsize $\gamma_1^*$ derived in \cref{theorem:S11-rankdeficient} and the optimal stepsize $\alpha^*$ for \ref{eq:HB} are the same. Indeed, with \cref{assumption:BWO} and \cref{lemma:rho-GD-eigs-C}, we have
%     \begin{equation*}
%         \begin{split}
%             \gamma_1^* &= \frac{2}{1+ \sqrt{1-\lambda_1(\bC \bC^\top)}} \\
%             &= \Bigg( \frac{2}{\sqrt{\lambda_1(\bA^\top \bA)} + \sqrt{\lambda_n(\bA^\top \bA)}} \Bigg)^2=\alpha^*.
%         \end{split}
%     \end{equation*} 
%     This offers some intuition as to why \ref{eq:BGD} is \textit{twice} as fast under \cref{assumption:BWO}: \ref{eq:BGD} and \ref{eq:HB} take the same stepsize, but \ref{eq:BGD} takes \textit{two} descent steps---with partial gradients and without momentum---in a single iteration, entailing a cost barely comparable to one descent step of \ref{eq:HB}.
% \end{itemize}

\section{Related Work}\label{section:related-work}
\subsection{BGD Versus GD and HB}
\cref{theorem:informal} shows that \ref{eq:BGD} can be faster than the heavy ball method (HB), and the latter has thus far been one of the theoretically fastest variants of accelerated GD for least-squares. Interestingly, though, \ref{eq:BGD} has no momentum at all! In hindsight, our justification is that \ref{eq:BGD} is a \textit{two-step} method in the sense of \citet{Polyak-1987}: It has two stepsizes $\gamma_1,\gamma_2$, comparable to the stepsize and momentum coefficient in HB, so, similarly to HB, choosing the stepsizes $\gamma_1,\gamma_2$ in an ``optimal'' way results in acceleration. While several methods that accelerate without momentum exist (e.g., Young's method \citep{Young-JMP1953}, GD with cyclic stepsizes, see \cref{subsection:cyclic-stepsizes}), our approach is unique in the sense that it achieves acceleration by operating on \textit{partial} gradients without momentum in a cyclic, block-wise manner.

\subsection{BGD Versus GD with Cyclic Stepsizes}\label{subsection:cyclic-stepsizes}
There has been some recent interest in using GD with stepsizes selected in a cyclic fashion from a sequence of $M$ positive numbers \citep{Smith-WACV2017,Oymak-SPL2021,Baptiste-AISTATS2022,Grimmer-arXiv2023}; see also \textit{fractal} stepsizes of \citet{Agarwal-ICML2021} and silver stepsizes of \citet{Altschuler-arXiv2023,Altschuler-arXiv2023b}. While \citet{Grimmer-arXiv2023} and \citet{Altschuler-arXiv2023,Altschuler-arXiv2023b} considered GD in the more general setting of (strongly-)convex smooth optimization, more relevant to ours is the work of \citet{Oymak-SPL2021,Baptiste-AISTATS2022}. \citet{Oymak-SPL2021} analyzed  GD with cyclic stepsizes in the least-squares context and \citet{Baptiste-AISTATS2022} analyzed HB with cyclic stepsizes. Both \citet{Oymak-SPL2021} and \citet{Baptiste-AISTATS2022} assume the spectrum of $\bA^\top\bA$ is \textit{clustered} into two or more disjoint intervals; this is different from \cref{assumption:BWO}, so the results of \citet{Oymak-SPL2021,Baptiste-AISTATS2022} are not directly comparable to ours. It should be noted that the use of cyclic stepsizes is inherent in \ref{eq:BGD}. In fact,  \ref{eq:BGD} not only cycles between stepsizes, but also between blocks. One more difference is this: \ref{eq:BGD} cycles \textit{more frequently}, meaning that it leverages two stepsizes in a single iteration with a cost comparable to one GD step, while GD with cyclic stepsizes needs $M$ iterations to make full use of $M$ stepsizes.

% In the least-squares context, \citet{Oymak-SPL2021,Baptiste-AISTATS2022} analyzed GD with cyclic stepsizes. They assume the spectrum of $\bA^\top\bA$ is \textit{clustered} into two or more disjoint intervals; this is different from \cref{assumption:BWO}, so the results of \citet{Oymak-SPL2021,Baptiste-AISTATS2022} are not directly comparable to ours. It should be noted that the use of cyclic stepsizes is inherent in \ref{eq:BGD}. In fact,  \ref{eq:BGD} not only cycles between stepsizes, but also between blocks. One more difference is this: \ref{eq:BGD} cycles \textit{more frequently}, meaning that it leverages two stepsizes in a single iteration with a cost comparable to one GD step, while GD with cyclic stepsizes needs $M$ iterations to make full use of $M$ stepsizes.

\subsection{BGD Versus Lower Bounds} 
As shown in Assumption 2.1.4 \& Theorem 2.1.13 of \citet{Nesterov-2018}, for constants $L>\mu>0$, there is a $\mu$-strongly convex and $L$-smooth function $F(\bx)$, globally minimized at $\bx^*$, so that the iterates $\bx^t$ of any first-order method created from linear combinations of any initialization $\bx^0$ and previous gradient evaluations $\nabla F(\bx^0),\dots,\nabla F(\bx^{t-1})$ satisfy the lower bound
    \begin{equation}\label{eq:LB}
        \|\bx^t - \bx^* \|_2 \geq \Big( \frac{\sqrt{\kappa}-1}{\sqrt{\kappa}+1} \Big)^t \cdot  \|\bx^0 - \bx^* \|_2, \quad \kappa:= L/\mu.
    \end{equation}
    See also Theorem 3.1 of \citet{Sun-MP2021}, which proves a similar lower bound for the multi-block version of BGD with block sizes $1$ and canonical stepsizes $1/L$ applied to least-squares. While \ref{eq:HB} asymptotically attains lower bound \cref{eq:LB}, our \cref{theorem:informal} suggests that, under \cref{assumption:BWO}, the iterates $\bx^t$ of \ref{eq:BGD} with optimal stepsizes would satisfy
    \begin{equation*}%\label{eq:LB}
        \|\bx^t - \bx^* \|_2 \leq \Big( \frac{\sqrt{\kappa}-1}{\sqrt{\kappa}+1} \Big)^{2t} \cdot  \|\bx^0 - \bx^* \|_2
    \end{equation*}
    asymptotically (as $t\to \infty$). Is this a contradiction or does that mean \ref{eq:BGD} breaches lower bound \cref{eq:LB}? The answer is no and the reason is two-fold. First, \ref{eq:BGD} is not a first-order method in the precise sense of Assumption 2.1.4 of \citet{Nesterov-2018}, e.g., its new iterates do not arise as linear combinations of previous ones or their gradients. Second, the objective function we optimize is least-squares with \cref{assumption:BWO}, which does not necessarily belong to the class of worst-case functions in Theorem 2.1.13 of \citet{Nesterov-2018} or in Theorem 3.1 of \citet{Sun-MP2021}. One more subtlety is that the lower bound of \citet{Sun-MP2021} is derived for \textit{canonical}, rather than optimal, stepsizes.
    
    % First, the objective function we optimize is least-squares with \cref{assumption:BWO}, which does not necessarily belong to the class of worst-case functions in Theorem 2.1.13 of \citet{Nesterov-2018}. 

\subsection{BGD Versus Block Exact Minimization} Under \cref{assumption:BWO}, one verifies  \ref{eq:BGD} with stepsizes $\gamma_1=1,\gamma_2=1$ is equivalent\footnote{Under \cref{assumption:BWO}, the Lipschitz smoothness constant $L_1$ (or $L_2$) of \cref{eq:LS2} in variable $\bx_1$ (or $\bx_2$) is $1$, so stepsizes $(\gamma_1,\gamma_2)=(1,1)$ also correspond to the commonly used choice $(1/L_1,1/L_2)$. } to \textit{block exact minimization} (BEM), 
\begin{equation}\label{eq:BEM}
    \begin{split}
        \bx_1^+ \in \argmin_{\bx_1\in\bbR^{n_1} } F([\bx_1; \bx_2]), \\ \bx_2^+ \in \argmin_{\bx_2\in\bbR^{n_2} } F([\bx_1^+; \bx_2]),
    \end{split}  \tag{BEM}
\end{equation}
where $F$ is the least-squares objective in \cref{eq:LS2}, so \ref{eq:BGD} under \cref{assumption:BWO} generalizes \ref{eq:BEM} by allowing for different stepsizes. While \ref{eq:BEM} ensures the largest possible decrease of $F$ when updating either $\bx_1^+$ or $\bx_2^+$, we have seen the optimal stepsizes of \ref{eq:BGD} are not as simple as $\gamma_1=1,\gamma_2=1$; \ref{eq:BGD} can do better than \ref{eq:BEM} under \cref{assumption:BWO}.

\subsection{BGD and Generalized Alternating Projection}\label{subsection:review-GAP}
As shown in \cref{lemma:BGD-AP}, \ref{eq:BGD} is highly related to alternating projection between two subspaces, to which our results, say \cref{corollary:main-result}, apply. Hence, in this subsection, we review existing theoretical results on alternating projection and its variants, and then later in \cref{section:close-gap} we show how our \cref{corollary:main-result} applies and improves prior work.

% As will be shown in \cref{section:assumption}, \cref{assumption:BWO} makes a natural connection between \ref{eq:BGD} and \textit{generalized alternating projection} (GAP) \citep{Falt-CDC2017}. While the optimal convergence rate of GAP was claimed to be proved using the stepsizes of \citet{Falt-CDC2017}, we show this claim is inaccurate as applying our results on \ref{eq:BGD} yields better stepsizes and faster convergence (see \cref{table:compare-methods} and \cref{section:close-gap}).

% In this section, we review the classic problem of projecting a point onto the intersection of two linear subspaces as well as the convergence rates of the associated methods proposed in the literature to solve this problem.

Given two linear subspaces $\cH_1,\cH_2$ of $\bbR^{m}$ of dimension $m-n_1,m-n_2$ respectively, we are interested in projecting some point $\bz^{(0)}\in \bbR^{m}$ onto the intersection $\cH_1\cap \cH_2$. To avoid trivial cases, we assume $\cH_1$ does not contain $\cH_2$ (and vice versa), and none of the two subspaces is $\{0\}$.

We consider iterative algorithms of the form 
\begin{equation}\label{eq:T-generic}
	\bz^{(t+1)} = \bT \bz^{(t)},
\end{equation}
where $\bT\in \bbR^{m\times m}$ is some matrix that transforms the current iterate $\bz^{(t)}$ into the next, $\bz^{(t+1)}$. The choice of $\bT$ typically depends on the two subspaces $\cH_1,\cH_2$. We next review several algorithms that differ in how $\bT$ is chosen and that converge at different rates (see \cref{table:compare-methods} for a summary).

\subsubsection{Project Versus Reflect }\label{subsubsection:project-reflect}
Let us first review two classic methods to project a point onto $\cH_1\cap \cH_2$: alternating projection and Douglas-Rachford.

\myparagraph{Alternating Projection} Let $\bP_{\cH_1},\bP_{\cH_2}$ be matrices representing orthogonal projections onto $\cH_1,\cH_2$ respectively. The \textit{alternating projection} algorithm sets $\bT$ to be
\begin{equation}\label{eq:AP}
	\bT:=\bP_{\cH_2} \bP_{\cH_1}. \tag{AP}
\end{equation}
It follows immediately that for every $t$ we have
\begin{equation*}
	\| \bz^{(t)} - \bP_{\cH_1\cap \cH_2} \bz^{(0)} \|_2 \leq  \| \bT^k - \bP_{\cH_1\cap \cH_2} \|_2 \cdot \| \bz^{(0)} \|_2.
\end{equation*}
As \citet{Kayalar-MCSS1988} reviewed, Von Neumann established that \ref{eq:AP} converges to $\bP_{\cH_1\cap \cH_2} \bz^{(0)}$ in the early 1930s and his result was not in print until almost two decades after; see, e.g., Theorem 13.7 of \citet{Von-1951}. Around the same time, \citet{Aronszajn-1950} showed $\bT$ of \ref{eq:AP} satisfies
\begin{equation}\label{eq:rate-AP}
	\| \bT^k - \bP_{\cH_1\cap \cH_2} \|_2\leq \cos^{2k-1}\theta_1,
\end{equation}
where $\theta_1\in[0,\pi/2]$ is defined to be the minimum nonzero \textit{principal angle} between $\cH_1$ and $\cH_2$, also widely known as the Friedrichs angle \citep{Friedrichs-1937}.\footnote{ Alternating projection is also related to a popular line of research called \textit{continual learning} \citep{Peng-ICML2023,Elenter-arXiv2023,Cai-arXiv2024}. There, convergence results similar to \cref{eq:rate-AP} can be found; see, e.g., Theorem 8 of \citet{Evron-COLT2022}. Our stepsizes and theory can be applied there to improve the rate $\cos^2(\theta_1)$ in Theorem 8 of \citet{Evron-COLT2022} to $\frac{\sin (\theta_{r}) - \sin (\theta_{1 }) }{\sin (\theta_{r}) + \sin (\theta_{1 }) }$. Here, $\theta_r$ and $\theta_1$ denote the largest and nonzero smallest principal angles between $\cH_1$ and $\cH_2$, respectively. See also \cref{table:compare-methods}. } While \citet{Deutsch-1984,Franchetti-1986} announced bounds that are looser than \cref{eq:rate-AP}, Theorem 2 of \citet{Kayalar-MCSS1988} proved that \cref{eq:rate-AP} always holds with equality, which implies this linear rate $\cos^{2}\theta_1$ is actually optimal for \ref{eq:AP}.

\begin{table*}[t]
	\centering
	\caption{Comparison of GAP++ (\cref{section:close-gap}) with existing rates. A smaller linear rate means better, and $\theta_r$ and $\theta_1$ denote the largest and nonzero smallest principal angles between $\cH_1$ and $\cH_2$, respectively. We find novel stepsizes (GAP++) that guarantee faster convergence. (The stepsizes $\gamma_1^*,\gamma_2^*$ of GAP++ are shown in \cref{remark:opt-stepsizes-AP} of the appendix.)}
	\vspace{-0.4em}
	\footnotesize
	\renewcommand{\tabcolsep}{3.5pt}
	\renewcommand\arraystretch{1.4}
	\begin{tabular}{l|cccc}
		\hline
		Method  & Linear Rate Factor &  Stepsizes &  Asymptotic? & Reference\\ 
		\hline
		\ref{eq:AP}  & $\cos^2 (\theta_{1 } ) $  & N.A. &  No & Theorem 2 \citep{Kayalar-MCSS1988}\\
		\ref{eq:DR} &  $\cos (\theta_{ 1 }) $ & N.A.  &   No & Theorem 4.3 \citep{Bauschke-JAT2014} \\
		\hline
		\ref{eq:RAP}  &   $\frac{1 - \sin^2 (\theta_{1 } ) }{1 + \sin^2 (\theta_{1 } ) }$ & $\gamma=\frac{2}{1+\sin^2(\theta_1)}$ & Yes & Theorem 3.6 \citep{Bauschke-NA2016}   \\
		\ref{eq:PRAP} & $\frac{\sin^2 (\theta_{r }) - \sin^2 (\theta_{1 }) }{\sin^2 (\theta_{r }) + \sin^2 (\theta_{1 }) }$ & $\gamma_1=\frac{ 2 }{\sin^2 (\theta_{r }) + \sin^2 (\theta_{1 }) }$ & Yes  & Theorem 3.7 \citep{Bauschke-NA2016} \\
		\ref{eq:GDR} &  $\cos (\theta_{ 1 }) $  & $\gamma=1$ &  No & Theorem 3.10 \citep{Bauschke-NA2016}  \\ 
		\ref{eq:GAP}   & $\frac{1 - \sin (\theta_{1 }) }{1 + \sin (\theta_{1 }) }$ & $\gamma=1, \gamma_1=\gamma_2=\frac{2}{1+\sin(\theta_1)}$ & Yes & Theorem 3, Remark 3 \citep{Falt-CDC2017} \\
		\hline
		GAP++ & $\frac{\sin (\theta_{r}) - \sin (\theta_{1 }) }{\sin (\theta_{r}) + \sin (\theta_{1 }) }$ & $\gamma=1,\gamma_1=\gamma_1^*, \gamma_2=\gamma_2^*$ &   Yes  & This Paper, \cref{corollary:GAP} \\
		\hline
	\end{tabular}
	\label{table:compare-methods}
\end{table*}

\myparagraph{Douglas-Rachford} Another classic iterative scheme proposed (implicitly) by \citet{Douglas-1956} is this:
\begin{equation}\label{eq:DR}
    \begin{split}
        \bT :=&\ \bP_{\cH_2} ( 2\bP_{\cH_1} - \bI) + \bI - \bP_{\cH_1} \\
            =&\ \frac{1}{2}\bI + \frac{1}{2}( 2\bP_{\cH_2} - \bI)( 2\bP_{\cH_1} - \bI).
    \end{split}    \tag{DR}
\end{equation}
In words, \ref{eq:DR} consists of applying in cascade two reflections, $2\bP_{\cH_1} - \bI$ and then $2\bP_{\cH_2} - \bI$, to the current iterate $\bz^{(t)}$ and then taking an average. \citet{Hesse-TSP2014} showed in their Theorem IV.6 that \ref{eq:DR} converges linearly without providing any explicit rate, while \citet{Bauschke-JAT2014} later made the rate precise: Their Theorem 4.3 proves $\bT$ of \ref{eq:DR} satisfies
\begin{equation*}
	\|\bP_{\cH_1} \bT^k - \bP_{\cH_1\cap \cH_2} \|_2 = \|\bP_{\cH_2} \bT^k - \bP_{\cH_1\cap \cH_2} \|_2 = \cos^{k}\theta_1 .
\end{equation*}
We can then say that \ref{eq:DR}  converges linearly with rate $ \cos\theta_1$ and this is optimal for \ref{eq:DR}. Note though that since $\cos^2\theta_1\leq \cos\theta_1$, \ref{eq:DR} is provably slower than \ref{eq:AP}.

% It comes to our attention that the above methodologies have been extended to more general settings involving multiple projections onto (non-)convex sets. Nevertheless, we 

\subsubsection{Four Generalized Methods}\label{subsubsection:generalized-method}
With two positive numbers $\gamma_1,\gamma_2$ and the $m\times m$ identity matrix $\bI$, define the relaxed projections $\bP_{\cH_1}^{\gamma_1}$ and $\bP_{\cH_2}^{\gamma_2}$ as
\begin{equation}\label{eq:relaxed-P1P2}
	\begin{split}
		\bP_{\cH_1}^{\gamma_1}&:= (1-\gamma_1)\bI + \gamma_1 \bP_{\cH_1},\\
		\bP_{\cH_2}^{\gamma_2}&:= (1-\gamma_2)\bI + \gamma_1 \bP_{\cH_2}.
	\end{split}
\end{equation}
Note that $\bP_{\cH_1}^{1}$ represents the projection $\bP_{\cH_1}$ and $\bP_{\cH_1}^{2}$ the reflection, so $\bP_{\cH_1}^{\gamma_1}$ generalizes the operators in \ref{eq:AP} and \ref{eq:DR}. For clarity, we will always write $\bP_{\cH_1}^{\gamma_1}$ to mean definitions in \cref{eq:relaxed-P1P2} and it never means the matrix $\bP_{\cH_1}$ raised to the power of $\gamma_1$; if ever needed, we will write $(\bP_{\cH_1})^{\gamma_1}$ for the latter.

% Since $\bP_{\cH_1}^{\gamma_1}$ might also imply  $\bP_{\cH_1}$ raised to the power of $\gamma_1$, what we mean by $\bP_{\cH_1}^{\gamma_1}$ would be clear from the context.

With the relaxed projections $\bP_{\cH_1}^{\gamma_1}$ and $\bP_{\cH_2}^{\gamma_2}$, we can now proceed to review methods that generalize \ref{eq:AP} and \ref{eq:DR}.

\myparagraph{Relaxed Alternating Projection} For some $\gamma\in(0,1]$, the  method we review now has the iterates $ \bz^{(t+1)} = \bT(\gamma) \bz^{(t)}$, where $\bT(\gamma)$ is a matrix function of $\gamma$, defined as
\begin{equation}\label{eq:RAP}
	\bT(\gamma) := (1-\gamma)\bI + \gamma \bP_{\cH_2}\bP_{\cH_1}.   \tag{RAP}
\end{equation}
Clearly, $\bT(1)$ corresponds to \ref{eq:AP}. How to choose $\gamma$ to obtain a faster rate than \ref{eq:AP}? The optimal choice of $\gamma$ would be $1$ if $\bP_{\cH_2}\bP_{\cH_1}$ were symmetric, while it is precisely the asymmetry of $\bP_{\cH_2}\bP_{\cH_1}$ that brings difficulty. \citet{Bauschke-NA2016} showed in their Theorem 3.6 that the choice $\gamma=\frac{2}{1+\sin^2(\theta_1)}$ is optimal, with which \ref{eq:RAP} converges asymptotically at a linear rate $\frac{1-\sin^2(\theta_1)}{1+\sin^2(\theta_1)}$ faster than \ref{eq:AP}.

%leading to a faster linear rate $\frac{1-\sin^2(\theta_1)}{1+\sin^2(\theta_1)}$ for \ref{eq:RAP}.

\myparagraph{Partial Relaxed Alternating Projection}
The method we present now relaxes \ref{eq:AP} \textit{partially} by setting
\begin{equation}\label{eq:PRAP}
	\bT(\gamma_1) := \bP_{\cH_2} \bP_{\cH_1}^{\gamma_1}. \tag{PRAP}
\end{equation}
Theorem 3.7 of \citet{Bauschke-NA2016} proves the optimal choice of $\gamma_1$ for \ref{eq:PRAP} is $\frac{ 2 }{\sin^2 (\theta_{r }) + \sin^2 (\theta_{1 }) }$, giving the asymptotic linear rate $\frac{\sin^2 (\theta_{r }) - \sin^2 (\theta_{1 }) }{\sin^2 (\theta_{r }) + \sin^2 (\theta_{1 }) }$. Moreover, one verifies this rate is faster than or the same as those of \ref{eq:AP}, \ref{eq:DR}, and \ref{eq:RAP}.

\myparagraph{Generalized Douglas-Rachford} In their Section 1.4, \citet{Demanet-MC2016} considered the following operator $T(\gamma)$ with iterate update $z^{k+1}= \bT(\gamma) z^k$:
\begin{equation}\label{eq:GDR}
	\bT(\gamma) =(1-\gamma)\bI + \gamma \bigg(\frac{1}{2}\bI + \frac{1}{2} \bP_{\cH_2}^{2} \bP_{\cH_1}^{2} \bigg).   \tag{GDR}
\end{equation}
\ref{eq:GDR} can be viewed as taking a convex combination of $\bI$ and the operator of \ref{eq:DR}, similarly to \ref{eq:RAP} and \ref{eq:AP}. However, choosing $\gamma\in[0,1]$ for \ref{eq:GDR} does not improve the rate: Theorem 3.10 of \citet{Bauschke-NA2016} shows  $\gamma=1$ is optimal for \ref{eq:GDR}, giving the same linear rate $\cos(\theta_1)$ as \ref{eq:DR}.

\myparagraph{Generalized Alternating Projection}
The generalized alternating projection algorithm is with the operator
\begin{equation}\label{eq:GAP}
	\bT(\gamma, \gamma_1,\gamma_2) =  (1-\gamma) \bI + \gamma \bP_{\cH_2}^{\gamma_2} \bP_{\cH_1}^{\gamma_1}. \tag{GAP}
\end{equation}
where $\gamma$ is assumed to lie in $(0,1]$. If $\gamma$ and $\gamma_2$ were taken to be $1$, then \ref{eq:GAP} recovers \ref{eq:PRAP}. \citet{Falt-CDC2017} proposed the stepsizes $\gamma=1, \gamma_1=\gamma_2=\frac{2}{1+\sin(\theta_1)}$, with which \ref{eq:GAP} converges asymptotically at linear rate $\frac{1 - \sin (\theta_{1 }) }{1 + \sin (\theta_{1 }) }$. Moreover, \citet{Falt-CDC2017} shows this choice of stepsizes is optimal as long as the largest principal angle $\theta_r$ between the two subspace $\cH_1,\cH_2$ is equal to $\pi/2$. Indeed, if $\theta_r$ were equal to $\pi/2$, then the rate $\frac{1 - \sin (\theta_{1 }) }{1 + \sin (\theta_{1 }) }$ of \citet{Falt-CDC2017} is the smallest among prior methods.

% Moreover, the rate $\frac{1 - \sin (\theta_{1 }) }{1 + \sin (\theta_{1 }) }$ is faster than if not the same as those of XXX.

\subsubsection{The Gap}\label{subsubsection:gap}
In light of the review in \cref{subsubsection:generalized-method}, a major gap to be bridged is as follows. The assumption $\theta_r=\pi/2$ of \citet{Falt-CDC2017} on the largest principal angle $\theta_r$ is violated with probability $1$ by any two subspaces of fixed dimensions randomly sampled from their respective Grassmannian manifolds; this means their stepsizes are sub-optimal with probability $1$. This sub-optimality manifests itself in comparison to \ref{eq:PRAP}: The stepsizes of \citet{Falt-CDC2017} do not necessarily yield a smaller rate than \ref{eq:PRAP} of \citet{Bauschke-NA2016}. More specifically, the rate $\frac{\sin^2 (\theta_{r }) - \sin^2 (\theta_{1 }) }{\sin^2 (\theta_{r }) + \sin^2 (\theta_{1 }) }$ of \ref{eq:PRAP} could be larger or smaller than the rate $\frac{1 - \sin (\theta_{1 }) }{1 + \sin (\theta_{1 }) }$ of \citet{Falt-CDC2017}, depending on problem configurations (e.g., the precise values of $\theta_r$).

\section{GAP++: Bridge the GAP}\label{section:close-gap}
A major contribution of this section is closing the gap discussed in \cref{subsubsection:gap}. This is achieved by invoking our theoretical results developed so far. Specifically, with the stepsizes that we propose, the \cref{eq:GAP} operator converges asymptotically at rate $\frac{\sin (\theta_{r}) - \sin (\theta_{1 }) }{\sin (\theta_{r}) + \sin (\theta_{1 }) }$. This rate is alway better than if not the same as the rates of \citet{Bauschke-NA2016} and  \citet{Falt-CDC2017}, for all possible choices of $\theta_1$ and $\theta_r$. On the conceptual level, our result can be viewed as a generalization of \citet{Falt-CDC2017} which dispenses with their assumption $\theta_r=\pi/2$. Therefore we call the proposed stepsizes GAP++.

In \cref{subsection:preliminaries} we build preliminary notations. In \cref{subsection:our-results-AP} we provide stepsizes that achieve the fastest rate in \cref{table:compare-methods}.

\subsection{Preliminaries}\label{subsection:preliminaries}
\myparagraph{Matrix Representations} For $j=1,2$, let $A_j\in \bbR^{m\times n_j}$ be an orthonormal basis matrix of the orthogonal complement $\cH_j^\perp$ of $\cH_j$. Then we can write $\bP_{\cH_j}:=\bI- A_j A_j^\top$ and the relaxed projections $\bP_{\cH_1}^{\gamma_2}, \bP_{\cH_1}^{\gamma_2}$  in \cref{eq:relaxed-P1P2} as 
\begin{equation}\label{eq:project-matrix-A}
	\bP_{\cH_j}^{\gamma_j}= \bI - \gamma_j A_j A_j^\top, \quad j=1, 2.
\end{equation}
Note that, here, we write the symbols $A_1,A_2$ without boldface to distinguish them from the notations $\bA_1,\bA_2$ of the main paper. However, the choices of the same letter $A, \bA$ should remind the reader that they will eventually be connected somehow.

The \ref{eq:GAP} operator can be expressed in terms of $A_1,A_2$:
\begin{equation*}
	\bT(\gamma, \gamma_1,\gamma_2)=  (1-\gamma) \bI + \gamma ( \bI - \gamma_2 A_2 A_2^\top) ( \bI - \gamma_1 A_1 A_1^\top).
\end{equation*}
We will set $\gamma=1$ and analyze $\bT(1, \gamma_1,\gamma_2)=( \bI - \gamma_2 A_2 A_2^\top) ( \bI - \gamma_1 A_1 A_1^\top)$, as doing so is sufficient to derive a faster rate than that of \citet{Falt-CDC2017}.

\myparagraph{Principal Angles} We recall the following definition of \textit{nonzero} principal angles.
\begin{definition}
	For  subspaces $\cH_1,\cH_2$ of $\bbR^m$, define $r:=\min\{ \dim(\cH_1) , \dim(\cH_2) \} - \dim (\cH_1\cap \cH_2)$ and $\cS_1:=\cH_1 \cap (\cH_1\cap \cH_2)^\perp$ and $\cZ_1:=\cH_2\cap (\cH_1\cap \cH_2)^\perp$. The \textit{nonzero principal angles} $\theta_1,\dots,\theta_r$ between $\cH_1$ and $\cH_2$ are then defined recursively for $i=1,\dots,r$ such that
	\begin{equation}
		(\bs_i,\bz_i)\in \argmax_{\substack{\bs\in \cS_{i},\bz\in \cZ_{i} \\ \| \bs\|_2=\| \bz\|_2=1}  } \bs^\top \bz, \quad  \theta_i = \arccos(\bs_i^\top \bz_i),
	\end{equation}
	where $\cS_{i}:= \cS_{i-1}\cap \Span(\bs_i)^\perp$, $\cZ_{i}:= \cZ_{i-1}\cap \Span(\bz_i)^\perp$.
\end{definition}
It is not hard to show the principal angles $\theta_1,\dots,\theta_r$ between $\cH_1$ and $\cH_2$ are indeed nonzero and satisfy $\theta_1\leq \cdots \leq \theta_r$. Moreover, these principal angles between $\cH_1,\cH_2$ correspond exactly to nonzero principal angles between the orthogonal complements $\cH_1^\perp$ and $\cH_2^\perp$; see Property 2.1 of \citet{Zhu-JNM2013} or Theorem 2.7 of \citet{Knyazev-SIAM-J-MAA2007}. More formally, we have the following statements:

\begin{lemma}\label{eq:singular-val=angle}
    Recall that $A_j\in \bbR^{m\times n_j}$ is an orthonormal basis matrix of the orthogonal complement $\cH_j^\perp$ of $\cH_j$ ($j=1,2$), and that $r:=\min\{ \dim(\cH_1) , \dim(\cH_2) \} - \dim (\cH_1\cap \cH_2)$. Define $\bC:=A_2^\top A_1$. Then we have $\rank(\bC)=r$. Denote by $\sigma_1\geq \cdots \geq \sigma_r $ the $r$ nonzero singular values of $\bC$. Then $\sigma_i=\cos(\theta_i)$ for every $i=1,\dots,r$.
% Let $\bB_1,\bB_2$ be  orthonormal basis matrices respectively for $\cH_1,\cH_2$. Let $\sigma_1\geq \cdots \geq \sigma_d $ be the $d$ singular values of $\bB_1^\top \bB_2$. Then $\sigma_i=\cos(\theta_i)$ for every $i=1,\dots,d$.
\end{lemma}

% singular values of 

\subsection{Asymptotic Convergence Rates}\label{subsection:our-results-AP}
Here we give stepsize choices $\gamma_1,\gamma_2$ and the associated asymptotic convergence rates for the iterates $\bz^+ = \bT(1, \gamma_1,\gamma_2) \bz$, where we recall $\bT(1, \gamma_1,\gamma_2)$ is defined as $\bT(1, \gamma_1,\gamma_2)=( \bI - \gamma_2 A_2 A_2^\top) ( \bI - \gamma_1 A_1 A_1^\top)$. 

Write $A:= [A_1\ A_2]\in\bbR^{m\times (n_1+n_2)}$ and assume $A$ has full column rank (necessarily, we assume $m\geq n_1+n_2=:n$). Let $\bQ$ be an orthonormal basis matrix for $\cH_1\cap \cH_2$. Then $\bQ^\top\bQ$ is the identity matrix, and we have $\bQ^\top A=0$. Furthermore, we can write $\bz$ as $\bz=A \bc_1 + \bQ\bc_2$ and then
\begin{equation*}
    \begin{split}
        &\ \bT(1, \gamma_1,\gamma_2)\bz= \bT(1, \gamma_1,\gamma_2) (A \bc_1 + \bQ\bc_2) \\ 
        =&\ ( \bI - \gamma_2 A_2 A_2^\top) ( \bI - \gamma_1 A_1 A_1^\top) ( A \bc_1 + \bQ\bc_2 )  \\ 
        =&\ ( \bI - \gamma_2 A_2 A_2^\top) ( \bI - \gamma_1 A_1 A_1^\top)  A \bc_1 + \bQ\bc_2  .
        % &=A A^\dagger  ( \bI - \gamma_2 A_2 A_2^\top) ( \bI - \gamma_1 A_1 A_1^\top)  A \bc_1 + \bQ\bc_2.
    \end{split}
\end{equation*}
Since the columns of $\bQ$ span $\cH_1\cap \cH_2$ and the goal is to reach a point in  $\cH_1\cap \cH_2$, we could assume $\bc_2=0$ without loss of generality, and analyze the update equation
\begin{equation*}
    \bd:=\bT(1, \gamma_1,\gamma_2) A \bc_1 = ( \bI - \gamma_2 A_2 A_2^\top) ( \bI - \gamma_1 A_1 A_1^\top)  A \bc_1.
\end{equation*}
Moreover, since $A$ has full column rank and $\bQ^\top \bd = \bQ^\top ( \bI - \gamma_2 A_2 A_2^\top) ( \bI - \gamma_1 A_1 A_1^\top)  A \bc_1=0$, the new iterate $\bd$ always lies in the range space of $A$, meaning that we can uniquely write $\bd$ as $\bd=\bA\bc^+_1$ for some $\bc^+_1$. Then the update equation that we need to analyze become
\begin{align}
    %\begin{split}
        &\ A \bc^+_1 = ( \bI - \gamma_2 A_2 A_2^\top) ( \bI - \gamma_1 A_1 A_1^\top)  A \bc_1 \label{eq:AP->BCD} \\
    \Leftrightarrow &\   \bc^+_1 = (A^\top A)^{-1} A^\top ( \bI - \gamma_2 A_2 A_2^\top) ( \bI - \gamma_1 A_1 A_1^\top)  A \bc_1. \nonumber
    %\end{split}
\end{align}
We need the following lemma to proceed.
% In the above, the last inequality holds as $AA^\dagger$ is the projection onto the range space of $A$, whereas $( \bI - \gamma_2 A_2 A_2^\top) ( \bI - \gamma_1 A_1 A_1^\top)  A \bc_1$ already lies in this range space (here $A^\dagger$ denotes the pseudoinverse of $A$). Then we raise $\bT(1, \gamma_1,\gamma_2)$ to the power of $t$, and obtain:
% \begin{equation}\label{eq:rewrite-T1}
%     \begin{split}
%         \bT(1, \gamma_1,\gamma_2)^t \bz&= \bT(1, \gamma_1,\gamma_2)^t (A \bc_1 + \bQ\bc_2) \\
%         &= \bT(1, \gamma_1,\gamma_2)^t A \bc_1 + \bQ\bc_2 \\
%         &=\Big( ( \bI - \gamma_2 A_2 A_2^\top) ( \bI - \gamma_1 A_1 A_1^\top) 
%  \Big)^t A \bc_1 + \bQ\bc_2 \\
%         &=\Big( A A^\dagger( \bI - \gamma_2 A_2 A_2^\top) ( \bI - \gamma_1 A_1 A_1^\top) 
%  \Big)^t A \bc_1 + \bQ\bc_2 \\
%     &=A \Big(  A^\dagger( \bI - \gamma_2 A_2 A_2^\top) ( \bI - \gamma_1 A_1 A_1^\top)A 
%  \Big)^t \bc_1 + \bQ\bc_2 \\
%     \end{split}
% \end{equation}
% Since the columns of $\bQ$ span $\cH_1\cap \cH_2$ and the goal is to reach a point in  $\cH_1\cap \cH_2$, we need only to analyze the convergence rate of $A \Big(  A^\dagger( \bI - \gamma_2 A_2 A_2^\top) ( \bI - \gamma_1 A_1 A_1^\top)A \Big)^t \bc_1$ to zero. To do so, let us make the following observation:
\begin{lemma}\label{lemma:BCD-AP-matrix}
If $A\in\bbR^{m\times (n_1+n_2)}$ has full column rank, then the matrix $(A^\top A)^{-1} A^\top ( \bI - \gamma_2 A_2 A_2^\top) ( \bI - \gamma_1 A_1 A_1^\top)A $ is equal to $M(\gamma_1,\gamma_2)$, where $M(\gamma_1,\gamma_2)$ is defined to be
\begin{equation*}
    \left(\bI - \begin{bmatrix}
            0 & 0 \\ 
            \gamma_2 A_2^\top A_1 & \gamma_2 \bI_2 
        \end{bmatrix} \Bigg) \Bigg(\bI - \begin{bmatrix}
            \gamma_1 \bI_1 & \gamma_1 A_1^\top A_2 \\ 
            0 & 0
        \end{bmatrix} \right). % =:M(\gamma_1,\gamma_2).
\end{equation*}
% where we recall $\bM(\gamma_1,\gamma_2)$ is defined in \cref{lemma:error-decrease}.

\end{lemma}
With \cref{eq:AP->BCD} and \cref{lemma:BCD-AP-matrix}, we conclude that we need to analyze the convergence rate of the iterate update
\begin{equation}\label{eq:AP-iterates}
    \bc_1^+ = M(\gamma_1,\gamma_2) \bc_1,
\end{equation}
where $M(\gamma_1,\gamma_2)$ is defined in \cref{lemma:BCD-AP-matrix}. Note that the definition of $M(\gamma_1,\gamma_2)$ here corresponds precisely to the definition of $\bM(\gamma_1,\gamma_2)$ in \cref{lemma:error-decrease} if $\bA_1,\bA_2$ satisfy the block-wise orthogonality assumption. Moreover, $A_1$ and $A_2$ already satisfy this assumption, therefore we can now make the following conclusion by invoking \cref{corollary:main-result}.
\begin{corollary}\label{corollary:GAP}
    Assume $A$ has full column rank and the $n_2\times n_1$ matrix $C:=A_2^\top A_1$ is nonzero. The following hold.
    \begin{itemize}[parsep=-1pt,topsep=0pt]
        \item If $C$ has full rank (i.e., $\rank(C)=\min\{n_1,n_2\}$), then
        \begin{equation*}
        \min_{\gamma_1>0,\gamma_2>0} \rho\big( M(\gamma_1,\gamma_2) \big) = \frac{\sin(\theta_r)-\sin(\theta_1)}{\sin(\theta_r)+\sin(\theta_1)}.
    \end{equation*}
        \item If $C$ is rank-deficient, then 
        \begin{equation*}
        \min_{\gamma_1>0,\gamma_2>0} \rho\big( M(\gamma_1,\gamma_2) \big) = \frac{1-\sin(\theta_1)}{1+\sin(\theta_1)}.
    \end{equation*}
    \end{itemize}    
    % Define $r:=\rank(C)$. Then we have
    % \begin{equation*}
    %     \min_{\gamma_1>0,\gamma_2>0} \rho\big( M(\gamma_1,\gamma_2) \big) = \begin{cases}
    %        \frac{\sin(\theta_r)-\sin(\theta_1)}{\sin(\theta_r)+\sin(\theta_1)}  & r = \min\{n_1,n_2\}; \\
    %        \frac{1-\sin(\theta_1}{1+\sin(\theta_1)} & r < \min\{n_1,n_2\}.
    %     \end{cases}
    % \end{equation*}
\end{corollary}
Note that for two \textit{generic} subspaces $\cH_1, \cH_2$, with probability $1$ the following hold: The matricx $A=[A_1,\ A_2]$ has full column rank if $n_1+n_2\leq m$; and $A_2^\top A_1$ has full rank.

\begin{figure*}
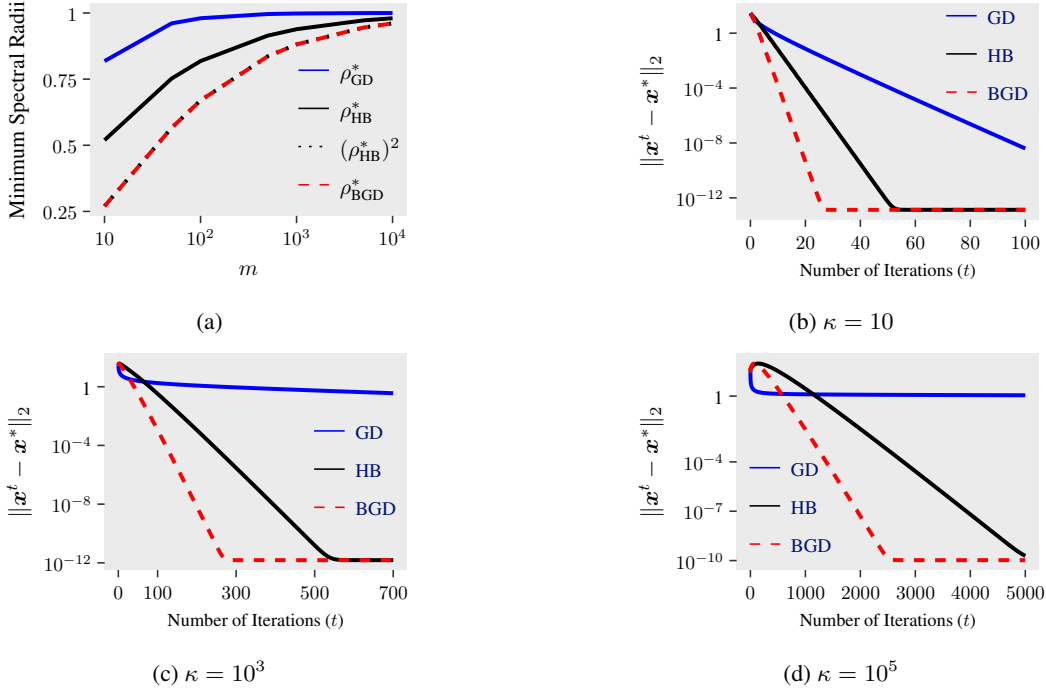

	\centering
	\begin{subfigure}{.49\textwidth}
		\centering
		% Created by tikzDevice version 0.12.3.1 on 2023-07-22 23:43:06
% !TEX encoding = UTF-8 Unicode
\begin{tikzpicture}[x=1pt,y=1pt]
\definecolor{fillColor}{RGB}{255,255,255}
\path[use as bounding box,fill=fillColor,fill opacity=0.00] (0,0) rectangle (162.61,115.63);
\begin{scope}
\path[clip] (  0.00,  0.00) rectangle (162.61,115.63);
\definecolor{drawColor}{RGB}{255,255,255}
\definecolor{fillColor}{RGB}{255,255,255}

\path[draw=drawColor,line width= 0.6pt,line join=round,line cap=round,fill=fillColor] ( -0.00,  0.00) rectangle (162.61,115.63);
\end{scope}
\begin{scope}
\path[clip] ( 37.26, 29.71) rectangle (157.11,110.13);
\definecolor{fillColor}{gray}{0.92}

\path[fill=fillColor] ( 37.26, 29.71) rectangle (157.11,110.13);
\definecolor{drawColor}{RGB}{0,0,255}

\path[draw=drawColor,line width= 1.5pt,line join=round] ( 42.71, 88.29) --
	( 68.09,102.57) --
	( 79.03,104.51) --
	(104.41,106.10) --
	(115.34,106.30) --
	(140.73,106.46) --
	(151.66,106.48);
\definecolor{drawColor}{RGB}{0,0,0}

\path[draw=drawColor,line width= 1.5pt,line join=round] ( 42.71, 58.37) --
	( 68.09, 81.68) --
	( 79.03, 88.29) --
	(104.41, 97.92) --
	(115.34,100.36) --
	(140.73,103.70) --
	(151.66,104.51);

\path[draw=drawColor,line width= 1.5pt,dash pattern=on 1pt off 3pt ,line join=round] ( 42.71, 33.37) --
	( 68.09, 63.01) --
	( 79.03, 73.39) --
	(104.41, 90.08) --
	(115.34, 94.59) --
	(140.73,100.99) --
	(151.66,102.57);
\definecolor{drawColor}{RGB}{255,0,0}

\path[draw=drawColor,draw opacity=0.50,line width= 1.5pt,dash pattern=on 4pt off 4pt ,line join=round] ( 42.71, 33.37) --
	( 68.09, 63.01) --
	( 79.03, 73.39) --
	(104.41, 90.08) --
	(115.34, 94.59) --
	(140.73,100.99) --
	(151.66,102.57);
\end{scope}
\begin{scope}
\path[clip] (  0.00,  0.00) rectangle (162.61,115.63);
\definecolor{drawColor}{RGB}{0,0,0}

\node[text=drawColor,anchor=base east,inner sep=0pt, outer sep=0pt, scale=  0.73] at ( 32.31, 28.35) {$0.25$};

\node[text=drawColor,anchor=base east,inner sep=0pt, outer sep=0pt, scale=  0.73] at ( 32.31, 53.39) {$0.5$};

\node[text=drawColor,anchor=base east,inner sep=0pt, outer sep=0pt, scale=  0.73] at ( 32.31, 78.43) {$0.75$};

\node[text=drawColor,anchor=base east,inner sep=0pt, outer sep=0pt, scale=  0.73] at ( 32.31,103.47) {$1$};
\end{scope}
\begin{scope}
\path[clip] (  0.00,  0.00) rectangle (162.61,115.63);
\definecolor{drawColor}{gray}{0.20}

\path[draw=drawColor,line width= 0.6pt,line join=round] ( 34.51, 31.38) --
	( 37.26, 31.38);

\path[draw=drawColor,line width= 0.6pt,line join=round] ( 34.51, 56.42) --
	( 37.26, 56.42);

\path[draw=drawColor,line width= 0.6pt,line join=round] ( 34.51, 81.46) --
	( 37.26, 81.46);

\path[draw=drawColor,line width= 0.6pt,line join=round] ( 34.51,106.50) --
	( 37.26,106.50);
\end{scope}
\begin{scope}
\path[clip] (  0.00,  0.00) rectangle (162.61,115.63);
\definecolor{drawColor}{gray}{0.20}

\path[draw=drawColor,line width= 0.6pt,line join=round] ( 42.71, 26.96) --
	( 42.71, 29.71);

\path[draw=drawColor,line width= 0.6pt,line join=round] ( 79.03, 26.96) --
	( 79.03, 29.71);

\path[draw=drawColor,line width= 0.6pt,line join=round] (115.34, 26.96) --
	(115.34, 29.71);

\path[draw=drawColor,line width= 0.6pt,line join=round] (151.66, 26.96) --
	(151.66, 29.71);
\end{scope}
\begin{scope}
\path[clip] (  0.00,  0.00) rectangle (162.61,115.63);
\definecolor{drawColor}{RGB}{0,0,0}

\node[text=drawColor,anchor=base,inner sep=0pt, outer sep=0pt, scale=  0.73] at ( 42.71, 18.70) {$10$};

\node[text=drawColor,anchor=base,inner sep=0pt, outer sep=0pt, scale=  0.73] at ( 79.03, 18.70) {$10^2$};

\node[text=drawColor,anchor=base,inner sep=0pt, outer sep=0pt, scale=  0.73] at (115.34, 18.70) {$10^3$};

\node[text=drawColor,anchor=base,inner sep=0pt, outer sep=0pt, scale=  0.73] at (151.66, 18.70) {$10^4$};
\end{scope}
\begin{scope}
\path[clip] (  0.00,  0.00) rectangle (162.61,115.63);
\definecolor{drawColor}{RGB}{0,0,0}

\node[text=drawColor,anchor=base,inner sep=0pt, outer sep=0pt, scale=  0.83] at ( 97.18,  7.42) {$m$};
\end{scope}
\begin{scope}
\path[clip] (  0.00,  0.00) rectangle (162.61,115.63);
\definecolor{drawColor}{RGB}{0,0,0}

\node[text=drawColor,rotate= 90.00,anchor=base,inner sep=0pt, outer sep=0pt, scale=  0.83] at ( 12.32, 69.92) {Minimum Spectral Radii};
\end{scope}
\begin{scope}
\path[clip] (  0.00,  0.00) rectangle (162.61,115.63);
\definecolor{drawColor}{RGB}{0,0,255}

\path[draw=drawColor,line width= 0.6pt,line join=round] (116.05, 84.83) -- (127.61, 84.83);
\end{scope}
\begin{scope}
\path[clip] (  0.00,  0.00) rectangle (162.61,115.63);
\definecolor{drawColor}{RGB}{0,0,0}

\path[draw=drawColor,line width= 0.6pt,line join=round] (116.05, 70.38) -- (127.61, 70.38);
\end{scope}
\begin{scope}
\path[clip] (  0.00,  0.00) rectangle (162.61,115.63);
\definecolor{drawColor}{RGB}{0,0,0}

\path[draw=drawColor,line width= 0.6pt,dash pattern=on 1pt off 3pt ,line join=round] (116.05, 55.93) -- (127.61, 55.93);
\end{scope}
\begin{scope}
\path[clip] (  0.00,  0.00) rectangle (162.61,115.63);
\definecolor{drawColor}{RGB}{255,0,0}

\path[draw=drawColor,draw opacity=0.50,line width= 0.6pt,dash pattern=on 4pt off 4pt ,line join=round] (116.05, 41.47) -- (127.61, 41.47);
\end{scope}
\begin{scope}
\path[clip] (  0.00,  0.00) rectangle (162.61,115.63);
\definecolor{drawColor}{RGB}{0,0,0}

\node[text=drawColor,anchor=base west,inner sep=0pt, outer sep=0pt, scale=  0.83] at (131.66, 81.42) {$\rho_{\textnormal{GD}}^*$};
\end{scope}
\begin{scope}
\path[clip] (  0.00,  0.00) rectangle (162.61,115.63);
\definecolor{drawColor}{RGB}{0,0,0}

\node[text=drawColor,anchor=base west,inner sep=0pt, outer sep=0pt, scale=  0.83] at (131.66, 66.97) {$\rho_{\textnormal{HB}}^*$};
\end{scope}
\begin{scope}
\path[clip] (  0.00,  0.00) rectangle (162.61,115.63);
\definecolor{drawColor}{RGB}{0,0,0}

\node[text=drawColor,anchor=base west,inner sep=0pt, outer sep=0pt, scale=  0.83] at (131.66, 52.52) {$(\rho_{\textnormal{HB}}^*)^2$};
\end{scope}
\begin{scope}
\path[clip] (  0.00,  0.00) rectangle (162.61,115.63);
\definecolor{drawColor}{RGB}{0,0,0}

\node[text=drawColor,anchor=base west,inner sep=0pt, outer sep=0pt, scale=  0.83] at (131.66, 38.06) {$\rho_{\textnormal{BGD}}^*$};
\end{scope}
\end{tikzpicture}
		\caption{ } \label{figa:sr}
	\end{subfigure}%
	\begin{subfigure}{.49\textwidth}
		\centering
		\input{./figures/experiment_convergence_10.tex}
		\caption{$\kappa=10$} \label{fig:convergence10}
	\end{subfigure}%
	
	\begin{subfigure}{.49\textwidth}
		\centering
	\input{./figures/experiment_convergence_1000.tex}
		\caption{$\kappa=10^3$} \label{fig:convergence1000}
	\end{subfigure}%
	\begin{subfigure}{.49\textwidth}
		\centering	\input{./figures/experiment_convergence_100000.tex}
		\caption{$\kappa=10^5$} \label{fig:convergence100000}
	\end{subfigure}
	\caption{Under \cref{assumption:BWO}, \cref{figa:sr} shows the numerical values of the minimum spectral radii, $\rho_{\textnormal{GD}}^*, \rho_{\textnormal{HB}}^*$, and $\rho_{\textnormal{BGD}}^*$, and  \cref{fig:convergence10,fig:convergence1000,fig:convergence100000} shows the errors of the three methods at every iteration $t$. \label{fig:sr-convergence}  }
\end{figure*}
\section{Numerical Experiments}\label{section:experiments}
In this section, we perform a simple experiment in order to validate our theory. The full MATLAB code for the experiment can be found in \cref{subsection:code-w.A1}. 

\myparagraph{Setup} We generate the linear regression data $\bA=[\bA_1\ \bA_2]$ and $\by$ randomly, via the function \texttt{gen\_data} (\cref{section:code}), such that \cref{assumption:BWO} is fulfilled and the condition number of $\bA$ can be specified as an input of the function. We implement \ref{eq:BGD} with stepsizes in \cref{eq:gamma1=gamma2-corollary}, and we implement \ref{eq:GD} and \ref{eq:HB} with their optimal stepsizes. Note that \cref{assumption:BWO} puts no restrictions on the spectrum of $\bA$, hence, even under \cref{assumption:BWO}, the optimal stepsizes of  \ref{eq:GD} and \ref{eq:HB} remain to be the same as reviewed in \cref{subsection:GD-HB}.

\myparagraph{Results} In \cref{figa:sr}, we plot under \cref{assumption:BWO} the numerical values of $\rho_{\textnormal{GD}}^*, \rho_{\textnormal{HB}}^*$, and $\rho_{\textnormal{BGD}}^*$. In \cref{fig:convergence10,fig:convergence1000,fig:convergence100000}, we plot the distances $\| \bx^t - \bx^* \|_2$ of the iterates $\{ \bx^t \}_t$ of \ref{eq:GD}, \ref{eq:HB}, and \ref{eq:BGD} for different condition numbers $\kappa:=\frac{\lambda_1(\bA^\top\bA)}{ \lambda_n(\bA^\top\bA) }$. Observe that the numerical convergence rates of  \ref{eq:HB} and \ref{eq:BGD} are following our theory and especially the inequality $\rho_{\textnormal{BGD}}^*\leq (\rho_{\textnormal{HB}}^*)^2$ (\cref{remark:rho-BGD<rho-HB}).

\myparagraph{Discussion} Can the proposed stepsizes for \ref{eq:BGD} be applied to general least-squares problems where $\bA=[\bA_1\ \bA_2]$ do not satisfy \cref{assumption:BWO}? As already indicated in \cref{section:assumption1-AP}, the answer is affirmative, and the idea is that one could always orthogonalize $\bA_1$ and $\bA_2$, respectively, e.g., via QR decomposition $\bA_1=\bQ_1\bR_1$ and $\bA_2=\bQ_2\bR_2$. Indeed, one could then consider minimizing $\| [\bQ_1\ \bQ_2] \bz - \by\|_2$ in variable $\bz$ via \ref{eq:BGD} with the proposed stepsizes as the matrix $[\bQ_1\ \bQ_2]$ now satisfies \cref{assumption:BWO}.  After such minimization, we obtain minimizer $\bz$ and just need to solve the upper triangular systems $\bz_j=\bR_j\bx_j$ for $\bx_j$ ($j=1,2$) by efficient back substitution. But how does this algorithm compare to \ref{eq:HB} applied directly to the data $\bA,\by$? Note that $[\bQ_1\ \bQ_2]$ is, in general, much better conditioned than $\bA$ and its block-wise orthogonality structure allows for faster matrix-vector multiplications $[\bQ_1\ \bQ_2]^\top [\bQ_1\ \bQ_2] \bz$ (e.g., already we have $\bQ_1^\top \bQ_1\bz_1=\bz_1$) than calculating $\bA^\top\bA \bx$. For these two reasons, we observe that combining block-wise QR orthogonalization with \ref{eq:BGD} takes much fewer iterations than \ref{eq:HB} to converge and uses less computation at each iteration; and eventually this method runs faster than \ref{eq:HB} when $\bA$ is highly ill-conditioned (so that \ref{eq:HB} converges slowly) and for medium problem sizes (so that block-wise QR decompositions are affordable). The code and experiments comparing the two approaches are in \cref{subsection:code-w.o.A1}.

% Moreover, $[\bQ_1\ \bQ_2]$ is, in general, much better conditioned than $\bA$ and its block-wise orthogonality structure allows for faster matrix-vector multiplications $[\bQ_1\ \bQ_2] \bz$ (e.g., $\bQ_1^\top \bQ_1\bz_1=\bz_1$ already).

\section{Conclusion and Future Work}\label{section:conclusion}
In this paper, we analyzed a basic algorithm, block gradient descent (\ref{eq:BGD}), applied to a basic setting, least-squares. Under a block-wise orthogonality assumption, we discover optimal stepsizes, with which \ref{eq:BGD} accelerates without any momentum and is provably twice as fast as the heavy ball method (\ref{eq:HB}). Moreover, our results apply to (generalized) alternating projection; Doing so not only leads to improvements over the bounds of \citet{Falt-CDC2017}, but also allows us to revise their inaccurate claim that an optimal convergence rate for (generalized) alternating projection was discovered. Among many interesting questions that might follow from our developments in this paper, we would like to elaborate on the following three points:
\begin{itemize}[wide,parsep=-1pt,topsep=0pt]
    \item Does adding some momentum term to \ref{eq:BGD} give faster rates? \citet{Baptiste-AISTATS2022} showed that GD with cyclic stepsizes and HB momentum is faster than HB if the spectrum of $\bA^\top \bA$ is clustered, so the answer might lie in combining the contributions of \citet{Baptiste-AISTATS2022} and ours.
    \item Can we generalize our analysis from two blocks to multiple blocks? If \ref{eq:BGD} with two blocks can run twice as fast as \ref{eq:HB}, would it be natural---or bold---to guess that \ref{eq:BGD} with $n$ blocks can be $n$ times faster? Moreover, in the case of $n$ blocks, our block-wise orthogonality assumption (\cref{assumption:BWO}) would generalize to the less restrictive requirement that the columns of $\bA$ are normalized; concerns regarding \cref{assumption:BWO} would automatically disappear.
    \item To what extent can we dispense with \cref{assumption:BWO}? Recall that, without this assumption, we need to analyze the spectral radius of the sophisticated matrix $\bM(\gamma_1,\gamma_2)$ defined in \cref{lemma:error-decrease}. While in this general case it is significantly more challenging to find a closed-form solution to the minimum spectral radius of $\bM(\gamma_1,\gamma_2)$, it is perhaps possible to compute a numerical or approximate solution to \cref{eq:minimize-sr}, e.g., via semidefinite relaxation techniques. Pushing this idea to its extreme constitutes an entirely different chapter, and is therefore left as a future endeavor. 
\end{itemize}
\section*{Impact Statement}
This paper presents work intending to advance the field of Machine Learning. There are many potential societal consequences of our work, none of which we feel must be specifically highlighted here.

\bibliography{Liangzu}

\begin{thebibliography}{43}
\providecommand{\natexlab}[1]{#1}
\providecommand{\url}[1]{\texttt{#1}}
\expandafter\ifx\csname urlstyle\endcsname\relax
  \providecommand{\doi}[1]{doi: #1}\else
  \providecommand{\doi}{doi: \begingroup \urlstyle{rm}\Url}\fi

\bibitem[Agarwal et~al.(2021)Agarwal, Goel, and Zhang]{Agarwal-ICML2021}
Agarwal, N., Goel, S., and Zhang, C.
\newblock Acceleration via fractal learning rate schedules.
\newblock In \emph{International Conference on Machine Learning}, pp.\  87--99.
  PMLR, 2021.

\bibitem[Altschuler \& Parrilo(2023{\natexlab{a}})Altschuler and
  Parrilo]{Altschuler-arXiv2023}
Altschuler, J.~M. and Parrilo, P.~A.
\newblock Acceleration by stepsize hedging i: Multi-step descent and the silver
  stepsize schedule.
\newblock Technical report, arXiv:2309.07879 [math.OC], 2023{\natexlab{a}}.

\bibitem[Altschuler \& Parrilo(2023{\natexlab{b}})Altschuler and
  Parrilo]{Altschuler-arXiv2023b}
Altschuler, J.~M. and Parrilo, P.~A.
\newblock Acceleration by stepsize hedging ii: Silver stepsize schedule for
  smooth convex optimization.
\newblock Technical report, arXiv:2309.16530 [math.OC], 2023{\natexlab{b}}.

\bibitem[Aronszajn(1950)]{Aronszajn-1950}
Aronszajn, N.
\newblock Theory of reproducing kernels.
\newblock \emph{Transactions of the American Mathematical Society}, 68\penalty0
  (3):\penalty0 337--404, 1950.

\bibitem[Bauschke et~al.(2014)Bauschke, Cruz, Nghia, Phan, and
  Wang]{Bauschke-JAT2014}
Bauschke, H.~H., Cruz, J.~B., Nghia, T.~T., Phan, H.~M., and Wang, X.
\newblock The rate of linear convergence of the {Douglas-Rachford} algorithm
  for subspaces is the cosine of the {Friedrichs} angle.
\newblock \emph{Journal of Approximation Theory}, 185:\penalty0 63--79, 2014.

\bibitem[Bauschke et~al.(2016)Bauschke, Bello~Cruz, Nghia, Pha, and
  Wang]{Bauschke-NA2016}
Bauschke, H.~H., Bello~Cruz, J., Nghia, T.~T., Pha, H.~M., and Wang, X.
\newblock Optimal rates of linear convergence of relaxed alternating
  projections and generalized {Douglas-Rachford} methods for two subspaces.
\newblock \emph{Numerical Algorithms}, 73:\penalty0 33--76, 2016.

\bibitem[Beck(2017)]{Beck-OptBook2017}
Beck, A.
\newblock \emph{First-Order Methods in Optimization}.
\newblock Society for Industrial and Applied Mathematics, 2017.

\bibitem[Beck \& Tetruashvili(2013)Beck and Tetruashvili]{Beck-SIAM-J-2013}
Beck, A. and Tetruashvili, L.
\newblock On the convergence of block coordinate descent type methods.
\newblock \emph{SIAM Journal on Optimization}, 23\penalty0 (4):\penalty0
  2037--2060, 2013.

\bibitem[Cai \& Diakonikolas(2024)Cai and Diakonikolas]{Cai-arXiv2024}
Cai, X. and Diakonikolas, J.
\newblock Last iterate convergence of incremental methods and applications in
  continual learning.
\newblock Technical report, arXiv:2403.06873 [math.OC], 2024.

\bibitem[Demanet \& Zhang(2016)Demanet and Zhang]{Demanet-MC2016}
Demanet, L. and Zhang, X.
\newblock Eventual linear convergence of the {Douglas-Rachford} iteration for
  basis pursuit.
\newblock \emph{Mathematics of Computation}, 85\penalty0 (297):\penalty0
  209--238, 2016.

\bibitem[Deutsch(1984)]{Deutsch-1984}
Deutsch, F.
\newblock Rate of convergence of the method of alternating projections.
\newblock \emph{Parametric Optimization and Approximation}, 72:\penalty0
  96--107, 1984.

\bibitem[Douglas \& Rachford(1956)Douglas and Rachford]{Douglas-1956}
Douglas, J. and Rachford, H.~H.
\newblock On the numerical solution of heat conduction problems in two and
  three space variables.
\newblock \emph{Transactions of the American mathematical Society}, 82\penalty0
  (2):\penalty0 421--439, 1956.

\bibitem[Elenter et~al.(2023)Elenter, NaderiAlizadeh, Javidi, and
  Ribeiro]{Elenter-arXiv2023}
Elenter, J., NaderiAlizadeh, N., Javidi, T., and Ribeiro, A.
\newblock Primal-dual continual learning: Stability and plasticity through
  lagrange multipliers.
\newblock Technical report, arXiv:2310.00154 [cs.LG], 2023.

\bibitem[Evron et~al.(2022)Evron, Moroshko, Ward, Srebro, and
  Soudry]{Evron-COLT2022}
Evron, I., Moroshko, E., Ward, R., Srebro, N., and Soudry, D.
\newblock How catastrophic can catastrophic forgetting be in linear regression?
\newblock In \emph{Conference on Learning Theory}, 2022.

\bibitem[F{\"a}lt \& Giselsson(2017)F{\"a}lt and Giselsson]{Falt-CDC2017}
F{\"a}lt, M. and Giselsson, P.
\newblock Optimal convergence rates for generalized alternating projections.
\newblock In \emph{IEEE Annual Conference on Decision and Control}, 2017.

\bibitem[Franchetti \& Light(1986)Franchetti and Light]{Franchetti-1986}
Franchetti, C. and Light, W.
\newblock On the {Von Neumann} alternating algorithm in {Hilbert} space.
\newblock \emph{Journal of Mathematical Analysis and Applications},
  114\penalty0 (2):\penalty0 305--314, 1986.

\bibitem[Friedrichs(1937)]{Friedrichs-1937}
Friedrichs, K.
\newblock On certain inequalities and characteristic value problems for
  analytic functions and for functions of two variables.
\newblock \emph{Transactions of the American Mathematical Society}, 41\penalty0
  (3):\penalty0 321--364, 1937.

\bibitem[Goujaud et~al.(2022)Goujaud, Scieur, Dieuleveut, Taylor, and
  Pedregosa]{Baptiste-AISTATS2022}
Goujaud, B., Scieur, D., Dieuleveut, A., Taylor, A.~B., and Pedregosa, F.
\newblock Super-acceleration with cyclical step-sizes.
\newblock In \emph{International Conference on Artificial Intelligence and
  Statistics}, 2022.

\bibitem[Grimmer(2023)]{Grimmer-arXiv2023}
Grimmer, B.
\newblock Provably faster gradient descent via long steps.
\newblock Technical report, arXiv:2307.06324 [math.OC], 2023.

\bibitem[Hesse et~al.(2014)Hesse, Luke, and Neumann]{Hesse-TSP2014}
Hesse, R., Luke, D.~R., and Neumann, P.
\newblock Alternating projections and {Douglas-Rachford} for sparse affine
  feasibility.
\newblock \emph{IEEE Transactions on Signal Processing}, 62\penalty0
  (18):\penalty0 4868--4881, 2014.

\bibitem[Hildreth(1957)]{Hildreth-1957}
Hildreth, C.
\newblock A quadratic programming procedure.
\newblock \emph{Naval Research Logistics Quarterly}, 4\penalty0 (1):\penalty0
  79--85, 1957.

\bibitem[Kayalar \& Weinert(1988)Kayalar and Weinert]{Kayalar-MCSS1988}
Kayalar, S. and Weinert, H.~L.
\newblock Error bounds for the method of alternating projections.
\newblock \emph{Mathematics of Control, Signals and Systems}, 1\penalty0
  (1):\penalty0 43--59, 1988.

\bibitem[Knyazev \& Argentati(2007)Knyazev and
  Argentati]{Knyazev-SIAM-J-MAA2007}
Knyazev, A.~V. and Argentati, M.~E.
\newblock Majorization for changes in angles between subspaces, {Ritz} values,
  and graph {Laplacian} spectra.
\newblock \emph{SIAM Journal on Matrix Analysis and Applications}, 29\penalty0
  (1):\penalty0 15--32, 2007.

\bibitem[Lin et~al.(2023)Lin, Song, and Diakonikolas]{Lin-ICML2023b}
Lin, C.~Y., Song, C., and Diakonikolas, J.
\newblock Accelerated cyclic coordinate dual averaging with extrapolation for
  composite convex optimization.
\newblock In \emph{International Conference on Machine Learning}, 2023.

\bibitem[Nesterov(2012)]{Nesterov-SIAM-J-Opt2012}
Nesterov, Y.
\newblock Efficiency of coordinate descent methods on huge-scale optimization
  problems.
\newblock \emph{SIAM Journal on Optimization}, 22\penalty0 (2):\penalty0
  341--362, 2012.

\bibitem[Nesterov(2018{\natexlab{a}})]{Nesterov-2004}
Nesterov, Y.
\newblock \emph{Introductory Lectures on Convex Optimization: A Basic Course}.
\newblock Springer, 2018{\natexlab{a}}.

\bibitem[Nesterov(2018{\natexlab{b}})]{Nesterov-2018}
Nesterov, Y.
\newblock \emph{Lectures on Convex Optimization}.
\newblock Springer, 2018{\natexlab{b}}.

\bibitem[Nutini et~al.(2015)Nutini, Schmidt, Laradji, Friedlander, and
  Koepke]{Nutini-ICML2015}
Nutini, J., Schmidt, M., Laradji, I., Friedlander, M., and Koepke, H.
\newblock Coordinate descent converges faster with the {G}auss-{S}outhwell rule
  than random selection.
\newblock In \emph{International Conference on Machine Learning}, 2015.

\bibitem[Oymak(2021)]{Oymak-SPL2021}
Oymak, S.
\newblock Provable super-convergence with a large cyclical learning rate.
\newblock \emph{IEEE Signal Processing Letters}, 28:\penalty0 1645--1649, 2021.

\bibitem[Peng \& Vidal(2023)Peng and Vidal]{Peng-arXiv2023}
Peng, L. and Vidal, R.
\newblock Block coordinate descent on smooth manifolds.
\newblock Technical report, arXiv:2305.14744 [math.OC], 2023.

\bibitem[Peng et~al.(2023)Peng, Giampouras, and Vidal]{Peng-ICML2023}
Peng, L., Giampouras, P., and Vidal, R.
\newblock The ideal continual learner: An agent that never forgets.
\newblock In \emph{International Conference on Machine Learning}, 2023.

\bibitem[Peng et~al.(2016)Peng, Wu, Xu, Yan, and Yin]{Peng-AMSA2016}
Peng, Z., Wu, T., Xu, Y., Yan, M., and Yin, W.
\newblock Coordinate-friendly structures, algorithms and applications.
\newblock \emph{Annals of Mathematical Sciences and Applications}, 1\penalty0
  (1):\penalty0 57--119, 2016.

\bibitem[Polyak(1964)]{Polyak-1964}
Polyak, B.~T.
\newblock Some methods of speeding up the convergence of iteration methods.
\newblock \emph{{USSR} Computational Mathematics and Mathematical Physics},
  4\penalty0 (5):\penalty0 1--17, 1964.

\bibitem[Polyak(1987)]{Polyak-1987}
Polyak, B.~T.
\newblock \emph{Introduction to Optimization}.
\newblock Optimization Software, 1987.

\bibitem[Shi et~al.(2016)Shi, Tu, Xu, and Yin]{Shi-arXiv2016v2}
Shi, H.-J.~M., Tu, S., Xu, Y., and Yin, W.
\newblock A primer on coordinate descent algorithms.
\newblock Technical report, arXiv:1610.00040v2 [math.OC], 2016.

\bibitem[Smith(2017)]{Smith-WACV2017}
Smith, L.~N.
\newblock Cyclical learning rates for training neural networks.
\newblock In \emph{IEEE Winter Conference on Applications of Computer Vision},
  2017.

\bibitem[Sun \& Ye(2021)Sun and Ye]{Sun-MP2021}
Sun, R. and Ye, Y.
\newblock Worst-case complexity of cyclic coordinate descent: $o(n^2)$ gap with
  randomized version.
\newblock \emph{Mathematical Programming}, 185:\penalty0 487--520, 2021.

\bibitem[Vidal et~al.(2016)Vidal, Ma, and Sastry]{Vidal-GPCA2016}
Vidal, R., Ma, Y., and Sastry, S.
\newblock \emph{Generalized Principal Component Analysis}.
\newblock Springer, 2016.

\bibitem[von Neumann(1951)]{Von-1951}
von Neumann, J.
\newblock \emph{Functional Operators (AM-22), Volume 2}.
\newblock Princeton University Press, Princeton, 1951.

\bibitem[Wright(2015)]{Wright-MP2015}
Wright, S.~J.
\newblock Coordinate descent algorithms.
\newblock \emph{Mathematical Programming}, 151\penalty0 (1):\penalty0 3--34,
  2015.

\bibitem[Xu \& Yin(2013)Xu and Yin]{Xu-SIAM-J-IS2013}
Xu, Y. and Yin, W.
\newblock A block coordinate descent method for regularized multiconvex
  optimization with applications to nonnegative tensor factorization and
  completion.
\newblock \emph{SIAM Journal on Imaging Sciences}, 6\penalty0 (3):\penalty0
  1758--1789, 2013.

\bibitem[Young(1953)]{Young-JMP1953}
Young, D.
\newblock On {R}ichardson's method for solving linear systems with positive
  definite matrices.
\newblock \emph{Journal of Mathematics and Physics}, 32\penalty0
  (1-4):\penalty0 243--255, 1953.

\bibitem[Zhu \& Knyazev(2013)Zhu and Knyazev]{Zhu-JNM2013}
Zhu, P. and Knyazev, A.~V.
\newblock Angles between subspaces and their tangents.
\newblock \emph{Journal of Numerical Mathematics}, 21\penalty0 (4):\penalty0
  325--340, 2013.

\end{thebibliography}
\bibliographystyle{icml2024}

%%%%%%%%%%%%%%%%%%%%%%%%%%%%%%%%%%%%%%%%%%%%%%%%%%%%%%%%%%%%%%%%%%%%%%%%%%%%%%%
%%%%%%%%%%%%%%%%%%%%%%%%%%%%%%%%%%%%%%%%%%%%%%%%%%%%%%%%%%%%%%%%%%%%%%%%%%%%%%%
% APPENDIX
%%%%%%%%%%%%%%%%%%%%%%%%%%%%%%%%%%%%%%%%%%%%%%%%%%%%%%%%%%%%%%%%%%%%%%%%%%%%%%%
%%%%%%%%%%%%%%%%%%%%%%%%%%%%%%%%%%%%%%%%%%%%%%%%%%%%%%%%%%%%%%%%%%%%%%%%%%%%%%%
\newpage
\appendix
\onecolumn

\section{Elementary Lemmas}
Here we first  give proofs to \cref{lemma:error-decrease,lemma:BGD-AP,lemma:BWO->M,lemma:rho-GD-eigs-C}:
\begin{proof}[Proof of \cref{lemma:error-decrease}]
    Note that the global minimizer $\bx^*$ of \cref{eq:LS2} satisfies the normal equation $\bA^\top \bA \bx^* = \bA^\top \by$, 
    which means $\bA_1^\top \bA \bx^* = \bA_1^\top \by$ and $\bA_2^\top \bA \bx^* = \bA_2^\top \by$, so we can write the first equation $\bx^+_1 = \bx_1 - \gamma_1 \cdot \bA_1^\top  \bA (\bx - \bx^*)$ of \ref{eq:BGD} as
    \begin{equation*}
        \begin{bmatrix}
            \bx^+_1 \\ 
            \bx_2 
        \end{bmatrix} - \bx^* = \Bigg(\bI - \begin{bmatrix}
            \gamma_1 \bA_1^\top \bA_1 & \gamma_1 \bA_1^\top \bA_2 \\ 
            0 & 0
        \end{bmatrix} \Bigg)  (\bx - \bx^*). 
    \end{equation*}
    We finish by re-writing the second equation of \ref{eq:BGD} similarly and combining.
\end{proof}
\begin{proof}[Proof of \cref{lemma:BGD-AP}]
    Left multiplying the first equation of \ref{eq:BGD} by $\bA_1$ and the second equation by $\bA_2$ gives
    \begin{equation*}
        \begin{split}
            \bA_1 \bx_1^+ &= \bA_1 \bx_1 - \gamma_1 \cdot  \bA_1 \bA_1^\top (\bA_1 \bx_1 + \bA_2 \bx_2  - \by ), \\
           \bA_2 \bx_2^+ &= \bA_2 \bx_2 - \gamma_2 \cdot \bA_2 \bA_2^\top ( \bA_1 \bx_1^+ + \bA_2 \bx_2  - \by ).
        \end{split}
    \end{equation*}
    Substitute the expression of $\bA_1 \bx_1^+$ into the second, and sum them up, and we obtain
    \begin{equation*}
        \begin{split}
            \bA \bx^+ &= \bA \bx - \gamma_1 \bA_1 \bA_1^\top (\bA\bx  - \by ) - \gamma_2 \cdot \bA_2 \bA_2^\top \big( \bA_1 \bx_1 - \gamma_1 \cdot  \bA_1 \bA_1^\top (\bA\bx  - \by ) + \bA_2 \bx_2  - \by \big) \\ 
            &=\bA \bx - \gamma_1 \bA_1 \bA_1^\top (\bA\bx  - \by ) - \gamma_2 \cdot \bA_2 \bA_2^\top \big( \bA \bx - \by - \gamma_1 \cdot  \bA_1 \bA_1^\top (\bA\bx  - \by ) \big) \\
            &=\bA \bx - \gamma_1 \bA_1 \bA_1^\top (\bA\bx  - \by ) - \gamma_2 \cdot \bA_2 \bA_2^\top \big( \bI -\gamma_1 \cdot  \bA_1 \bA_1^\top  \big) (\bA\bx  - \by ).
        \end{split}
    \end{equation*}
    As shown in the proof of \cref{lemma:error-decrease} we have $\bA_1^\top \bA \bx^* = \bA_1^\top \by$ and $\bA_2^\top \bA \bx^* = \bA_2^\top \by$. Subtracting $\bA \bx^*$ the above equation  yields
    \begin{equation*}
        \bA (\bx^+ - \bx^*) = \bA(\bx - \bx^*)- \gamma_1 \bA_1 \bA_1^\top (\bA\bx  - \bA \bx^* ) - \gamma_2 \cdot \bA_2 \bA_2^\top \big( \bI -\gamma_1 \cdot  \bA_1 \bA_1^\top  \big) (\bA\bx  - \bA \bx^* ).
    \end{equation*}
    The proof is finished by simplifying the above equation.
\end{proof}

\begin{proof}[Proof of \cref{lemma:BWO->M}]
    Under \cref{assumption:BWO} and with $\bC:=\bA_2^\top \bA_1$, we have
    \begin{equation*}
        \begin{split}
            \bM(\gamma_1,\gamma_2) &= \left(\bI - \begin{bmatrix}
            0 & 0 \\ 
            \gamma_2 \bA_2^\top \bA_1 & \gamma_2 \bA_2^\top \bA_2 
        \end{bmatrix} \Bigg) \Bigg(\bI - \begin{bmatrix}
            \gamma_1 \bA_1^\top \bA_1 & \gamma_1 \bA_1^\top \bA_2 \\ 
            0 & 0
        \end{bmatrix} \right) \\
        &= \left(\bI - \begin{bmatrix}
            0 & 0 \\ 
            \gamma_2 \bC & \gamma_2 \bI_2 
        \end{bmatrix} \Bigg) \Bigg(\bI - \begin{bmatrix}
            \gamma_1 \bI_1 & \gamma_1 \bC^\top \\ 
            0 & 0
        \end{bmatrix} \right) \\
        &= \bI - \begin{bmatrix}
            0 & 0 \\ 
            \gamma_2 \bC & \gamma_2 \bI_2 
        \end{bmatrix} -  \begin{bmatrix}
            \gamma_1 \bI_1 & \gamma_1 \bC^\top \\ 
            0 & 0
        \end{bmatrix} 
        + \begin{bmatrix}
            0 & 0 \\ 
            \gamma_2 \bC & \gamma_2 \bI_2 
        \end{bmatrix} \begin{bmatrix}
            \gamma_1 \bI_1 & \gamma_1 \bC^\top \\ 
            0 & 0
        \end{bmatrix} \\
        &= \begin{bmatrix}
            (1-\gamma_1)\bI_1 & -\gamma_1 \bC^\top \\ 
            -\gamma_2 \bC & (1-\gamma_2) \bI_2 
        \end{bmatrix} + \begin{bmatrix}
            0 & 0 \\ 
            \gamma_1\gamma_2 \bC & \gamma_1\gamma_2\bC \bC^\top 
        \end{bmatrix} \\ 
        &=\begin{bmatrix}
            (1-\gamma_1)\bI_1 & -\gamma_1 \bC^\top \\ 
            -\gamma_2(1-\gamma_1) \bC & (1-\gamma_2) \bI_2 + \gamma_1\gamma_2 \bC \bC^\top 
        \end{bmatrix},
        \end{split}
    \end{equation*}
    and the proof is complete.
\end{proof}

\begin{proof}[Proof of \cref{lemma:rho-GD-eigs-C}]
    By \cref{assumption:BWO} we write $\bA^\top \bA$ as
    \begin{equation*}
        \bA^\top \bA = \begin{bmatrix}
            \bI_1 & \bC^\top \\
            \bC & \bI_2
        \end{bmatrix} = \bI + \bD, \quad \bD:=\begin{bmatrix}
        0 & \bC^\top \\
        \bC & 0
    \end{bmatrix}.
    \end{equation*}
    We will first find the eigenvalues of the $(n_1+n_2)\times (n_1+n_2)$ matrix $\bD$. Since $\bC$ has rank $r$, $\bD$ has eigenvalue $0$ of geometric multiplicity $n_1+n_2 - 2r$, and so the algebraic multiplicity of eigenvalue $0$ is at least $n_1+n_2 - 2r$. With a variable $\lambda\neq 0$ we now look at the characteristic polynomial $\det ( \lambda \bI - \bD )$. By \cref{lemma:BMD} we have
    \begin{equation*}
        \det ( \lambda \bI - \bD ) = \det ( \lambda \bI_1 ) \cdot \det \Big( \lambda \bI_2 - \frac{1}{\lambda} \bC \bC^\top  \Big).
    \end{equation*}
    Therefore, $\det ( \lambda \bI - \bD )=0$ if and only if $\det \big( \lambda^2 \bI_2 - \bC \bC^\top  \big)=0$. Since $\lambda \neq 0$, we have $\det \big( \lambda^2 \bI_2 - \bC \bC^\top  \big)=0$ if and only if $\lambda^2$ is a positive eigenvalue of $\bC\bC^\top$. Since $\bC\bC^\top$ has $r$ such eigenvalues $\lambda_1(\bC\bC^\top),\dots,\lambda_r(\bC\bC^\top)$, we have $2r$ such choices for  $\lambda$, namely $\pm \sqrt{\lambda_1(\bC\bC^\top)},\dots, \pm \sqrt{\lambda_r(\bC\bC^\top)}$. Eigenvalues of $\bA^\top \bA$ are eigenvalues of $\bD$ plus $1$. We finish by recalling \cref{section:background} and simple calculation.
\end{proof}

\cref{lemma:BMD,lemma:basic-minmax,lemma:basic-minmax-2,lemma:basic-root-compare,lemma:basic-2->r,lemma:basic-quadratic} below are elementary, so their proofs are omitted.
\begin{lemma}[Block Matrix Determinant]\label{lemma:BMD}
    We have
    \begin{equation*}
        \det \begin{bmatrix}
        \bQ_1 & \bQ_2 \\ 
        \bQ_3 & \bQ_4
    \end{bmatrix} = \det(\bQ_1) \cdot \det (\bQ_4 - \bQ_3 \bQ_1^{-1} \bQ_2 ),
    \end{equation*}
    where  $\bQ_1$ and $\bQ_4$ are square matrices with $\bQ_1$ invertible.
\end{lemma}
\begin{lemma}\label{lemma:basic-minmax}
    Let $\xi_1\geq \xi_2\geq \dots\geq \xi_r>0$. Then we have
    \begin{equation*}
        \min_{\gamma>0} \Big( \max\big\{ |1- \gamma \xi_1|,\dots, |1- \gamma \xi_r| \big\} \Big) = \begin{cases}
            \frac{\xi_1 - \xi_r}{\xi_1 + \xi_r} & r > 1 ;\\ 
            \frac{1}{\xi_1} & r=1.
        \end{cases}
    \end{equation*}
    Here the minimum is attained at $\gamma=\frac{2}{\xi_1 + \xi_r}$ if $r>1$ or at $\gamma=\frac{1}{\xi_1}$ if $r=1$.
\end{lemma}
\begin{lemma}\label{lemma:basic-minmax-2}
    For functions $f_i:S\to \bbR$ ($\forall i=1,\dots,r$), we have
    \begin{equation*}
        \hat{\gamma} \in \argmin_{\gamma\in S } f_i(\gamma),\forall i=1,\dots,r  \Rightarrow \hat{\gamma}\in \argmin_{\gamma\in S } \Big( \max_{i=1,\dots,r} f_i(\gamma) \Big).
    \end{equation*}
\end{lemma}

\begin{lemma}\label{lemma:basic-root-compare}
    With $\phi\in(0,1)$, $\gamma_1\in(0,1)$, and $\hat{\gamma}_2:= \frac{4-2\gamma_1}{ 2 - \gamma_1\phi}$, the quadratic equation
    \begin{equation*}
        z^2- ( 2-\gamma_1-\hat{\gamma}_2 +\gamma_1\hat{\gamma}_2 \phi )z + (1-\gamma_1)(1-\hat{\gamma}_2) =0
    \end{equation*}
    in variable $z$ has two real roots, $z=\gamma_1-1$ and $z=\hat{\gamma}_2-1$ with $\hat{\gamma}_2-1>1-\gamma_1$.
\end{lemma}

\begin{lemma}\label{lemma:basic-2->r}
    For functions $f_i:S\to \bbR$ ($\forall i=0,1,\dots,r$), if
    \begin{equation*}
        \hat{\gamma} \in \argmin_{\gamma\in S} \Big( \max\{ f_0(\gamma),f_1(\gamma) \}  \Big)
    \end{equation*}
    and $f_i(\hat{\gamma}) \leq \max\{ f_0(\hat{\gamma}),f_1(\hat{\gamma}) \}$ for every $i=1,\dots, r$, then
    \begin{equation*}
        \hat{\gamma} \in \argmin_{\gamma\in S} \Big( \max_{i=0,\dots,r}f_i(\gamma)   \Big).
    \end{equation*}
\end{lemma}

\begin{lemma}\label{lemma:basic-quadratic}
Consider the two (perhaps complex) roots $z_1$ and $z_2$ to the quadratic equation
\begin{equation*}
    z^2 + bz + c = 0
\end{equation*}
 in variable $z$ with $c\geq 0$ and $b\in\bbR $. We have 
 \begin{equation*}
     |z_1|\leq \sqrt{c} \textnormal{\ \ and \ } |z_2|\leq \sqrt{c} \Leftrightarrow |z_1|=|z_2|= \sqrt{c} \Leftrightarrow b^2 -4c\leq 0.
 \end{equation*}
% \begin{proof}
%     Without loss of generality, assume $|z_1|\geq |z_2$. Note that $z_1z_2=c$, which implies the first equivalence. We proceed by proving the second equivalence.

%     Assume $b^2 - 4c\leq 0$. If $b^2 - 4c=0$, then $z_1=z_2$, and therefore $|z_1|=|z_2| = \sqrt{c}$. If $b^2 - 4c<0$, then $z_1$ and $z_2$ are complex numbers, are the conjugate of each other, and have the same magnitude. Therefore, $|z_1|=|z_2| = \sqrt{c}$.

%     On the other hand, assume $|z_1|=|z_2| = \sqrt{c}$. Since $z_1z_2=c$, we have that either $z_1,z_2$ are complex roots (that is $b^2 - 4c < 0$) or $z_1=z_2$ (that is $b^2 - 4c= 0$). 
% \end{proof}
\end{lemma}

\section{Proof of The Main Result (\cref{corollary:main-result}) }

\subsection{A Spectrum Lemma}\label{section:spectrum}
The first step towards solving \cref{eq:minimize-sr} is to analyze the spectrum of $\bM(\gamma_1,\gamma_2)$. We do so in the lemma below. 

\begin{lemma}[Spectrum of $\bM(\gamma_1,\gamma_2)$]\label{lemma:spectrum-M}
    Suppose \cref{assumption:BWO} holds and $\bC=\bA_2^\top\bA_1\in\bbR^{n_2\times n_1}$ has rank $r$. Then, besides the eigenvalues $1-\gamma_1$ and $1-\gamma_2$ shown in \cref{table:eigenvalues}, the remaining eigenvalues of $\bM(\gamma_1,\gamma_2)$ are given as follows:
    \begin{itemize}[wide,parsep=-1pt,topsep=0pt]
        \item (Case 1: $\gamma_1=1$) The remaining $r$ eigenvalues of $\bM(\gamma_1,\gamma_2)$ are $1-\gamma_2 + \gamma_2 \lambda_i(\bC\bC^\top)$, $\forall i=1,\dots,r$;
        \item (Case 2: $\gamma_2=1$) The remaining $r$ eigenvalues of $\bM(\gamma_1,\gamma_2)$ are $1-\gamma_1 + \gamma_1 \lambda_i(\bC\bC^\top)$, $\forall i=1,\dots,r$;
        \item (Case 3: $\gamma_1\neq 1,\gamma_2\neq 1$) The remaining $2r$ eigenvalues of $\bM(\gamma_1,\gamma_2)$ arise as the roots of the following quadratic equation in variable $z$ ($\forall i=1,\dots,r$):
    \begin{equation}\label{eq:quadratic-eigs}
        z^2- \big( 2-\gamma_1-\gamma_2 +\gamma_1\gamma_2 \lambda_i(\bC \bC^\top)  \big)z + (1-\gamma_1)(1-\gamma_2) = 0.
    \end{equation}
   % in variable $z$ are eigenvalues of $\bM(\gamma_1,\gamma_2)$. % ($2r$ eigenvalues arise in this way.)
    \end{itemize}   
\end{lemma}
\begin{table}[t]
    \centering
    \ra{1.1}
    \caption{Algebraic multiplicities of eigenvalues $1-\gamma_1$ and $1-\gamma_2$ of $\bM(\gamma_1,\gamma_2)$. The rank of $\bC$ is denoted by $r$. The table accompanies \cref{lemma:spectrum-M}. }
    \begin{tabular}{cccc}
        \toprule
        & &\multicolumn{2}{c}{Eigenvalues} \\ 
        \cmidrule{3-4}
        & & $1-\gamma_1$ & $1-\gamma_2$ \\
        \midrule
      \multirow{3}{*}{ \begin{turn}{90}Cases\end{turn}}  & $\gamma_1=1$  & $n_1$ & $n_2-r$ \\
        &$\gamma_2=1$  & $n_1-r$ & $n_2$ \\ 
        &$\gamma_1\neq 1,\gamma_2\neq 1$ & $n_1-r$ & $n_2-r$ \\
        \bottomrule
    \end{tabular}
    \label{table:eigenvalues}
\end{table}
\begin{remark}
    The proof of \cref{lemma:spectrum-M} is lengthy, and the difficulty lies in the fact that stepsizes $\gamma_1,\gamma_2$ could change the eigenvalues of $\bM(\gamma_1,\gamma_2)$ in a ``discontinuous'' way (cf. \cref{table:eigenvalues}). And we need to handle all such discontinuities.
\end{remark}
\begin{remark}
    \cref{lemma:spectrum-M}, in conjunction with \cref{table:eigenvalues}, characterizes the spectrum of $\bM(\gamma_1,\gamma_2)$ in different cases. Case 1 ($\gamma_1=1$) and Case 2 ($\gamma_2=1$) turn out to be symmetric, and they coincide if $\gamma_1=\gamma_2=1$. In either case, the eigenvalues of $\bM(\gamma_1,\gamma_2)$ are given explicitly, in terms of $\gamma_1,\gamma_2$, $\lambda_i(\bC\bC^\top)$. In Case 3 ($\gamma_1\neq 1,\gamma_2\neq 1$), we have $2r$ eigenvalues implicitly defined by  quadratic equations \cref{eq:quadratic-eigs}.
\end{remark}
\begin{proof}[Proof of \cref{lemma:spectrum-M}]
    Under \cref{assumption:BWO}, $\bM(\gamma_1,\gamma_2)$ is of the form in \cref{lemma:BWO->M}. We approach the proof by considering the three cases corresponding to \cref{table:eigenvalues} separately.

    \myparagraph{Case 1: $\gamma_1=1$} In this case, $1-\gamma_1=0$ and $\bM(\gamma_1,\gamma_2)$ is given as
    \begin{equation*}
        \bM(1,\gamma_2) = \begin{bmatrix}
            0 & -\bC^\top \\ 
            0 & (1- \gamma_2) \bI_2  + \gamma_2 \bC \bC^\top 
        \end{bmatrix}.
    \end{equation*}
    The eigenvalues of $(1- \gamma_2) \bI_2  + \gamma_2 \bC \bC^\top$ give $n_2$ eigenvalues to $\bM(1,\gamma_2)$; they are $1- \gamma_2$ and $(1-\gamma_2)+\gamma_2 \lambda_i(\bC \bC^\top)$, $i=1,\dots,r$. Note that $1- \gamma_2$ has geometric multiplicity $n_2-r$ matching the dimension of the nullspace of $\bC \bC^\top$. Finally, $0$ is an eigenvalue of $\bM(1,\gamma_2)$ with geometric multiplicity $n_1$. Since all $n_1+n_2$ eigenvalues are counted. the abovementioned geometric multiplicity coincides with algebraic multiplicity. Similarly, in the rest of the proof, the eigenvalues will be counted using either geometric or algebraic ways, but by finding all eigenvalues (counting geometric or algebraic multiplicity) we will eventually reveal that in our case the geometric multiplicity turns out to coincide with algebraic multiplicity. Hence, in what follows we will just write \textit{multiplicity} for simplicity.

    \myparagraph{Case 2: $\gamma_2=1$} In this case, $1-\gamma_2=0$ and $\bM(\gamma_1,\gamma_2)$ is given as
    \begin{equation*}
        \bM(\gamma_1,1) = \begin{bmatrix}
            (1-\gamma_1) \bI_1 & -\gamma_1 \bC^\top \\ 
            -(1-\gamma_1) \bC &  \gamma_1 \bC\bC^\top
        \end{bmatrix}.    
    \end{equation*}
    % We need to consider whether $\gamma_1=1$ or not. 
    If $\gamma_1=1$, by \textit{Case 1} we know that $0=1-\gamma_1=1-\gamma_2$ is an eigenvalue of $\bM(\gamma_1,1)$ with multiplicity $n_1+n_2-r$, and the remaining $r$ eigenvalues are given as $\lambda_i(\bC \bC^\top)$, $i=1,\dots,r$; \cref{table:eigenvalues} is correct. So in the remaining proof of \textit{Case 2} we assume $\gamma_1\neq 1$. 
    
    \myparagraph{Case 2.1: Eigenvalue $1-\gamma_1$}  Let us now test whether $\bM(\gamma_1,1)$ has eigenvalue $1-\gamma_1$. With an $(n_1+n_2)$-dimensional vector $\bv:=[\bv_1;\bv_2]$ we have 
    \begin{equation*}
        \begin{split}
            \bM(\gamma_1,1) \bv = (1-\gamma_1)\bv &\Leftrightarrow \begin{cases}
                (1-\gamma_1)\bv_1 - \gamma_1\bC^\top \bv_2 = (1-\gamma_1)\bv_1 \\
                -(1-\gamma_1)\bC\bv_1  + \gamma_1\bC\bC^\top \bv_2 =(1-\gamma_1)\bv_2
            \end{cases} \\
            &\Leftrightarrow \begin{cases}
                \bC^\top \bv_2 = 0 \\
                \bC \bv_1 =-\bv_2
            \end{cases}
        \end{split} 
    \end{equation*}
    and this implies $\bC^\top \bC\bv_1=-\bC^\top \bv_2=0$. Hence $\bv_2=\bC\bv_1=0$. On the other hand, if $\bC\bv_1=\bv_2=0$ then the above is satisfied. We now conclude that $\bM(\gamma_1,1) \bv = (1-\gamma_1)\bv$ if and only if $\bC\bv_1=\bv_2=0$. The nullspace of $\bC$ is of dimension $n_1-r$, so $(1-\gamma_1)$ is an eigenvalue of $\bM(\gamma_1,1)$ of multiplicity $n_1-r$.

    \myparagraph{Case 2.2: Eigenvalue 0} Following a similar proof of \textit{Case 2.1}, we have
    \begin{equation*}
        \begin{split}
            \bM(\gamma_1,1) \bv = 0 &\Leftrightarrow \begin{cases}
                (1-\gamma_1)\bv_1 - \gamma_1\bC^\top \bv_2 = 0 \\
                -(1-\gamma_1)\bC\bv_1  + \gamma_1\bC\bC^\top \bv_2 =0
            \end{cases} \\
            &\Leftrightarrow (1-\gamma_1)\bv_1 - \gamma_1\bC^\top \bv_2 = 0.
        \end{split} 
    \end{equation*}
    Since $\gamma_1\neq 1$, for any $\bv_2\in\bbR^{n_2}$ we can set $\bv_1$ to be $\gamma_1 \bC^\top \bv_2/(1-\gamma_1)$ so that the equality $\bM(\gamma_1,1) \bv = 0$ holds. Therefore $0$ is an eigenvalue of $\bM(\gamma_1,1)$ of multiplicity $n_2$. 
    
    \myparagraph{Case 2.3: The Remaining $r$ Eigenvalues} Consider  the characteristic polynomial $\det \big(z\bI - \bM(\gamma_1,1) \big)$ in variable $z$, with $z\neq 1-\gamma_1$ and $z\neq 0$. \cref{lemma:BMD} implies
    \begin{equation*}
        \begin{split}
            \det\big(z \bI - \bM(\gamma_1,1) \big) = \det \Big( z \bI_1 - (1-\gamma_1) \bI_1  \Big) \cdot & \det \Big( z \bI_2 - \gamma_1 \bC \bC^\top \\ 
        &- (1-\gamma_1)\bC \cdot \frac{1}{z - (1-\gamma_1)} \cdot \gamma_1 \bC^\top \Big).
        \end{split}
    \end{equation*}
    With $z\neq 1-\gamma_1$ and $z\neq 0$, the above implies  
    \begin{equation}\label{eq:det=0-case2.3}
         \det\big(z \bI - \bM(\gamma_1,1) \big)=0 \Leftrightarrow \det \Big( z \bI_2 - (1-\gamma_1)\bI_2- \gamma_1 \bC \bC^\top \Big) = 0.
    \end{equation}
    We now see that for each $i=1,\dots,r$, $(1-\gamma_1)+ \gamma_1\lambda_i (\bC\bC^\top)$ is a solution to \cref{eq:det=0-case2.3}, and is therefore an eigenvalue of $\bM(\gamma_1,1)$.

    \myparagraph{Case 3: $\gamma_1\neq 1,\gamma_2\neq 1$} Suppose for the sake of contradiction that $\bM(\gamma_1,\gamma_2) \bv=0$ for some $\bv=[\bv_1;\bv_2]\in\bbR^{n_1+n_2}$. Then, using the formula of $\bM(\gamma_1,\gamma_2)$ in \cref{lemma:BWO->M} we can write
    \begin{equation*}
        \begin{split}
            &\begin{cases}
                (1-\gamma_1)\bv_1 - \gamma_1\bC^\top \bv_2 = 0\\
            -\gamma_2(1-\gamma_1)\bC\bv_1 + (1-\gamma_2) \bv_2 + \gamma_1\gamma_2 \bC\bC^\top \bv_2 =0
            \end{cases} \\ 
            \Leftrightarrow& \begin{cases}
                (1-\gamma_1)\bv_1 - \gamma_1\bC^\top \bv_2 = 0\\
            -\gamma_2(1-\gamma_1)\bC\bv_1 + (1-\gamma_2) \bv_2 + \gamma_2 \bC (1-\gamma_1)\bv_1 =0
            \end{cases} \\
            \Leftrightarrow& \begin{cases}
                (1-\gamma_1)\bv_1 + \gamma_1\bC^\top \bv_2 = 0\\
            (1-\gamma_2) \bv_2 =0
            \end{cases}
        \end{split}
    \end{equation*}
    and since $\gamma_1\neq 1$ and $\gamma_2\neq 1$, it is now clear that $\bv_1=\bv_2=0$. So $\bM(\gamma_1,\gamma_2)$ is full rank.
    
    We will show there are $2r$ eigenvalues different than $1-\gamma_1$ and $1-\gamma_2$ and given by the roots of \cref{eq:quadratic-eigs}. Consider the characteristic polynomial $\det\big(z\bI - \bM(\gamma_1,\gamma_2)  \big)$ in variable $z$, with $z\neq 1-\gamma_1$, $z\neq 1-\gamma_2$, and $z\neq 0$. By \cref{lemma:BMD} we have
        \begin{equation*}
        \begin{split}
            \det\big(z \bI - \bM(\gamma_1,\gamma_2) \big) = \det \Big( z \bI_1 - (1-\gamma_1) \bI_1  \Big) \cdot &\det \Big( z \bI_2 - (1-\gamma_2) \bI_2 - \gamma_1\gamma_2 \bC \bC^\top \\ 
        &- \gamma_2(1-\gamma_1)\bC \frac{1}{z - (1-\gamma_1)} \cdot \gamma_1 \bC^\top \Big).
        \end{split}
    \end{equation*}
    Simplifying the above equation yields
    \begin{equation*}% \label{eq:det-reduced}
        \begin{split}
            \det\big(z \bI - \bM(\gamma_1,\gamma_2) \big)=0 \Leftrightarrow \det\Big( \frac{\big(z - (1-\gamma_2) \big) \cdot \big(z - (1-\gamma_1) \big) }{z}\bI_2 -  \gamma_1\gamma_2 \bC\bC^\top \Big) = 0.
        \end{split}
    \end{equation*}
    The characteristic polynomial $\det (z'\bI_2 - \gamma_1\gamma_2 \bC\bC^\top  ) $ in variable $z'$ has $r$ nonzero roots, namely $\gamma_1\gamma_2 \lambda_1(\bC \bC^\top), \dots, \gamma_1\gamma_2 \lambda_r(\bC \bC^\top)$. Moreover, for each root $z'$, the equation 
    \begin{equation*}
        \frac{\big(z - (1-\gamma_2) \big) \cdot \big(z - (1-\gamma_1) \big) }{z} = z' \Leftrightarrow z^2- \big( 2-\gamma_1 - \gamma_2+ z'\big)z + (1-\gamma_1)(1-\gamma_2)= 0 
    \end{equation*}
    always has two nonzero (potentially complex) solutions, and there are $2r$ such solutions in total; none of these solutions is equal to $1-\gamma_1$ or $1-\gamma_2$. In particular, by construction, these $2r$ solutions must  be eigenvalues of $\bM(\gamma_1,\gamma_2)$.
    
    It remains to show  $1-\gamma_1$ and $1-\gamma_2$ are eigenvalues of $\bM(\gamma_1,\gamma_2)$ with mulpliticity $n_1 - r$ and $n_2 -r$, respectively; cf. \cref{table:eigenvalues}. We consider two cases, $\gamma_1=\gamma_2$ and $\gamma_1 \neq \gamma_2$.

    \myparagraph{Case 3.1: $\gamma_1\neq \gamma_2$} With some vector  $\bv:=[\bv_1;\bv_2]\in\bbR^{n_1+n_2}$ we have 
    \begin{equation*}
        \begin{split}
            \bM(\gamma_1,\gamma_2) \bv = (1-\gamma_1)\bv 
            % &\Leftrightarrow \begin{cases}
            %     (1-\gamma_1)\bv_1 + \gamma_1\bC^\top \bv_2 = (1-\gamma_1)\bv_1 \\
            %     \gamma_2(1-\gamma_1)\bC\bv_1 + (1-\gamma_2) \bv_2 + \gamma_1\gamma_2 \bC\bC^\top \bv_2 =(1-\gamma_1)\bv_2
            % \end{cases} \\ 
            &\Leftrightarrow \begin{cases}
                -\gamma_1\bC^\top \bv_2 = 0\\
                -\gamma_2(1-\gamma_1)\bC\bv_1 + (\gamma_1-\gamma_2) \bv_2 + \gamma_1\gamma_2 \bC\bC^\top \bv_2 =0
            \end{cases} \\
            &\Leftrightarrow \begin{cases}
                \bC^\top \bv_2 = 0; \\
                -\gamma_2(1-\gamma_1)\bC\bv_1 + (\gamma_1-\gamma_2) \bv_2 =0.
            \end{cases}
        \end{split} 
    \end{equation*}
    Left multiplying the last equation by $\bC^\top$ yields $\gamma_2(1-\gamma_1)\bC^\top\bC\bv_1=0$. Since $\gamma_1\neq 1$ and $\gamma_2>0$, we obtain $\bC \bv_1=0$. Substituting it back gives $(\gamma_1-\gamma_2)\bv_2=0$, but $\gamma_1\neq \gamma_2$, so it must be that $\bv_2=0$. Therefore, for $\bM(\gamma_1,\gamma_2) \bv = (1-\gamma_1)\bv$ to hold, it is necessary that $\bv_2=0$ and $\bC \bv_1=0$; one verifies this is also sufficient, and thus % If $\bv_2=0$ and $\bC \bv_1=0$, then we have $\bM(\gamma_1,\gamma_2) \bv = (1-\gamma_1)\bv$. To summarize, we have 
    \begin{equation*}
        \bM(\gamma_1,\gamma_2) \bv = (1-\gamma_1)\bv \Leftrightarrow \bC \bv_1=0, \bv_2=0.
    \end{equation*}
    Since the nullspace of $\bC$ has dimension $n_1-r$, so $1-\gamma_1$ is an eigenvalue of $\bM(\gamma_1,\gamma_2)$ of multiplicity $n_1-r$. Proving $1-\gamma_2$ is an eigenvalue of $\bM(\gamma_1,\gamma_2)$ of multiplicity $n_2-r$ follows a similar route, detailed next. With $\bv:=[\bv_1;\bv_2]\in\bbR^{n_1+n_2}$ we have
    \begin{equation*}
        \begin{split}
            \bM(\gamma_1,\gamma_2) \bv = (1-\gamma_2)\bv 
            % &\Leftrightarrow \begin{cases}
            %     (1-\gamma_1)\bv_1 + \gamma_1\bC^\top \bv_2 = (1-\gamma_2)\bv_1 \\
            %     \gamma_2(1-\gamma_1)\bC\bv_1 + (1-\gamma_2) \bv_2 + \gamma_1\gamma_2 \bC\bC^\top \bv_2 =(1-\gamma_2)\bv_2
            % \end{cases} \\ 
            &\Leftrightarrow \begin{cases}
                -\gamma_1\bC^\top \bv_2 = (\gamma_1 -\gamma_2)\bv_1\\
                -(1-\gamma_1)\bC\bv_1  + \gamma_1 \bC\bC^\top \bv_2 =0
            \end{cases} \\
            &\Leftrightarrow \begin{cases}
                -\gamma_1\bC^\top \bv_2 = (\gamma_1 -\gamma_2)\bv_1\\
                -(1-\gamma_1)\bC\bv_1  - \bC \cdot (\gamma_1-\gamma_2)\bv_1 =0
            \end{cases} \\ 
            &\Leftrightarrow \begin{cases}
                \gamma_1\bC^\top \bv_2 = (\gamma_1 -\gamma_2)\bv_1; \\
                (1-\gamma_2)\bC\bv_1 =0.
            \end{cases}
        \end{split} 
    \end{equation*}
    Since $\gamma_2\neq 1$, the last equation implies $\bC\bv_1 =0$. Left multiply the equality $\gamma_1\bC^\top \bv_2 = (\gamma_1 -\gamma_2)\bv_1$ by $\bC$ and use $\bC\bv_1 =0$, and we obtain $\bC \bC^\top \bv_2=0$, that is $\bC^\top \bv_2=0$. Substitute this back and we get $(\gamma_1 -\gamma_2)\bv_1=0$. But $\gamma_1\neq \gamma_2$, so $\bv_1=0$. Thus, $\bM(\gamma_1,\gamma_2) \bv = (1-\gamma_2)\bv$ implies $\bv_1=0$ and $\bC^\top \bv_2=0$, and one verifies that the converse is also true. Therefore, $1-\gamma_2$ is an eigenvalue of $\bM(\gamma_1,\gamma_2)$ with multiplicity $n_2-r$, that is the dimension of the nullspace of $\bC^\top$.
    
    \myparagraph{Case 3.2: $\gamma_1=\gamma_2$} In this case, we consider the linear equations $\bM(\gamma_1,\gamma_1) \bv = (1-\gamma_1)\bv$ in variable $\bv:=[\bv_1;\bv_2]\in\bbR^{n_1+n_2}$. And we have
    \begin{equation*}
        \begin{split}
            \bM(\gamma_1,\gamma_1) \bv = (1-\gamma_1)\bv 
            % &\Leftrightarrow \begin{cases}
            %     (1-\gamma_1)\bv_1 + \gamma_1\bC^\top \bv_2 = (1-\gamma_1)\bv_1 \\
            %     \gamma_1(1-\gamma_1)\bC\bv_1 + (1-\gamma_1) \bv_2 + \gamma_1^2 \bC\bC^\top \bv_2 =(1-\gamma_1)\bv_2
            % \end{cases} \\ 
            &\Leftrightarrow \begin{cases}
                \bC^\top \bv_2 = 0  \\
                \gamma_1(1-\gamma_1) \bC \bv_1 + \gamma_1^2 \bC\bC^\top \bv_2 = 0
            \end{cases} \\
            &\Leftrightarrow 
            \bC_2 \bv = 0, \quad \textnormal{where\ \ } \bC_2:=\begin{bmatrix}
                \bC & 0 \\ 
                0 & \bC^\top 
            \end{bmatrix}.
        \end{split}
    \end{equation*}
    The nullspace of $\bC_2$ is of dimension $n_1+n_2 - 2r$, so $1-\gamma_1$ is an eigenvalue of $\bM(\gamma_1,\gamma_2)$ of multiplicity $n_1+n_2 - 2r$. The proof is now complete.
\end{proof}

\cref{lemma:spectrum-M} is important as it allows us to proceed for studying the spectral radius $\rho \big( \bM(\gamma_1,\gamma_2) \big)$. That is not to say analyzing $\rho \big( \bM(\gamma_1,\gamma_2) \big)$ is made into a trivial task. In fact, \cref{table:eigenvalues} already suggests an obstacle: If $n_1=r$ and $\gamma_1\neq 1$, then $1-\gamma_1$ is not an eigenvalue of $\bM(\gamma_1,\gamma_2)$ and should not be taken into account when studying $\rho \big( \bM(\gamma_1,\gamma_2) \big)$; similarly for $1-\gamma_2$. In other words, this suggests we have to further break our analysis into subcases that classify the relations between the rank $r$ and block size $n_1$ (or $n_2$). Performing such analysis carefully is our main duty in \cref{section:simplified-case,section:general-case}.

% To better appreciate the underlying challenge of analyzing $\rho \big( \bM(\gamma_1,\gamma_2) \big)$ and to provide intuition to do so, we consider simplified examples in  \cref{section:simplified-case}. Then, in \cref{section:general-case}, we tackle the problem of minimizing $\rho \big( \bM(\gamma_1,\gamma_2) \big)$ in its full generality and discover the optimal stepsizes for block coordinate descent.

% Specifically, in \cref{section:simplified-case} we study $\rho \big( \bM(\gamma_1,\gamma_2) \big)$ through examples, and in \cref{section:general-case} we minimize $\rho \big( \bM(\gamma_1,\gamma_2) \big)$ and discover the optimal stepsizes for block coordinate descent.

\subsection{Simplified Cases}\label{section:simplified-case} 
In this section we make some basic calculations of the spectral radius in simplified cases. The purpose is to build some intuition and preliminary results for the sequel. The first easy case we consider is where $\gamma_1=1$ or $\gamma_2=1$. In this case, \cref{lemma:spectrum-M} reveals the whole spectrum of $\bM(\gamma_1,\gamma_2)$, so we can minimize $\rho \big( \bM(\gamma_1, \gamma_2) \big)$ with relative ease:
\begin{proposition}\label{prop:gamma1=1-gamma2=1}
    Recall $\bC:=\bA_2^\top\bA_1$ and $r:=\textnormal{rank}(\bC)$. Under \cref{assumption:BWO}, $\rho \big( \bM(\gamma_1, 1) \big)$ is minimized at $\gamma_1= \frac{2}{2 - \lambda_1(\bC \bC^\top)}$ if $r<n_1$, at  $\gamma_1 =\frac{2}{2 - \lambda_1(\bC \bC^\top)-\lambda_r(\bC \bC^\top)} $ if $r=n_1>1$, and its minimum values are given as
    \begin{equation*}
        \begin{split}
            \min_{\gamma_1>0} \rho \big( \bM(\gamma_1, 1) \big) = \begin{cases}
            \frac{\lambda_1(\bC \bC^\top)}{2 - \lambda_1(\bC \bC^\top)} & \textnormal{if } r<n_1; \\ 
            \frac{\lambda_1(\bC \bC^\top)- \lambda_r(\bC \bC^\top) }{2 - \lambda_1(\bC \bC^\top)- \lambda_r(\bC \bC^\top) } & \textnormal{if } r=n_1.
        \end{cases} 
        \end{split}
    \end{equation*}
\end{proposition}
\begin{remark}\label{remark:gamma1=1-gamma2=1}
    Similarly, we can solve $\min_{\gamma_2>0} \rho \big( \bM(1,\gamma_2) \big)$.
\end{remark}
\begin{remark}\label{remark:sp=0}
    If $r=n_1=1$, then $\lambda_1(\bC \bC^\top) = \lambda_r(\bC \bC^\top)$, in which case the minimum of $\rho \big( \bM(\gamma_1, 1) \big)$ is $0$. 
    % Note though that this does not necessarily mean \ref{eq:BGD} converges in a single iteration; it only implies the behavior of \ref{eq:BGD} in the limit: If $r=n_1=1$ and \cref{assumption:BWO} holds, then
    % \begin{equation*}
    % \tiny{
    %     \begin{split}
    %     \lim_{t\to \infty} \Big\| \bM\Big(\frac{1}{1-\lambda_1(\bC\bC^\top)},1\Big)^t \Big\|_2^{1/t} =  \rho\Big(\bM\Big(\frac{1}{1-\lambda_1(\bC\bC^\top)},1\Big)\Big) 
    %     = 0.
    %     \end{split} 
    %     }
    % \end{equation*}
    % The first equality is a property of spectral radii.
\end{remark}
\begin{proof}[Proof of \cref{prop:gamma1=1-gamma2=1}] 
\cref{lemma:spectrum-M} implies
    \begin{equation*}
        \begin{split}
            \rho\big( \bM(\gamma_1,1) \big) = \begin{cases}
                \max \Big\{ |1-\gamma_1|, \max_{i=1,\dots,r} \big\{ |1 - \gamma_1 \big( 1 - \lambda_i(\bC \bC^\top) \big) | \big\} \Big\} & \textnormal{if } r<n_1; \\
                \max_{i=1,\dots,r} \big\{ |1 - \gamma_1 \big( 1 - \lambda_i(\bC \bC^\top) \big) | \big\} & \textnormal{if } r=n_1.
            \end{cases}
        \end{split}
    \end{equation*}
    The proof finishes with a basic and standard argument; see, e.g., \cref{lemma:basic-minmax}.
\end{proof}

Another easy case is where  $\gamma_1=\gamma_2$. In this case, we need to minimize $\rho \big( \bM(\gamma_1,\gamma_2) \big)$ only in a single variable $\gamma_1$: 
\begin{proposition}\label{prop:sr-gamma1=gamma2}
    Define $\gamma_1^* := \frac{2}{1+ \sqrt{1-\lambda_1(\bC \bC^\top)}}$. Under \cref{assumption:BWO}, we have $\gamma_1^*-1 < \rho \big( \bM(1,1)\big)$ and 
    \begin{equation}\label{eq:sr-gamma1=gamma2}
        \min_{\gamma_1>0} \rho \big( \bM(\gamma_1,\gamma_1) \big) = \rho \big( \bM(\gamma_1^*,\gamma_1^*) \big)= \gamma_1^*-1.
    \end{equation}
\end{proposition}
\begin{remark}
    While in \cref{prop:gamma1=1-gamma2=1} the minima can vary with the rank $r$, the minimum of \cref{eq:sr-gamma1=gamma2} is attained at $\gamma_1^*$ for any $r$.
\end{remark}
\begin{proof}[Proof of \cref{prop:sr-gamma1=gamma2}]
    Under \cref{assumption:BWO}, $\gamma_1^*$ is real-valued, and with some calculation, one further verifies $0<\gamma_1^*-1<\lambda_1(\bC\bC^\top)=\rho \big( \bM(1,1)\big)$; the last equality is due to \cref{example:BEM}.
    We proceed in two steps. In Step 1 we prove $\rho \big( \bM(\gamma_1^*,\gamma_1^*) \big)= \gamma_1^*-1$. In Step 2 we prove $\gamma_1^*-1 = \min_{\gamma_1>0} \rho \big( \bM(\gamma_1,\gamma_1) \big)$.

    \myparagraph{Step 1} First observe that for every $i=1,\dots,r$ we have  
    \begin{equation*}
        \begin{split}
            &\ \Big( 2-2\gamma_1^* +(\gamma_1^*)^2 \lambda_i(\bC \bC^\top)  \Big)^2 - 4(\gamma_1^*-1)^2 \leq 0  \\ 
            % \Leftrightarrow&\ \big( 2-2\gamma_1^* +(\gamma_1^*)^2 \lambda_i(\bC \bC^\top) + 2 ( \gamma_1^*-1)   \big) \big( 2-2\gamma_1^* +(\gamma_1^*)^2 \lambda_i(\bC \bC^\top) - 2 ( \gamma_1^*-1) \big) \leq 0 \\ 
            \Leftrightarrow&\ (\gamma_1^*)^2 \lambda_i(\bC \bC^\top)     \cdot \Big( 4-4\gamma_1^* +(\gamma_1^*)^2 \lambda_i(\bC \bC^\top)  \Big) \leq 0 \\
            \Leftrightarrow&\ 4-4\gamma_1^* +(\gamma_1^*)^2 \lambda_i(\bC \bC^\top) \leq 0. 
        \end{split}
    \end{equation*}
    But as one can verify, the definition of $\gamma_1^*$ implies $\gamma_1^*$ is the smaller root of the equation $4-4\xi +\xi^2 \lambda_1(\bC \bC^\top) = 0$ in variable $\xi$, that is $4-4\gamma_1^* +(\gamma_1^*)^2 \lambda_1(\bC \bC^\top) = 0$. As a consequence,  for every $i=1,\dots,r$, the roots of the quadratic equation
    \begin{equation*}
        z^2- \Big( 2-2\gamma_1^* +(\gamma_1^*)^2 \lambda_i(\bC \bC^\top)  \Big)z + (\gamma_1^*-1)^2 =0
    \end{equation*} 
    in variable $z$ have magnitudes equal to $\gamma_1^*-1$. Invoking \cref{lemma:spectrum-M}, we see that these roots are precisely $2r$ eigenvalues $\bM(\gamma_1^*,\gamma_1^*)$, while \cref{lemma:spectrum-M} further suggests that  $\bM(\gamma_1^*,\gamma_1^*)$ might have eigenvalue $\gamma_1^*-1$ with multiplicity $n_1+n_2-2r$ (if $n_1+n_2-2r\neq 0$). We can now conclude \textit{Step 1} with $\rho \big( \bM(\gamma_1^*,\gamma_1^*) \big)= \gamma_1^*-1$.

    \myparagraph{Step 2} To prove $\gamma_1^*-1 = \min_{\gamma_1>0} \rho \big( \bM(\gamma_1,\gamma_1) \big)$, we consider two cases.

    \myparagraph{Case 2.1: $\gamma_1>\gamma_1^*$} In this case, $\gamma_1>\gamma_1^*>1$. Note that for every $i=1,\dots,r$ the equation \cref{eq:quadratic-eigs} must has a root whose magnitude is larger than or equal to $\sqrt{(1-\gamma_1)(\gamma_2-1)}=\gamma_1-1$. Therefore $\rho \big( \bM(\gamma_1,\gamma_1) \big)\geq \gamma_1 - 1 > \gamma_1^*-1$.
    % \cref{lemma:spectrum-M} suggests there must be one eigenvalue of $\bM(\gamma_1,\gamma_1)$ arising as a root of \cref{eq:quadratic-eigs} that has magnitude larger than or equal to $\gamma_1-1$. Therefore $\rho \big( \bM(\gamma_1,\gamma_1) \big)\geq \gamma_1 - 1 > \gamma_1^*-1$.

    \myparagraph{Case 2.2: $0<\gamma_1\leq\gamma_1^*$} The choice of $\gamma_1^*$ indicates $ 4-4\gamma_1 +\gamma_1^2 \lambda_1(\bC \bC^\top) \geq 0$. We then prove $2-2\gamma_1 +\gamma_1^2 \lambda_1(\bC \bC^\top)\geq 0$: This is true if $\gamma_1\in(0,1]$; otherwise, if $\gamma_1>1$, this is implied by $ 4-4\gamma_1 +\gamma_1^2 \lambda_1(\bC \bC^\top) >0$. As a result, the quadratic equation $ z^2- \big( 2-2\gamma_1 +\gamma_1^2 \lambda_1(\bC \bC^\top)  \big)z + (\gamma_1-1)^2 =0$ has two different non-negative real roots and they are eigenvalues of $\bM(\gamma_1,\gamma_1)$. The larger root, denoted by $z(\gamma_1)$, is given as
    \begin{equation*}
        \begin{split}
          z(\gamma_1) =&\ \frac{2-2\gamma_1 +\gamma_1^2 \lambda_1(\bC \bC^\top)+ \sqrt{\big(2-2\gamma_1 +\gamma_1^2 \lambda_1(\bC \bC^\top) \big)^2 - 4(\gamma_1-1)^2}  }{2} \\ 
       % =&\ \frac{4-4\gamma_1 +2 \gamma_1^2 \lambda_1(\bC \bC^\top)+ 2\sqrt{  \gamma_1^2 \lambda_1(\bC \bC^\top) \cdot  \big(4-4\gamma_1 +\gamma_1^2 \lambda_1(\bC \bC^\top) \big) }  }{4} \\ 
        =&\ \frac{\Big( \gamma_1 \sqrt{\lambda_1(\bC \bC^\top)} + \sqrt{4-4\gamma_1 +\gamma_1^2 \lambda_1(\bC \bC^\top)} \Big)^2  }{4}.
        \end{split}
    \end{equation*}
    We need to prove $z(\gamma_1)\geq \gamma_1^*-1$ for every $\gamma_1\in(0,\gamma_1^*]$. Since $z(\gamma_1^*)=\gamma_1^*-1$, it suffices to show that the function 
    \begin{equation*}
        f(\xi) = \xi \sqrt{\lambda_1(\bC \bC^\top)} + \sqrt{4-4\xi +\xi^2 \lambda_1(\bC \bC^\top)}
    \end{equation*}
    is non-increasing in $[0,\gamma_1^*]$, i.e., its derivative is non-positive. Noting that, with $\xi\in [0,\gamma_1^*]$ we have $2-\xi \lambda_1(\bC \bC^\top)  \geq 2- \gamma_1^* \lambda_1(\bC \bC^\top)  \geq 0$. We can then verify
    \begin{equation*}
        \begin{split}
            f'(\xi) \leq 0 &\Leftrightarrow \sqrt{\lambda_1(\bC \bC^\top)} + \frac{ \xi \lambda_1(\bC \bC^\top) -2}{\sqrt{4-4\xi +\xi^2 \lambda_1(\bC \bC^\top)}} \leq 0 \\ 
            &\Leftrightarrow \lambda_1(\bC \bC^\top) \cdot \Big( 4-4\xi +\xi^2 \lambda_1(\bC \bC^\top) \Big) \leq \Big(2 - \xi \lambda_1(\bC \bC^\top) \Big)^2 \\ 
            &\Leftrightarrow \lambda_1(\bC \bC^\top) \leq 1.
        \end{split}
    \end{equation*}
    We have thus completed the proof.
\end{proof}

From \cref{prop:gamma1=1-gamma2=1,prop:sr-gamma1=gamma2}, we see that, under \cref{assumption:BWO}, stepsizes better than $\gamma_1=\gamma_2=1$ do exist, suggesting that it is possible for \ref{eq:BGD} to converge faster than \ref{eq:BEM}. 

Note that the minima in \cref{prop:gamma1=1-gamma2=1,prop:sr-gamma1=gamma2} are both attained at some stepsize larger than $1$, behind which the intuition is as follows. Suppose \cref{assumption:BWO} holds, then the maximum eigenvalue of $\bC\bC^\top$ is smaller than or equal to $1$. If furthermore $\gamma_1=1$, then $\bM(\gamma_1,\gamma_2)$ can be written as $\bM(1,\gamma_2)=\begin{bmatrix}
            0 & -\bC^\top \\ 
            0 & (1- \gamma_2) \bI_2  + \gamma_2 \bC \bC^\top 
        \end{bmatrix}$. 
If $\gamma_2$ were restricted to lie in $(0,1]$, then there would be no doubt that the convex combination $(1- \gamma_2) \bI_2  + \gamma_2 \bC \bC^\top$ would be minimized at $\gamma_2=1$ under \cref{assumption:BWO}. However, inspect that setting $\gamma_2$ to be appropriately larger than $1$ would actually further reduce the magnitudes of $(1- \gamma_2) \bI_2  + \gamma_2 \bC \bC^\top$ and therefore of its eigenvalues. Optimal stepsizes that we show in \cref{section:general-case} are in fact all larger than $1$.  
\subsection{General Cases}\label{section:general-case}
\begin{figure*}[h]
    \centering
    \begin{tikzpicture}
    
    \draw [-to] (0,0) -> (2.5,0) node [at end, below right] {$\gamma_1$};
    
    \draw [-to] (0,0) -> (0,2.5) node [at end, left] {$\gamma_2$};

    \draw (0,1) -- (2,1);
    \draw (0,2) -- (2,2);

    \draw (1,0) -- (1,2);
    \draw (2,0) -- (2,2);

    \draw (1,0) node [below] {$1$};
    \draw (2,0) node [below] {$\infty$};
    
    \draw (0,1) node [left] {$1$};
    \draw (0,2) node [left] {$\infty$};
    \draw (-0.165,-0.16) node {$0$};

    \draw (0.5,0.5) node {$S_{00}$};
    \draw (0.5,1.5) node {$S_{01}$};
    \draw (1.5,0.5) node {$S_{10}$};
    \draw (1.5,1.5) node {$S_{11}$};

    \draw (-2.5, 0.5) node {$S_{00}:=(0,1]\times (0,1]$};
    \draw (-2.5, 1.5) node {$S_{01}:=(0,1]\times [1,\infty)$};
    \draw (4.5, 0.5) node {$S_{10}:=[1,\infty)\times (0,1]$};
    \draw (4.5, 1.5) node {$S_{11}:=[1,\infty)\times [1,\infty)$};
    
    \end{tikzpicture}
    \caption{We divide the quadrant $\{ (\gamma_1,\gamma_2):\gamma_1>0,\gamma_2>0 \}$ of all possible stepsizes into $4$ regions, $S_{00},S_{01}, S_{10}$, and $S_{11}$. We minimize the spectral radius $\rho\big(\bM(\gamma_1,\gamma_2) \big)$ over each region separately, which will give a solution to \cref{eq:minimize-sr}. }  
    \label{fig:stepsizes-S}
\end{figure*}
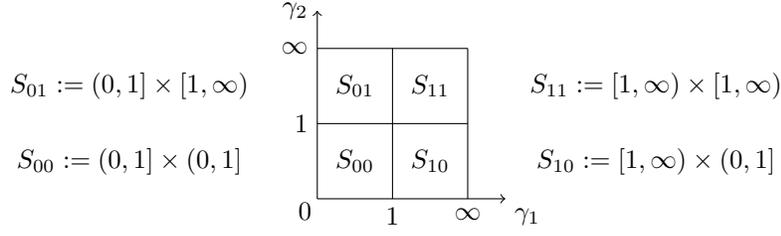

The limitation of \cref{prop:gamma1=1-gamma2=1,prop:sr-gamma1=gamma2} is that, there, we search stepsizes on two rays, namely $\{ (\gamma_1,1):\gamma_1>0 \}$ and $\{ (\gamma_1,\gamma_1):\gamma_1>0 \}$, while a better choice might lie on the \textit{quadrant} $\{(\gamma_1,\gamma_2):\gamma_1>0,\gamma_2>0\}$ yet on neither of the two rays. To address this point, we divide the quadrant into four regions, $S_{00},S_{01},S_{10},S_{11}$ as shown in \cref{fig:stepsizes-S}, and we search for stepsizes that leading to smaller spectral radii over each region. Combining yields a solution to \cref{eq:minimize-sr}.

To proceed, we develop a technical lemma, shown below, which will be used as a common sub-routine to prove our main results in \cref{subsection:stepsizes-S00-S01-S10,subsection:stepsizes-S11}. 
\begin{lemma}\label{lemma:distriminant-gamma1>1-gamma1<1}
    Fix $\phi\in(0,1)$. Consider the following quadratic equation in variable $z$ and its discriminant: 
    \begin{align}\label{eq:equation+discriminant}
        %\begin{split}
           & z^2 - (2 - \gamma_1 - \gamma_2 + \gamma_1 \gamma_2 \phi)z + (1-\gamma_1) (1-\gamma_2) = 0, \nonumber \\ 
        &\Delta(\gamma_1,\gamma_2):=(2 - \gamma_1 - \gamma_2 + \gamma_1 \gamma_2 \phi)^2 \\ 
        &\quad \quad \quad \quad \quad \ \  - 4(1-\gamma_1) (1-\gamma_2). \nonumber
        %\end{split}
    \end{align}
    Assume $\Delta(\gamma_1,\gamma_2)\geq 0$. Then
    \begin{equation*}
        \begin{cases}
            2\sqrt{\Delta(\gamma_1,\gamma_2)} \cdot |\gamma_1\phi-1| \leq \big| \frac{\partial \Delta(\gamma_1,\gamma_2)}{ \partial \gamma_2 } \big| & \textnormal{if } \gamma_1 \geq 1; \\ 
            2\sqrt{\Delta(\gamma_1,\gamma_2)} \cdot |\gamma_1\phi-1| \geq \big| \frac{\partial \Delta(\gamma_1,\gamma_2)}{ \partial \gamma_2 } \big| & \textnormal{if } \gamma_1 \in(0,1).
        \end{cases}
    \end{equation*}
\end{lemma}
\begin{remark}
    In words, \cref{lemma:distriminant-gamma1>1-gamma1<1} analyzes a single quadratic equation  \cref{eq:equation+discriminant} which is of the form \cref{eq:quadratic-eigs}, highlighting some properties useful for analyzing the eigenvalues as roots of \cref{eq:quadratic-eigs}  in the sequel.
\end{remark}
\begin{proof}[Proof of \cref{lemma:distriminant-gamma1>1-gamma1<1}]
        We only write down the proof for the case $\gamma_1\geq 1$; the case for $\gamma_1\in(0,1)$ is similar. First note that 
    \begin{equation*}
            2\sqrt{\Delta(\gamma_1,\gamma_2)} \cdot |\gamma_1\phi-1|  \leq \big|\frac{\partial \Delta(\gamma_1,\gamma_2)}{ \partial \gamma_2 } \big| \Leftrightarrow 
            4 \Delta(\gamma_1,\gamma_2) \cdot (\gamma_1\phi-1)^2 \leq \Big( \frac{\partial \Delta(\gamma_1,\gamma_2)}{ \partial \gamma_2 } \Big)^2. 
    \end{equation*}
    Canceling a common additive term yields the following equivalent inequality:
    \begin{equation*}
        -(\gamma_1-1)(\gamma_2-1) (\gamma_1\phi-1)^2 \leq  (\gamma_1-1)^2 - (2-\gamma_1-\gamma_2 + \gamma_1\gamma_2\phi)(\gamma_1 \phi - 1) (\gamma_1-1).
    \end{equation*}
    If $\gamma_1= 1$, we are done. If $\gamma_1>1$, we can divide both sides of the above by $\gamma_1-1$ to get an equivalent inequality: % (if $\gamma_1<1$ we need to switch the direction of the inequality)
    \begin{equation*}
        \begin{split}
            &\ -(\gamma_2-1) (\gamma_1\phi-1)^2 \leq  (\gamma_1-1) - (2-\gamma_1-\gamma_2 + \gamma_1\gamma_2\phi)(\gamma_1 \phi - 1)  \\ 
            \Leftrightarrow&\ (\gamma_1 \phi - 1)\cdot \big( (2-\gamma_1-\gamma_2 + \gamma_1\gamma_2\phi) - (\gamma_2-1) (\gamma_1\phi-1)  \big) \leq (\gamma_1-1) \\ 
            \Leftrightarrow&\ (\gamma_1 \phi - 1)\cdot ( 1-\gamma_1 + \gamma_1\phi ) \leq (\gamma_1-1) \\ 
            \Leftrightarrow&\ \gamma_1^2\phi(\phi-1) \leq 0.
        \end{split}
    \end{equation*}
    But $\phi\in(0,1)$, so $\gamma_1^2\phi(\phi-1) \leq 0$, and we finished the proof.
    \end{proof}

\subsubsection{Optimal Stepsizes on $S_{00}\cup S_{01}\cup S_{10}$}\label{subsection:stepsizes-S00-S01-S10}
Let us first acquire a better understanding of how the roots of a single quadratic equation in \cref{eq:equation+discriminant} behave:
\begin{lemma}\label{lemma:1quad-gamma12}
    Fix $\phi\in(0,1)$. Consider the quadratic equation in \cref{eq:equation+discriminant} and let
    $z_1(\gamma_1,\gamma_2)$, $z_2(\gamma_1,\gamma_2)$ be its two roots. Define  
    \begin{equation*}%\label{eq:single-f-gamma12}
        f(\gamma_1,\gamma_2):= \max \big\{ |z_1 (\gamma_1,\gamma_2)| , |z_2 (\gamma_1,\gamma_2)| \big\}.
    \end{equation*}
    
    With $S_{00}$, $S_{01}$, $S_{10}$ defined in \cref{fig:stepsizes-S}, the following hold:
    \begin{itemize}[wide,parsep=-1pt,topsep=0pt]
        \item Part 1: $S_{00}$. With $\gamma_1\in(0,1]$ fixed, $\gamma_2\mapsto f(\gamma_1,\gamma_2)$ is non-increasing in $(0,1]$. With $\gamma_2\in(0,1]$ fixed, $\gamma_1\mapsto f(\gamma_1,\gamma_2)$ is non-increasing in $(0,1]$. Therefore
        \begin{equation*}%\label{eq:[01]x[01]}
            \min_{[\gamma_1;\gamma_2]\in S_{00}} f(\gamma_1,\gamma_2)=f(1,1).
        \end{equation*}
        \item Part 2: $S_{01}\cup S_{10}$. With $\gamma_2\in[1,\infty)$ fixed, $\gamma_1\mapsto f(\gamma_1,\gamma_2)$ is non-increasing in $(0,1]$. With $\gamma_1\in[1,\infty)$ fixed, $\gamma_2\mapsto f(\gamma_1,\gamma_2)$ is non-increasing in $(0,1]$. Hence we have
        \begin{equation*}
            \begin{split}
                \min_{[\gamma_1;\gamma_2]\in S_{01}} f(\gamma_1,\gamma_2)=\min_{\gamma_2\geq 1} f(1,\gamma_2), \\ \min_{[\gamma_1;\gamma_2]\in S_{10}} f(\gamma_1,\gamma_2)=\min_{\gamma_1\geq 1} f(\gamma_1,1).
            \end{split}
        \end{equation*}
       % \item (Part 3: $S_{10}$) With $\gamma_2\in(0,1]$ fixed, the function $\gamma_1\mapsto f(\gamma_1,\gamma_2)$ is non-increasing in $[1,\frac{2-\gamma_2}{1-\gamma_2\phi})$ and non-decreasing in $[\frac{2-\gamma_2}{1-\gamma_2\phi}, 2)$.
    \end{itemize}    
\end{lemma}

\begin{proof}[Proof of \cref{lemma:1quad-gamma12}]
    Recall $\Delta(\gamma_1,\gamma_2):=(2 - \gamma_1 - \gamma_2 + \gamma_1 \gamma_2 \phi)^2 - 4(1-\gamma_1) (1-\gamma_2)$ (\cref{lemma:distriminant-gamma1>1-gamma1<1}). We present the proofs below for \textit{Part 1} and \textit{Part 2}, one after another.
    
    \myparagraph{Part 1: $S_{00}$} Assume $\gamma_1,\gamma_2\in(0,1]$, and we derive the monotonicity of $f$. In this case we have $1-\gamma_1\geq 0$ and $1-\gamma_2\geq 0$, and therefore we can write
    \begin{equation*}
        \Delta(\gamma_1,\gamma_2) = \bigg( \Big(\sqrt{1-\gamma_1} - \sqrt{1-\gamma_2} \Big)^2 + \gamma_1 \gamma_2 \phi \bigg)  \bigg( \Big(\sqrt{1-\gamma_1} + \sqrt{1-\gamma_2} \Big)^2 + \gamma_1 \gamma_2 \phi \bigg), 
    \end{equation*}
    which implies $\Delta(\gamma_1,\gamma_2)>0$. Furthermore, since $2 - \gamma_1 - \gamma_2 + \gamma_1 \gamma_2 \phi\geq 0$, we have 
    \begin{equation*}
        \frac{\partial f(\gamma_1,\gamma_2)}{\partial\gamma_1} = \frac{ -1 +\gamma_2\phi }{2} + \frac{1 }{4 \sqrt{\Delta(\gamma_1,\gamma_2)}} \cdot \frac{\partial \Delta(\gamma_1,\gamma_2) }{\partial \gamma_1}.
    \end{equation*}
    Since $\gamma_2\in(0,1)$, applying \cref{lemma:distriminant-gamma1>1-gamma1<1} with $\gamma_1$ and $\gamma_2$ swapped yields $2\sqrt{\Delta(\gamma_1,\gamma_2)} \cdot |\gamma_1\phi-1| \geq \big| \frac{\partial \Delta(\gamma_1,\gamma_2)}{ \partial \gamma_1 } \big|$. Since $-1+\gamma_2\phi \leq 0$, this further implies $\frac{\partial f(\gamma_1,\gamma_2)}{\partial\gamma_1}\leq 0$. Similarly, using  \cref{lemma:distriminant-gamma1>1-gamma1<1} we can show  $\frac{\partial f(\gamma_1,\gamma_2)}{\partial\gamma_2}\leq 0$. This proves \textit{Part 1}.
    
    \myparagraph{Part 2: $S_{01}\cup S_{10}$}  Since $(1-\gamma_1)(1-\gamma_2)\leq 0$, the two roots of \cref{eq:equation+discriminant} are real-valued. It suffices to fix $\gamma_1\in[1,\infty)$ and prove $\gamma_2\mapsto f(\gamma_1,\gamma_2)$ is non-increasing in $(0,1]$. To do so, we first calculate $\frac{\partial \Delta(\gamma_1,\gamma_2)}{\partial \gamma_2}$ and verify $\frac{\partial \Delta(\gamma_1,\gamma_2)}{\partial \gamma_2} \leq 0$ for every $\gamma_2\in(0,1]$:
    \begin{equation*}
        \begin{split}
            \frac{\partial \Delta(\gamma_1,\gamma_2)}{\partial \gamma_2} \leq 0 \Leftrightarrow&\  2(\gamma_1\phi-1)^2\gamma_2 - 2( \gamma_1\phi+ 1 -2\phi )\gamma_1 \leq 0 \\ 
            \Leftarrow&\  (\gamma_1\phi-1)^2 - ( \gamma_1\phi+ 1 -2\phi )\gamma_1 \leq 0 \\
            \Leftrightarrow&\ \gamma_1^2\phi(\phi-1) + 1 - \gamma_1 \leq 0.
        \end{split}
    \end{equation*}
    This indeed holds as $\phi\in(0,1)$ and $\gamma_1\geq 1$. Then, by definition we can write
    \begin{equation*}
        \frac{\partial f(\gamma_1,\gamma_2)}{ \partial \gamma_2 } = \begin{cases}
            \frac{-1 + \gamma_1\phi}{2}+ \frac{1}{4\sqrt{\Delta(\gamma_1,\gamma_2)}} \cdot \frac{\partial \Delta(\gamma_1,\gamma_2)}{ \partial \gamma_2 } & \quad  \textnormal{if \ } 2-\gamma_1-\gamma_2 + \gamma_1\gamma_2\phi > 0; \\
            \frac{1 - \gamma_1\phi}{2}+ \frac{1 }{4\sqrt{\Delta(\gamma_1,\gamma_2)}} \cdot \frac{\partial \Delta(\gamma_1,\gamma_2)}{ \partial \gamma_2 } & \quad \textnormal{if \ } 2-\gamma_1-\gamma_2 + \gamma_1\gamma_2\phi < 0; \\ 
            \sqrt{(\gamma_1-1)(1-\gamma_2)}  & \quad \textnormal{if \ } 2-\gamma_1-\gamma_2 + \gamma_1\gamma_2\phi = 0.
        \end{cases}
    \end{equation*}
    Here, $\gamma_2\mapsto (\gamma_1-1)(1-\gamma_2)$ is clearly a decreasing function. For the other two cases, we can also prove $\frac{\partial f(\gamma_1,\gamma_2)}{ \partial \gamma_2 }\leq 0$ by using \cref{lemma:distriminant-gamma1>1-gamma1<1} and the fact $\frac{\partial \Delta(\gamma_1,\gamma_2)}{\partial \gamma_2} \leq 0$.
\end{proof}

Having studied a single quadratic equation in \cref{lemma:1quad-gamma12}, we can now analyze $r$ such equations \cref{eq:quadratic-eigs}. Specifically, armed with \cref{lemma:distriminant-gamma1>1-gamma1<1,lemma:1quad-gamma12}, we reach the following result:
\begin{theorem}\label{theorem:S00-S01-S10} Under \cref{assumption:BWO}, the following hold.
    \begin{itemize}[wide,parsep=-1pt,topsep=0pt]
        \item (Part 1: $S_{00}$) On $S_{00}$, we have % the spectral radius $\rho \big( \bM(\gamma_1,\gamma_2) \big)$ is minimized at $[1;1]$:
        \begin{equation*}
            \min_{[\gamma_1;\gamma_2]\in S_{00}} \rho \big( \bM(\gamma_1,\gamma_2) \big) = \rho \big( \bM(1,1) \big).
        \end{equation*}
        \item (Part 2: $S_{01}\cup S_{10}$) We have
        \begin{equation*}
            \begin{split}
                 \min_{[\gamma_1;\gamma_2]\in S_{01}} \rho \big( \bM(\gamma_1,\gamma_2) \big) = \min_{\gamma_2\geq 1} \rho \big( \bM(1,\gamma_2) \big), \\ 
            \min_{[\gamma_1;\gamma_2]\in S_{10}} \rho \big( \bM(\gamma_1,\gamma_2) \big) = \min_{\gamma_1\geq 1} \rho \big( \bM(\gamma_1,1) \big) .
            \end{split}
        \end{equation*}
    \end{itemize}
\end{theorem}
\begin{remark}
    The situation on $S_{00}$ is now clear  (\textit{Part 1}): The minimum spectral radius is attained at $(1,1)$. The situation on $S_{01} \cup S_{10}$ is clearly only partially (\textit{Part 2}), and it seems that one still needs to minimize over a single stepsize $\gamma_2$ (\textit{resp}. $\gamma_1$) to find the minimum spectral radius on $S_{01}$ (\textit{resp.} $S_{10}$), with the other stepsize set to $1$. This is in fact an easier task and has been addressed already in \cref{prop:gamma1=1-gamma2=1}.
\end{remark}
\begin{proof}[Proof of \cref{theorem:S00-S01-S10}]
    The proof of \cref{theorem:S00-S01-S10} relies on two technical lemmas, \cref{lemma:distriminant-gamma1>1-gamma1<1,lemma:1quad-gamma12}.

    Let $z_{i1}(\gamma_1,\gamma_2)$ and $z_{i2}(\gamma_1,\gamma_2)$ be the two roots of \cref{eq:quadratic-eigs}. Define 
    \begin{equation}\label{eq:fi}
        f_i(\gamma_1,\gamma_2):= \max\big\{ |z_{i1}(\gamma_1,\gamma_2)|, |z_{i2}(\gamma_1,\gamma_2)| \big\}.
    \end{equation}
    If $r<\min\{ n_1,n_2 \}$, then, by \cref{lemma:spectrum-M}, $1-\gamma_1$ and $1-\gamma_2$ are eigenvalues of $\bM(\gamma_1,\gamma_2)$, and \cref{lemma:1quad-gamma12,lemma:basic-minmax-2} further imply
    \begin{equation*}
        \begin{split}
            \min_{[\gamma_1;\gamma_2]\in S_{00}} \rho \big( \bM(\gamma_1,\gamma_2) \big) &= \min_{[\gamma_1;\gamma_2]\in S_{00} } \bigg( \max \Big\{ 1-\gamma_1, 1-\gamma_2, \max_{i,\dots,r} f_i(\gamma_1,\gamma_2) \Big\} \bigg) \\ 
        &= \rho \big( \bM(1,1) \big).
        \end{split}
    \end{equation*}
    The case where $r=n_1$ or $r=n_2$ is similar, and this finishes the proof of \textit{Part 1}. \textit{Part 2} follows similarly from \cref{lemma:basic-minmax-2} and \textit{Part 2} of \cref{lemma:1quad-gamma12,}.
\end{proof}

\subsubsection{Optimal Stepsizes on $S_{11}$}\label{subsection:stepsizes-S11}
Solving \cref{eq:minimize-sr} on $S_{11}$ is harder than on $S_{00}\cup S_{01} \cup S_{10}$, as on $S_{11}$ it is more challenging to quantify the change of the roots (eigenvalues) defined in \cref{eq:quadratic-eigs} with respect to stepsize $\gamma_1$ or $\gamma_2$. We tackle this challenge by further dividing our analysis into two separate cases: $\bC$ has full rank and $\bC$ is rank-deficient.

The proofs for the two cases are much more sophisticated. Hence we only state the theorems below. The proofs are presented separately in \cref{proof-theorem-S11}.

%we address the case where $\bC$ is full rank (i.e., $r=\min\{n_1,n_2\}$):

\myparagraph{The Full Rank Case} If $\bC$ is full rank (i.e., $r=\min\{n_1,n_2\}$), we have:
\begin{theorem}\label{theorem:S11-fullrank}
    Suppose $r=\min\{n_1,n_2\}$ and \cref{assumption:BWO} holds. Let $S_{11}$ be defined in \cref{fig:stepsizes-S}. Then the minimum value  $\min_{[\gamma_1;\gamma_2]\in S_{11} } \rho\big( \bM(\gamma_1,\gamma_2) \big)$ is equal to 
    \begin{equation}\label{eq:min-sr-S11-fullrank}
        \frac{\sqrt{1-\lambda_r(\bC\bC^\top)} - \sqrt{1-\lambda_1(\bC\bC^\top)}  }{\sqrt{1-\lambda_r(\bC\bC^\top)} + \sqrt{1-\lambda_1(\bC\bC^\top)}}.
    \end{equation}    
    % If $r=1$ then
    % \begin{equation*}
    %     \min_{[\gamma_1;\gamma_2]\in S_{11} } \rho\big( \bM(\gamma_1,\gamma_2) \big) = 0.
    % \end{equation*}
    % If $r>1$ then $\min_{[\gamma_1;\gamma_2]\in S_{11} } \rho\big( \bM(\gamma_1,\gamma_2) \big)$ is equal to
    % \begin{equation*}
    %     \frac{\sqrt{1-\lambda_r(\bC\bC^\top)} - \sqrt{1-\lambda_1(\bC\bC^\top)}  }{\sqrt{1-\lambda_r(\bC\bC^\top)} + \sqrt{1-\lambda_1(\bC\bC^\top)}}.
    % \end{equation*}
    % \begin{equation*}
    %     \min_{[\gamma_1;\gamma_2]\in S_{11} } \rho\big( \bM(\gamma_1,\gamma_2) \big) = \begin{cases}
    %         0 & r = 1 \\ 
    %         \frac{\sqrt{1-\lambda_r(\bC\bC^\top)} - \sqrt{1-\lambda_1(\bC\bC^\top)}  }{\sqrt{1-\lambda_r(\bC\bC^\top)} + \sqrt{1-\lambda_1(\bC\bC^\top)}} & r > 1
    %     \end{cases}
    % \end{equation*}
\end{theorem}
\begin{remark}
    If $r=1$, then we have $\lambda_1(\bC\bC^\top)= \lambda_r(\bC\bC^\top)$ and \cref{theorem:S11-fullrank} suggests the minimum spectral radius is $0$, coinciding with \cref{prop:gamma1=1-gamma2=1,remark:sp=0}.   
\end{remark}

\myparagraph{The Rank-Deficient Case} It remains to address the case where $\bC$ is rank-deficient (i.e., $r<\min\{n_1,n_2\}$). We do this in \cref{theorem:S11-rankdeficient}:
\begin{theorem}\label{theorem:S11-rankdeficient}
    Suppose $r < \min\{n_1,n_2\}$ and \cref{assumption:BWO} holds. With $S_{11}$ defined in \cref{fig:stepsizes-S} and $\gamma_1^*$ defined in \cref{prop:sr-gamma1=gamma2}, we have $\rho \big( \bM(\gamma_1^*,\gamma_1^*) \big) = \gamma_1^*-1$ and
    \begin{equation}\label{eq:min-sr-S11-rank-deficient}
        \begin{split}
            \min_{[\gamma_1;\gamma_2]\in S_{11}} \rho \big( \bM(\gamma_1,\gamma_2) \big) &= \gamma_1^*-1 \\
        &= \frac{1 - \sqrt{1-\lambda_1(\bC\bC^\top)}}{1+ \sqrt{1-\lambda_1(\bC\bC^\top)} }.
        \end{split}
    \end{equation}
\end{theorem}
There is a basic connection between \cref{eq:min-sr-S11-fullrank} and \cref{eq:min-sr-S11-rank-deficient}: They would be identical if $\lambda_r(\bC\bC^\top)$ were equal to $0$. Despite this direct connection,  their proofs are very different.

\section{Proofs of \cref{theorem:S11-fullrank,theorem:S11-rankdeficient}}\label{proof-theorem-S11}
\subsection{Proof of \cref{theorem:S11-fullrank}}\label{subsection:proof-S11-fullrank}
For the proof of \cref{theorem:S11-fullrank}, we need the following lemma.
\begin{lemma}\label{lemma:HB}
    Assume $s>1$ and $\zeta_1\geq \cdots \geq \zeta_s>0$.  For $i=1,\dots,s$, let $z_{i1}(\alpha,\beta)$, $z_{i2}(\alpha,\beta)$ be the two roots of the following quadratic equation in variable $z$:
    \begin{equation}\label{eq:quadratic-alpha-beta}
        z^2 - (\beta + 1 - \alpha \zeta_i )z + \beta = 0.
    \end{equation}
    Define $g_i(\alpha,\beta):=\max\{ |z_{i1}(\alpha,\beta)|, |z_{i2}(\alpha,\beta)| \}$. Then
    \begin{equation*}
        \min_{\alpha>0,\beta\geq 0} \Big( \max_{i=1,\dots,s} g_i(\alpha,\beta) \Big) = \frac{ \sqrt{\zeta_1} - \sqrt{\zeta_s} }{ \sqrt{\zeta_1} + \sqrt{\zeta_s} },
    \end{equation*}
    where the minimum is attained at
    \begin{equation*}
        \begin{split}
            \alpha = \Big( \frac{ 2 }{ \sqrt{\zeta_1} + \sqrt{\zeta_s} } \Big)^2, \quad \beta = \Big( \frac{ \sqrt{\zeta_1} - \sqrt{\zeta_s} }{ \sqrt{\zeta_1} + \sqrt{\zeta_s} } \Big)^2.
        \end{split}
    \end{equation*}
\end{lemma} 
\begin{proof}[Proof of \cref{lemma:HB}]
    By the definition of \cref{eq:quadratic-alpha-beta} we have 
    \begin{equation*}
        |z_{i1}(\alpha,\beta)| \cdot |z_{i2}(\alpha,\beta)| = \beta.
    \end{equation*}
    Without loss of generality, we assume $|z_{i1}(\alpha,\beta)| \geq |z_{i2}(\alpha,\beta)|$ If $z_{i1}(\alpha,\beta)$ and $z_{i1}(\alpha,\beta)$ are two complex roots, then their magnitudes are equal to $\sqrt{\beta}$. If they are real roots, then we must have $|z_{i1}(\alpha,\beta)|\geq \sqrt{\beta}$, where the equality is attained if and only if the two real roots are equal. Therefore $g_i(\alpha,\beta)\geq \sqrt{\beta}$. As a consequence, for any $\beta\geq 0$ the optimization problem  
    \begin{equation}\label{eq:min-alpha-only}
        \min_{\alpha>0} \Big( \max_{i=1,\dots,s} g_i(\alpha,\beta) \Big)
    \end{equation}
    is lower bounded by $\sqrt{\beta}$. Clearly, the lower bound $\sqrt{\beta}$ of \cref{eq:min-alpha-only} is attained at $\alpha$ if and only if $\alpha$ satisfies the following for all $i=1,\dots,s$:
    \begin{equation}\label{eq:alpha_optimal_condition}
        \begin{split}
           &\ (\beta + 1 - \alpha \zeta_i )^2 - 4\beta \leq 0 \\ 
           \Leftrightarrow&\ \big(\beta + 1 - \alpha \zeta_i +2 \sqrt{\beta} \big)\big(\beta + 1 - \alpha \zeta_i  - 2 \sqrt{\beta} \big) \leq 0 \\
           \Leftrightarrow&\ \big( (\sqrt{\beta}+1)^2 - \alpha \zeta_i \big)\big( (\sqrt{\beta}-1)^2
 - \alpha \zeta_i\big) \leq 0 \\ 
            \Leftrightarrow&\  \big| 1 - \sqrt{\alpha \zeta_i} \big| \leq  \sqrt{\beta} \leq  1 + \sqrt{\alpha \zeta_i}.
        \end{split}
    \end{equation}
    Next, for $\beta$ to be as small as possible while satisfying \cref{eq:alpha_optimal_condition}, we need to choose $\alpha$ to be the global minimizer of 
    \begin{equation*}
        \min_{\alpha>0} \Big\{ \max_{i=1,\dots,s}  \big| 1 - \sqrt{\alpha \zeta_i} \big| \Big\} = \min_{\alpha>0} \Big\{ \max \big\{ \big| 1 - \sqrt{\alpha \zeta_1} \big|, \big| 1 - \sqrt{\alpha \zeta_s} \big|  \big\} \Big\}.
    \end{equation*}
    A simple geometric argument shows that the minimum of the above problem is attained at $\alpha = \frac{4}{(\sqrt{\zeta_1} + \sqrt{\zeta_s})^2}$, with the associated minimum value being $\frac{ \sqrt{\zeta_1} - \sqrt{\zeta_s} }{ \sqrt{\zeta_1} + \sqrt{\zeta_s} }$. To summarize, with $\alpha^* = \frac{4}{(\sqrt{\zeta_1} + \sqrt{\zeta_s})^2}$ and $\beta^* = \frac{ ( \sqrt{\zeta_1} - \sqrt{\zeta_s})^2 }{ (\sqrt{\zeta_1} + \sqrt{\zeta_s})^2 }$, we have  $\max_{i=1,\dots,s} g_i(\alpha^*,\beta^*)$ equal to $\sqrt{\beta^*}$.

    To prove optimality, we need to show that for any $\alpha, \beta$ that satisfies $\sqrt{\beta} < \big| 1 - \sqrt{\alpha \zeta_i} \big|$ for some $i$, the corresponding objective is larger than $\sqrt{\beta^*}$. To do so, it suffices to consider the following cases:
    \begin{itemize}
        \item If $\alpha>\alpha^*$, then $|1- \sqrt{\alpha \zeta_s}| < |1-\sqrt{\alpha \zeta_1}| $. Assume $\sqrt{\beta} <|1-\sqrt{\alpha \zeta_1}|$. The equation $z^2 - (\beta + 1 - \alpha \zeta_1 )z + \beta = 0$ has two real roots, and the larger one is
        \begin{equation*}
            \frac{ |\beta + 1 - \alpha \zeta_1| + \sqrt{ (\beta + 1 - \alpha \zeta_1)^2 - 4\beta}  }{2}. 
        \end{equation*}
        The derivative of this root with respect to $\beta$ is 
        \begin{equation*}
            \pm \frac{1}{2} + \frac{\beta  - \alpha \zeta_1 - 1}{2 \sqrt{ (\beta + 1 - \alpha \zeta_1)^2 - 4\beta} },
        \end{equation*}
        which is negative regardless of the sign on $1/2$ (use the assumption $\sqrt{\beta} <|1-\sqrt{\alpha \zeta_i}|$ to verify). Hence the root is minimized when $\sqrt{\beta} = |1-\sqrt{\alpha \zeta_1}|$, at which point this root can be written as
        \begin{equation*}
            \frac{ (1-\sqrt{\alpha \zeta_1})^2 + 1 - \alpha \zeta_1 }{2} = |1-\sqrt{\alpha \zeta_1}|.
        \end{equation*}
        Since $\alpha>\alpha^*$, one verifies $|1-\sqrt{\alpha \zeta_1}|$ is larger than $|1-\sqrt{\alpha^* \zeta_1}|= \sqrt{\beta^*}$. 
        \item The case $\alpha< \alpha^*$ can be proved similarly. Specifically, in this case we have $|1- \sqrt{\alpha \zeta_1}| < |1-\sqrt{\alpha \zeta_s}| $, so we can assume $\sqrt{\beta} <|1-\sqrt{\alpha \zeta_s}|$ and proceed with a similar argument.
    \end{itemize}
\end{proof}

We are then ready to prove \cref{theorem:S11-fullrank}.
\begin{proof}[Proof of \cref{theorem:S11-fullrank}]
    We can easily verify from \cref{lemma:spectrum-M} that
    \begin{equation*}
        \begin{split}
            \rho\Big( \bM\Big( \frac{1}{1-\lambda_1(\bC\bC^\top)},1\Big) \Big) = 0 &\quad \textnormal{if \ }  r=n_1=1; \\ 
            \rho\Big( \bM\Big(1, \frac{1}{1-\lambda_1(\bC\bC^\top)}\Big) \Big) = 0 &\quad \textnormal{if \ } r=n_2=1.
        \end{split}
    \end{equation*}
    Next we assume $r>1$. To proceed, we consider the function $\tau:S_{11}\to \bbR^2$ defined as
    \begin{equation}\label{eq:tau}
        \tau(\gamma_1,\gamma_2):=[\gamma_1\gamma_2; (\gamma_1-1)(\gamma_2-1)].
    \end{equation}
    Then we can write the image $\tau(S_{11})$ as 
    \begin{equation*}
        \tau(S_{11}):=\big\{[\gamma_1\gamma_2; (\gamma_1-1)(\gamma_2-1)]\in\bbR^2: [\gamma_1;\gamma_2]\in S_{11} \big\}.
    \end{equation*}
    Note that we have $\alpha>0$ and $\beta\geq 0$ for every $[\alpha,\beta]\in \tau(S_{11})$. 
    
    On the other hand, observe that \cref{eq:quadratic-eigs} can be written as
    \begin{equation}\label{eq:quadratic-eigs-rewrite}
        z^2 -\Big( (\gamma_1-1)(\gamma_2-1) + 1 -\gamma_1\gamma_2 \big(1-\lambda_i(\bC\bC^\top) \big) \Big)z + (\gamma_1-1)(\gamma_2-1) = 0.
    \end{equation}
    By inspecting \cref{eq:quadratic-alpha-beta} of \cref{lemma:HB} and \cref{eq:quadratic-eigs-rewrite}, we find we can invoke \cref{lemma:HB} with $s=r$ and $\zeta_i=1-\lambda_{r+1-i}(\bC\bC^\top)$ and obtain
    \begin{equation}\label{eq:f>g}
        \begin{split}
            \min_{[\gamma_1;\gamma_2]\in S_{11} } \Big( \max_{i=1,\dots,r} f_i(\gamma_1,\gamma_2) \Big) &= \min_{[\alpha,\beta]\in \tau(S_{11}) } \Big(\max_{i=1,\dots,r} 
 g_i(\alpha,\beta)\Big)  \\ 
            &\geq \min_{\alpha>0,\beta\geq 0 } \Big(\max_{i=1,\dots,r} 
 g_i(\alpha,\beta)\Big) \\ 
            &= \max_{i=1,\dots,r} 
 g_i(\alpha^*,\beta^*) \\ 
            &= \frac{\sqrt{1-\lambda_r(\bC\bC^\top)} - \sqrt{1-\lambda_1(\bC\bC^\top)}  }{\sqrt{1-\lambda_r(\bC\bC^\top)} + \sqrt{1-\lambda_1(\bC\bC^\top)}},
        \end{split} 
    \end{equation}
    where the last two steps follow from \cref{lemma:HB} with $\alpha^*$ and $\beta^*$ defined as
    \begin{equation*}
        \begin{split}
            \alpha^* = \Big( \frac{ 2 }{ \sqrt{\zeta_1} + \sqrt{\zeta_r} } \Big)^2 &= \Bigg( \frac{ 2 }{ \sqrt{1-\lambda_r(\bC\bC^\top)} + \sqrt{1-\lambda_1(\bC\bC^\top)} } \Bigg)^2, \\ 
            \beta^*= \Big( \frac{ \sqrt{\zeta_1} - \sqrt{\zeta_r} }{ \sqrt{\zeta_1} + \sqrt{\zeta_r} } \Big)^2 &= \Bigg( \frac{ \sqrt{1-\lambda_r(\bC\bC^\top)} - \sqrt{1-\lambda_1(\bC\bC^\top)} }{ \sqrt{1-\lambda_r(\bC\bC^\top)} + \sqrt{1-\lambda_1(\bC\bC^\top)} } \Bigg)^2.
        \end{split}
    \end{equation*}
    We will show $[\alpha^*;\beta^*]\in \tau(S_{11})$, as this will prove the inequality in \cref{eq:f>g} is in fact an equality. Indeed, we can always solve the equations
    \begin{equation*}
        \gamma_1\gamma_2=\alpha^*, \quad (\gamma_1-1)(\gamma_2-1) = \beta^*
    \end{equation*}
    for $\gamma_1$ and $\gamma_2$ and obtain two solutions $[\gamma_1^*;\gamma_2^*]\in S_{11}$ and $[\gamma_2^*;\gamma_1^*]\in S_{11}$ (as the reader could verify), where $\gamma_1^*$ and $\gamma_2^*$ are defined as
    \begin{equation*}
        \begin{split}
            \gamma_1^*= \Bigg( \frac{ \sqrt{(1 + \sqrt{\zeta_r} )( 1 + \sqrt{\zeta_1})} + \sqrt{(1 - \sqrt{\zeta_r} )( 1 - \sqrt{\zeta_1})}   }{ \sqrt{\zeta_1}+\sqrt{\zeta_r} } \Bigg)^2, \\ 
            \gamma_2^*= \Bigg( \frac{ \sqrt{(1 + \sqrt{\zeta_r} )( 1 + \sqrt{\zeta_1})} - \sqrt{(1 - \sqrt{\zeta_r} )( 1 - \sqrt{\zeta_1})}   }{ \sqrt{\zeta_1}+\sqrt{\zeta_r} } \Bigg)^2.
        \end{split} 
    \end{equation*}
    Of course $(\gamma_1^*-1)(\gamma_2^*-1) = \beta^*$ by construction, so we have proved 
    \begin{equation*}
        \min_{[\gamma_1;\gamma_2]\in S_{11} } \Big( \max_{i=1,\dots,r} f_i(\gamma_1,\gamma_2) \Big) = \max_{i=1,\dots,r} f_i(\gamma_1^*,\gamma_2^*) = \max_{i=1,\dots,r} f_i(\gamma_2^*,\gamma_1^*) = \sqrt{\beta^*}.
    \end{equation*}
    We can now finish the proof by inspecting \cref{lemma:spectrum-M} and the values of $\gamma_1^*$ and $\gamma_2^*$. Specifically, we have the following three cases.
    \begin{itemize}
        \item If $r=n_1=n_2$, then neither $1-\gamma_1^*$ nor $1-\gamma_2^*$ is an eigenvalue of $\bM(\gamma_1^*,\gamma_2)$, therefore both $[\gamma_1^*;\gamma_2^*]$ and $[\gamma_2^*;\gamma_1^*]$ minimize $\rho\big(\bM(\gamma_1,\gamma_2\big)$ on $S_{11}$.
        \item If $r=n_1<n_2$, then $\bM(\gamma_1^*,\gamma_2^*)$ has one extra eigenvalue $1-\gamma_2^*$. Since $\gamma_1^*\geq \gamma_2^*\geq 1$ and $\beta^*=(\gamma_1^*-1)(\gamma_2^*-1) \geq (\gamma_2^*-1)^2 $, we have 
        \begin{equation*}
            \min_{[\gamma_1;\gamma_2]\in S_{11} } \Big( \max_{i=1,\dots,r} f_i(\gamma_1,\gamma_2) \Big)= \max_{i=1,\dots,r} f_i(\gamma_1^*,\gamma_2^*) = \sqrt{\beta^*}  \geq \gamma_2^*-1,
        \end{equation*}
        and therefore $\rho\big( \bM(\gamma_1,\gamma_2) \big)$ is minimized at $[\gamma_1^*;\gamma_2^*]$  on $S_{11}$.
        \item If $r=n_2 < n_1$, similarly, $\rho\big( \bM(\gamma_1,\gamma_2) \big)$ is minimized at $[\gamma_2^*;\gamma_1^*]$ on $S_{11}$.
    \end{itemize}
\end{proof}

% \begin{proof}[Proof of \cref{theorem:S11-rankdeficient}]
%     Since $r < \min\{n_1,n_2\}$, \cref{lemma:spectrum-M} implies  $1-\gamma_1$ and $1-\gamma_2$ are eigenvalues of $\bM(\gamma_1,\gamma_2)$. For the sake of convenience, define $\lambda_{r+1}(\bC\bC^\top)=0$. We can now observe that $1-\gamma_1$ and $1-\gamma_2$ are the roots of the quadratic equation 
%     \begin{equation}\label{eq:quadratic-eigs-rewrite2}
%         \ z^2 -\Big( (\gamma_1-1)(\gamma_2-1) + 1 -\gamma_1\gamma_2 \big(1- \lambda_{r+1}(\bC\bC^\top) \big) \Big)z + (\gamma_1-1)(\gamma_2-1) = 0 
%     \end{equation}
%     in variable $z$. The roots of \cref{eq:quadratic-eigs-rewrite} and \cref{eq:quadratic-eigs-rewrite2} furnish all eigenvalues of $\bM(\gamma_1,\gamma_2)$. Moreover, comparing \cref{eq:quadratic-eigs-rewrite2}, \cref{eq:quadratic-eigs-rewrite}, and \cref{eq:quadratic-alpha-beta}, we can finish the proof by invoking
%     \cref{lemma:HB} with $s=r+1$, $\xi_i=1-\lambda_{r+2-i}(\bC\bC^\top)$ for $i=1,\dots,r+1$ and proceeding in a similar fashion to the proof of \cref{theorem:S11-fullrank} (e.g., making use of the map $\tau$ in \cref{eq:tau}).
% \end{proof}

\subsection{Proof of \cref{theorem:S11-rankdeficient}}\label{section:alternative-proof}

Here we prove \cref{theorem:S11-rankdeficient}. This is achieved by extending \cref{lemma:distriminant-gamma1>1-gamma1<1} into an analysis of the monotonicity of the determinant $\Delta(\gamma_1,\gamma_2)$ in \cref{eq:equation+discriminant}. We begin with stating and proving \cref{lemma:discriminant,lemma:1quad-gamma12-S11}.
\begin{lemma}\label{lemma:discriminant}
    Let $\phi\in(0,1)$ and $\gamma_1\in[1,2)$ be fixed with $\gamma_1\phi\neq 1$. Consider
    \begin{equation*}
        \begin{split}
            \Delta(\gamma_2) &= (2-\gamma_1-\gamma_2 + \gamma_1\gamma_2\phi)^2 - 4(\gamma_1-1)(\gamma_2-1). 
        \end{split}  
    \end{equation*}
    Let $\xi_1$ and $\xi_2$ be the two roots of $\Delta(\gamma_2)=0$ with $|\xi_1|\leq |\xi_2|$. Then $1\leq \xi_1$. Moreover, we have:
    \begin{itemize}[wide,parsep=-1pt,topsep=0pt]
        \item $\Delta'(\gamma_2)\leq 0$ for every $\gamma_2\in[1,\gamma_1]$.
        \item If $\gamma_1\in [1, \frac{2}{1+\sqrt{1-\phi}}]$, then $\Delta(\gamma_2)\geq 0$ for every $\gamma_2\in[1,\gamma_1]$.
        \item If $\gamma_1\in[\frac{2}{1+\sqrt{1-\phi}},2)$, then $\xi_1\leq \gamma_1$ and 
        \begin{equation*}
            \begin{split}
                \Delta(\gamma_2)&\geq 0, \quad \forall \gamma_2\in[1,\xi_1], \\ 
                \Delta(\gamma_2)&\leq 0, \quad \forall \gamma_2\in[\xi_1,\gamma_1].
            \end{split}
        \end{equation*}
    \end{itemize}
\end{lemma}
\begin{proof}
    Since $\gamma_1\phi\neq 1$, we know $\Delta(\gamma_2)$ is a quadratic function in $\gamma_2$, and we can simplify its expression and calculate its derivative as follows:
    \begin{equation*}
        \begin{split}
            \Delta(\gamma_2) &= (\gamma_1\phi-1)^2\gamma_2^2 + 2(2-\gamma_1)(\gamma_1\phi-1)\gamma_2 -4(\gamma_1-1)\gamma_2 +  (2-\gamma_1)^2 + 4(\gamma_1-1) \\ 
            &= (\gamma_1\phi-1)^2\gamma_2^2 - 2(\gamma_1\phi+ 1 -2\phi)\gamma_1\gamma_2 +  \gamma_1^2,  \\ 
            \Delta'(\gamma_2) &= 2(2-\gamma_1-\gamma_2 + \gamma_1\gamma_2\phi)(\gamma_1\phi-1) - 4(\gamma_1-1) \\ 
            &=  2(\gamma_1\phi-1)^2\gamma_2 - 2( \gamma_1\phi+ 1 -2\phi )\gamma_1.
        \end{split}
    \end{equation*}
    Since $\gamma_1\phi\neq 1$, the minimum of  $\Delta(\gamma_2)$ is attained at $\gamma_2^*:=( \gamma_1\phi+ 1 -2\phi)\gamma_1 /(\gamma_1\phi-1)^2$. Next we verify $\Delta(\gamma_2^*)\leq  0$ and $\gamma_1\leq \gamma_2^*$: We have
    \begin{equation*}
        \begin{split}
            \Delta(\gamma_2^*)\leq  0  
            \Leftrightarrow &\  -\frac{( \gamma_1\phi+ 1 -2\phi)^2\gamma_1^2}{(\gamma_1\phi-1)^2} + \gamma_1^2 \leq  0 \\ 
            \Leftrightarrow &\  (\gamma_1\phi-1)^2 \leq ( \gamma_1\phi+ 1 -2\phi)^2 \\ 
            \Leftrightarrow &\ (1-\phi) (\gamma_1\phi -\phi) \geq 0 \\ 
            \Leftrightarrow &\   (1-\phi) (\gamma_1-1) \geq 0, \\
            \gamma_1\leq \gamma_2^* \Leftrightarrow &\ \gamma_1 \leq \frac{( \gamma_1\phi+ 1 -2\phi)\gamma_1}{(\gamma_1\phi-1)^2}  \\
            \Leftrightarrow &\ (\gamma_1\phi-1)^2 \leq ( \gamma_1\phi+ 1 -2\phi) \\ 
            \Leftrightarrow &\ \gamma_1^2\phi - 3\gamma_1+2\leq 0,
        \end{split}
    \end{equation*}
    and $(1-\phi) (\gamma_1-1) \geq 0$ because $\phi\in(0,1)$ and $\gamma_1\in[1,2)$, while $\gamma_1^2\phi - 3\gamma_1+2\leq 0$ holds because $\gamma_1^2\phi - 3\gamma_1+2\leq (\gamma_1-1)(\gamma_1-2)\leq 0$.

    With the above calculations, we immediately obtain  $\Delta'(\gamma_2)\leq 0$ for every $\gamma_2\in[1,\gamma_1]$. And with $\gamma_1\geq 1$ and $\Delta(1)\geq 0$, we see $\xi_1$ is real-valued, $\xi_1\geq 1$, and
    \begin{equation*}
        \begin{split}
            \gamma_1\leq \xi_1 \Leftrightarrow &\ \Delta(\gamma_1)\geq 0 \\ 
            \Leftrightarrow &\ (\gamma_1\phi-1)^2\gamma_1^2 - 2(\gamma_1\phi+ 1 -2\phi)\gamma_1^2 +  \gamma_1^2 \geq 0 \\ 
            \Leftrightarrow &\ (\gamma_1\phi-1)^2 - 2(\gamma_1\phi+ 1 -2\phi) +  1 \geq 0 \\ 
            \Leftrightarrow &\ \gamma_1^2\phi - 4\gamma_1 +4  \geq 0 \\ 
            \Leftrightarrow &\ (2-\gamma_1)^2 \geq \gamma_1^2(1-\phi) \\ 
            \Leftrightarrow &\ 2-\gamma_1 \geq \gamma_1\sqrt{1-\phi} \\ 
            \Leftrightarrow &\ \gamma_1 \leq  \frac{2}{1+\sqrt{1-\phi}}. 
        \end{split}
    \end{equation*}
    We now draw the other two statements. If $\gamma_1\in [1,\frac{2}{1+\sqrt{1-\phi}}]$, then $\gamma_1\leq \xi_1$, and therefore $\Delta(\gamma_2)\geq 0$ for every $\gamma_2\in[1,\gamma_1]$. On the other hand, assume $\gamma_1\in [\frac{2}{1+\sqrt{1-\phi}},2)$. Then $\xi_1\leq \gamma_1\leq \gamma_2^*$ and $\Delta(\gamma_1)\leq 0$, so $\Delta(\gamma_2)\geq 0$ for each $\gamma_2\in[1,\xi_1]$ and $\Delta(\gamma_2)\leq 0$ for each $\gamma_2\in[\xi_1,\gamma_1]$. 
\end{proof}
\begin{lemma}\label{lemma:1quad-gamma12-S11}
    Let $\phi\in(0,1)$ and $\gamma_1\in[1,2)$ be fixed with $\gamma_1\phi\neq 1$. Consider the following quadratic equation in variable $z$ and its discriminant $\Delta(\gamma_2)$:
    \begin{equation*}
        \begin{split}
            &\ z^2 - (2-\gamma_1-\gamma_2 + \gamma_1\gamma_2\phi)z +  (\gamma_1-1)(\gamma_2-1)=0, \\ 
            &\ \Delta(\gamma_2) := (2-\gamma_1-\gamma_2 + \gamma_1\gamma_2\phi)^2 - 4(\gamma_1-1)(\gamma_2-1).
        \end{split}  
    \end{equation*}
    Let $z_1(\gamma_2),z_2(\gamma_2)$ be the two roots of the above quadratic equation, and let $\xi_1,\xi_2$ be the two roots of $\Delta(\gamma_2)=0$ with $|\xi_1|\leq |\xi_2|$. Then $1\leq \xi_1$. Moreover, with $f(\gamma_2):=\big\{ |z_1(\gamma_2)|, |z_2(\gamma_2)| \big\}$, we have:
    \begin{itemize}
        \item If $\gamma_1\in[1,\frac{2}{1+\sqrt{1-\phi}}]$, then $f(\gamma_2)$ is non-increasing in $[1,\gamma_1]$.
        \item If $\gamma_1\in[\frac{2}{1+\sqrt{1-\phi}}, 2)$, then $\xi_1\leq \gamma_1$, $f(\gamma_2)$ is non-increasing in $[1,\xi_1]$, and
        \begin{equation*}
            f(\gamma_2)=\sqrt{(\gamma_1-1)(\gamma_2-1)}, \quad \forall \gamma_2\in[\xi_1,\gamma_1].
        \end{equation*}
    \end{itemize}
\end{lemma}
\begin{proof}[Proof of \cref{lemma:1quad-gamma12-S11}]
    The case with $\gamma_2\in[\xi_1,\gamma_1]$ follows directly from \cref{lemma:discriminant}. It remains to consider the case $\gamma_2\in [1,\min\{\xi_1,\gamma_1\}]$, where by \cref{lemma:discriminant} we know $\Delta(\gamma_2)\geq 0$ and $z_1(\gamma_2), z_2(\gamma_2)$ are real-valued. Moreover, in this case we have
    \begin{equation*}
        f(\gamma_2)= \frac{|2-\gamma_1-\gamma_2 + \gamma_1\gamma_2\phi| + \sqrt{\Delta(\gamma_2)}}{2},
    \end{equation*}
    and therefore the derivative of $f(\gamma_2)$ is given as
    \begin{equation*}
        f'(\gamma_2) = \begin{cases}
            \frac{-1 + \gamma_1\phi}{2}+ \frac{\Delta'(\gamma_2)}{4\sqrt{\Delta(\gamma_2)}} & \quad  2-\gamma_1-\gamma_2 + \gamma_1\gamma_2\phi \geq 0; \\
            \frac{1 - \gamma_1\phi}{2}+ \frac{\Delta'(\gamma_2)}{4\sqrt{\Delta(\gamma_2)}} & \quad  2-\gamma_1-\gamma_2 + \gamma_1\gamma_2\phi < 0.
        \end{cases}
    \end{equation*}
    Furthermore, since $\Delta'(\gamma_2)\leq 0$ for every $\gamma_2\in[1,\gamma_1]$ (\cref{lemma:discriminant}), to prove $f(\gamma_2)$ is non-increasing in $\gamma_2\in [1,\min\{\xi_1,\gamma_1\}]$ it suffices to show 
    \begin{equation*}
        2\sqrt{\Delta(\gamma_2)} \cdot |\gamma_1\phi-1|  \leq |\Delta'(\gamma_2)|,
    \end{equation*}
    but by \cref{lemma:distriminant-gamma1>1-gamma1<1} this indeed holds. The proof is complete.
\end{proof}
We are in a position to begin proving \cref{theorem:S11-rankdeficient}:
\begin{proof}[Proof of \cref{theorem:S11-rankdeficient}]
    We use the notations in the proof of \cref{theorem:S00-S01-S10}, recalled here for convenience: $z_{i1}(\gamma_1,\gamma_2)$ and $z_{i2}(\gamma_1,\gamma_2)$ denote the two roots of \cref{eq:quadratic-eigs} and  
    \begin{equation*}
        f_i(\gamma_1,\gamma_2):= \max\Big\{ |z_{i1}(\gamma_1,\gamma_2)|, |z_{i2}(\gamma_1,\gamma_2)| \Big\}.
    \end{equation*}
    Define $\Delta(\gamma_1,\gamma_2)$ to be the discriminant of \cref{eq:quadratic-eigs}, that is
    \begin{equation*}
        \Delta(\gamma_1,\gamma_2):= \Big(2-\gamma_1-\gamma_2 + \gamma_1\gamma_2\lambda_1(\bC\bC^\top) \Big)^2 - 4(\gamma_1-1)(\gamma_2-1).
    \end{equation*}
    In light of the symmetry in the spectrum of $\bM(\gamma_1,\gamma_2)$ (\cref{lemma:spectrum-M}), we will assume $\gamma_2\leq \gamma_1$ without loss of generality. We proceed by considering two cases, $\gamma_1 \lambda_1(\bC\bC^\top) \neq 1$ and $\gamma_1 \lambda_1(\bC\bC^\top) = 1$.

    \myparagraph{Case 1: $\gamma_1\lambda_1(\bC\bC^\top)\neq 1$} For this case, we roughly follow the proof logic in \textit{Part 2} of \cref{theorem:S00-S01-S10} with non-trivial modifications to handle the extra difficulty brought by optimization over $S_{11}$. First observe  $\Delta(\gamma_1,\gamma_2)$ is quadratic in $\gamma_2$, and solving $\Delta(\gamma_1,\gamma_2)=0$ for $\gamma_2$ gives two roots, $\xi_1(\gamma_1)$ and $\xi_2(\gamma_1)$. Without loss of generality we assume $|\xi_1|\leq |\xi_2|$. We fix $\gamma_1\in[1,2)$ and proceed by  addressing two sub-cases, $\gamma_1\in [1, \gamma_1^*]$ and $\gamma_1\in [\gamma_1^*, 2)$.

    \myparagraph{Case 1.1: $\gamma_1\in [1, \gamma_1^*]$} In this case, \cref{lemma:1quad-gamma12-S11} suggests that $\gamma_2\mapsto f_1(\gamma_1,\gamma_2)$ is non-increasing in $[1,\gamma_1]$. Furthermore, we have $|z_{11}(\gamma_1,\gamma_1)| \cdot |z_{12}(\gamma_1,\gamma_1)|= (\gamma_1-1)^2$, and therefore
    \begin{equation*}
        f_1(\gamma_1,\gamma_1) \geq \gamma_1 - 1=\gamma_2-1,
    \end{equation*}
    from which it follows that $\gamma_1$ is a global minimizer of 
    \begin{equation}\label{eq:S11-min-gamma2}
        \min_{\gamma_2\in [1, \gamma_1]} \Big( \max\big\{ \gamma_2-1,  f_1(\gamma_1,\gamma_2) \big\} \Big).
    \end{equation}
    We next prove
    \begin{equation}\label{eq:ub-gamma2-hat}
        f_i(\gamma_1,\hat{\gamma}_2) \leq \max\big\{ \hat{\gamma}_2-1,  f_1(\gamma_1,\hat{\gamma}_2) \big\}, \quad \forall i=1,\dots,r,
    \end{equation}
    where $\hat{\gamma}_2$ is a global minimizer of \cref{eq:S11-min-gamma2}, that is   $\hat{\gamma}_2=\gamma_1$. In other words, we will prove $f_i(\gamma_1,\gamma_1)\leq \gamma_1-1$ for every $i=1,\dots,r$. Note that $f_i(\gamma_1,\gamma_1)$ is associated with the quadratic equation
    \begin{equation}\label{eq:gamma2hat-lambdai-S11}
        z^2 - \big(2-2\gamma_1 + \gamma_1^2 \lambda_i(\bC\bC^\top) \big)z + (\gamma_1-1)^2 = 0.
    \end{equation}
    If its two roots $z_{i1}(\gamma_1,\gamma_1)$ and $z_{i2}(\gamma_1,\gamma_1)$ are complex, then $f_i(\gamma_1,\gamma_1)=\sqrt{(\gamma_1-1)^2}=\gamma_1-1$, so \cref{eq:ub-gamma2-hat} holds. Then we consider the case where \cref{eq:gamma2hat-lambdai-S11} admits real-valued roots, which means
    \begin{equation*}
        \begin{split}
            \big(2-2\gamma_1 + \gamma_1^2 \lambda_i(\bC\bC^\top) \big)^2 - 4 (\gamma_1-1)^2  \geq 0 &\Leftrightarrow \gamma_1^2 \lambda_i(\bC\bC^\top) \cdot \big(4-4\gamma_1 + \gamma_1^2 \lambda_i(\bC\bC^\top) \big) \geq 0 \\ &\Leftrightarrow 4-4\gamma_1 + \gamma_1^2 \lambda_i(\bC\bC^\top) \geq 0 \\ 
            &\Rightarrow 2-2\gamma_1 + \gamma_1^2 \lambda_i(\bC\bC^\top) \geq 0. \\ 
        \end{split}
    \end{equation*}
    Here the last step follows from the fact $\gamma_1\geq 1$. But then $2-2\gamma_1 + \gamma_1^2 \lambda_i(\bC\bC^\top)\leq 2-2\gamma_1 + \gamma_1^2 \lambda_1(\bC\bC^\top)$, it must be the case that $f_i(\gamma_1,\gamma_1)\leq f_1(\gamma_1,\gamma_1)$, and we have proved \cref{eq:ub-gamma2-hat}. As a result, we have obtained (see, e.g., \cref{lemma:basic-2->r})
    \begin{equation*}
        \gamma_1\in \argmin_{\gamma_2\in [1,\gamma_1]} \Big( \max\big\{ \gamma_2-1, \max_{i=1,\dots,r} f_i(\gamma_1,\gamma_2) \big\} \Big),
    \end{equation*}
    which, by the definition of $\rho\big(\bM (\gamma_1,\gamma_2) \big)$, leads to 
    \begin{equation*}
        \begin{split}
            \min_{\substack{\gamma_1\in [1,\gamma_1^*]\\ \gamma_2\in[1,\gamma_1] } } \rho\big(\bM (\gamma_1,\gamma_2) \big) &=  \min_{\gamma_1\in[1,\gamma_1^*]} \max\Big\{ \gamma_1-1, \min_{\gamma_2\in [1,\gamma_1]} \Big( \max\big\{ \gamma_2-1, \max_{i=1,\dots,r} f_i(\gamma_1,\gamma_2) \big\} \Big) \Big\}\\
        &=\min_{\gamma_1\in[1,\gamma_1^*]} \max\Big\{ \gamma_1-1,  \max_{i=1,\dots,r} f_i(\gamma_1,\gamma_1)\Big\} \\ 
        &=\min_{\gamma_1\in [1,\gamma_1^*] } \rho\big(\bM (\gamma_1,\gamma_1) \big).
        \end{split} 
    \end{equation*}
    We can now finish the proof for \textit{Case 1.1} by invoking \cref{prop:sr-gamma1=gamma2}.
    
    \myparagraph{Case 1.2: $\gamma_1\in [\gamma_1^*, 2)$} Similarly to \textit{Case 1.1}, from \cref{lemma:1quad-gamma12-S11}, we obtain $1\leq \xi_1(\gamma_1)\leq \gamma_1$ and $\gamma_2\mapsto f_1(\gamma_1,\gamma_2)$ is non-increasing in $[1,\xi_1(\gamma_1)]$. In particular, we have $\Delta(\gamma_1,\xi_1(\gamma_1))=0$, which means $z_{11}\big(\gamma_1,\xi_1(\gamma_1) \big) = z_{12}\big(\gamma_1,\xi_1(\gamma_1) \big)$ and therefore 
    \begin{equation*}
        f_1\big(\gamma_1,\xi_1(\gamma_1)\big)=\sqrt{(\gamma_1-1)\big(\xi_1(\gamma_1)-1\big)} \geq \xi_1(\gamma_1)-1.
    \end{equation*}
    \cref{lemma:1quad-gamma12-S11} implies $\Delta(\gamma_1,\gamma_2)\leq 0$ for every $\gamma_2\in[\xi_1(\gamma_1),\gamma_1]$, so $f(\gamma_1,\gamma_2)=\sqrt{(\gamma_1-1)(\gamma_2-1)}$ . We can now conclude  $\xi_1(\gamma_1)$ is a global minimizer of \cref{eq:S11-min-gamma2}, and moreover, we have
    \begin{equation*}
        \begin{split}
            \min_{ \substack{\gamma_1\in[\gamma_1^*,2)\\ \gamma_2\in [1,\gamma_1] } } \max\Big\{ \gamma_1-1, \gamma_2-1, f_1(\gamma_1,\gamma_2)  \Big\} &= \min_{\gamma_1\in[\gamma_1^*,2)} \max\Big\{ \gamma_1-1, \sqrt{(\gamma_1-1)\big(\xi_1(\gamma_1)-1\big)} \Big\} \\ 
        &= \gamma_1^*-1,
        \end{split}
    \end{equation*}
    where the last equality follows from the fact $\xi(\gamma_1^*)\leq \gamma_1^*$, proved in \cref{lemma:1quad-gamma12-S11}. We have thus obtained a lower bound of $\rho\big(\bM(\gamma_1, \gamma_2 \big)$, that is
    \begin{equation*}
        \begin{split}
            \rho\big(\bM (\gamma_1^*,\gamma_1^*) \big) &= \gamma_1^*-1\\ 
            &= \min_{\gamma_1\in[\gamma_1^*,2),\gamma_2\in [1,\gamma_1]} \max\Big\{ \gamma_1-1, \gamma_2-1, f_1(\gamma_1,\gamma_2)  \Big\} \\
            &\leq \min_{\gamma_1\in[\gamma_1^*,2),\gamma_2\in [1,\gamma_1]} \max\Big\{ \gamma_1-1, \gamma_2-1,  \max_{i=1,\dots,r} f_i(\gamma_1,\gamma_2)   \Big\} \\ 
        &= \min_{\gamma_1\in [\gamma_1^*,2), \gamma_2\in[1,\gamma_1]  } \rho\big(\bM (\gamma_1,\gamma_2) \big).
        \end{split}
    \end{equation*}
    We need to prove the above inequality is actually an equality. To do so, we can show $f_i(\gamma_1^*,\gamma_1^*) \leq \gamma_1^*-1$ for every $i=1,\dots,r$. Note that $f_i\big(\gamma_1^*,\gamma_1^*\big)$ is associated with the equation
    \begin{equation*}%\label{eq:gamma2hat-lambdai-S11-2}
        z^2 - \big(2-2\gamma_1^* + (\gamma_1^*)^2\cdot\lambda_i(\bC\bC^\top) \big)z + (\gamma_1^*-1)^2 = 0,
    \end{equation*}
    but the definition of $\gamma_1^*$ implies $4-4\gamma_1^* + (\gamma_1^*)^2 \lambda_1(\bC\bC^\top) =0$ and therefore
    \begin{equation*}
         4-4\gamma_1^* + (\gamma_1^*)^2 \lambda_i(\bC\bC^\top) \leq 0,\ \ \ \ \forall i=1,\dots,r.
    \end{equation*}
    This means $f_i(\gamma_1^*,\gamma_1^*) = \gamma_1^*-1$ for every $i=1,\dots,r$. We finished \textit{Case 1}.

    \myparagraph{Case 2: $\gamma_1\lambda_1(\bC\bC^\top)=1$} We will show that the minimum of $\rho\big( \bM(\gamma_1,\gamma_2) \big)$ in this case is larger than or equal to its minimum $\gamma_1^*-1$ in \textit{Case 1}, that is,
    \begin{equation}\label{eq:S11-case2>case1}
        \min_{\gamma_1\lambda_1(\bC\bC^\top)=1, \gamma_2\in[1,\gamma_1]} \rho\big( \bM(\gamma_1,\gamma_2) \big) \geq \gamma_1^*-1.
    \end{equation}
    
    We consider two subcases, $\gamma_1\geq 4/3$ and $\gamma_1 < 4/3$. 

    \myparagraph{Case 2.1: $\gamma_1\geq 4/3$} First, recall $\gamma_1^*=\frac{2}{1+ \sqrt{1-\lambda_1(\bC\bC^\top) }}$ and note that
    \begin{equation*}
        \begin{split}
            \gamma_1\geq \frac{4}{3} &\Leftrightarrow \lambda_1(\bC\bC^\top)  \leq \frac{3}{4} \\ 
            &\Leftrightarrow 2 \sqrt{1-\lambda_1(\bC\bC^\top) } \geq 1 \\ 
            &\Leftrightarrow 1 \geq 2 \bigg(1 - \sqrt{1-\lambda_1(\bC\bC^\top) } \bigg) \\ 
            &\Leftrightarrow \frac{1}{\lambda_1(\bC\bC^\top)} \geq \frac{2}{1+ \sqrt{1-\lambda_1(\bC\bC^\top) } }\\ 
            &\Leftrightarrow \gamma_1-1 \geq \gamma_1^*-1 .
        \end{split}
    \end{equation*}
    But we know from \cref{lemma:spectrum-M} that $\gamma_1-1$ is an eigenvalue of $\bM(\gamma_1,\gamma_2)$, therefore \cref{eq:S11-case2>case1} holds.

    \myparagraph{Case 2.2: $\gamma_1< 4/3$} Note that $1\leq \gamma_2\leq \gamma_1<2$, so we can write
    \begin{equation*}
        \begin{split}
            \gamma_1< \frac{4}{3} &\Rightarrow 2-\gamma_1 > 2(\gamma_1-1) \\
            &\Rightarrow (2-\gamma_1)^2 > 4(\gamma_1-1)^2 \\ 
            &\Rightarrow (2-\gamma_1)^2 > 4(\gamma_1-1)(\gamma_2-1).
        \end{split}
    \end{equation*}
    This means the quadratic equation
    \begin{equation*}
        z^2 - (2-\gamma_1) z + (\gamma_1-1)(\gamma_2-1) = 0
    \end{equation*}
    always has two different roots. Moreover, since $\gamma_1\lambda_1(\bC\bC^\top)=1$, \cref{lemma:spectrum-M} implies these two roots are eigenvalues of $\bM(\gamma_1,\gamma_2)$, and this quadratic equation coincides with \cref{eq:quadratic-eigs}, whose roots were denoted previously by $z_{11}(\gamma_1,\gamma_2)$ and $z_{12}(\gamma_1,\gamma_2)$. We can then write
    \begin{equation*}
        \begin{split}
            f_1(\gamma_1, \gamma_2) &=\max\{ z_{11}(\gamma_1,\gamma_2), z_{12}(\gamma_1,\gamma_2) \} \\
            &=\frac{2-\gamma_1}{2} + \frac{\sqrt{(2-\gamma_1)^2 - 4 (\gamma_1-1)(\gamma_2-1)}}{2}.
        \end{split}
    \end{equation*}
    Therefore $\gamma_2\mapsto f_1(\gamma_1,\gamma_2)$ is decreasing in $[1,\gamma_1]$, and we have thus obtained
    \begin{equation*}
        f(\gamma_1,\gamma_1)= \min_{\gamma_2\in[1,\gamma_1]} f_1(\gamma_1, \gamma_2). 
    \end{equation*}
    We finish the proof by observing
    \begin{equation*}
        \min_{\gamma_2\in[1,\gamma_1]} \rho\big( \bM(\gamma_1,\gamma_2) \big) \geq \max\big\{ \gamma_1-1, f_1(\gamma_1, \gamma_1)   \big\} \geq \gamma^*_1-1,
    \end{equation*}
    where the last ineequality follows from \cref{prop:sr-gamma1=gamma2} with the special case $r=1$.
\end{proof}

\section{Extra Discussions on Gauss–Seidel}\label{section:Gauss-Seidel}
The classic Gauss-Seidel method applied to normal equations $\bA^\top \bA \bx = \bA^\top \by$ is precisely a BGD method applied to \cref{eq:LS} with $n$ blocks, each of size $1$, and block $i$ is with stepsize $1/a_{ii}$, where $a_{ii}$ is the $i$-th diagonal entry of $\bA^\top \bA$. This stepsize rule also corresponds to block exact minimization. Again, this choice of stepsizes is sub-optimal. To see this,  consider a simple case where $\bA$ has only two columns $\ba_1$ and $\ba_2$ ($n=2$). By \cref{lemma:error-decrease}, the iterates of \ref{eq:BGD} satisfy $\bx^+ - \bx^*= \bM(\gamma_1,\gamma_2) \cdot (\bx - \bx^*)$ with $\bM(\gamma_1,\gamma_2)$ defined as
\begin{equation*}
    \bM(\gamma_1,\gamma_2):=\begin{bmatrix}
        1-\gamma_1 \ba_1^\top \ba_1 & -\gamma_1 \ba_1^\top \ba_2 \\ 
        -\gamma_2 (1-\gamma_1 \ba_1^\top \ba_1) \ba_2^\top \ba_1 & \gamma_1\gamma_2(\ba_1^\top \ba_2)^2 + 1 - \gamma_2 \ba_2^\top \ba_2 
    \end{bmatrix} .
\end{equation*}
% \begin{equation*}
%     \begin{split}
%         \bM(\gamma_1,\gamma_2) &= \Bigg(\bI - \begin{bmatrix}
%             0 & 0 \\ 
%             \gamma_2 \ba_2^\top \ba_1 & \gamma_2 \ba_2^\top \ba_2 
%         \end{bmatrix} \Bigg) \Bigg(\bI - \begin{bmatrix}
%             \gamma_1 \ba_1^\top \ba_1 & \gamma_1 \ba_1^\top \ba_2 \\ 
%             0 & 0
%         \end{bmatrix} \Bigg) \\ 
%         &= \begin{bmatrix}
%             1-\gamma_1 \ba_1^\top \ba_1 & -\gamma_1 \ba_1^\top \ba_2 \\ 
%             -\gamma_2 (1-\gamma_1 \ba_1^\top \ba_1) \ba_2^\top \ba_1 & \gamma_1\gamma_2(\ba_1^\top \ba_2)^2 + 1 - \gamma_2 \ba_2^\top \ba_2 
%         \end{bmatrix} 
%     \end{split}
% \end{equation*}
The Gauss-Seidel method takes stepsizes $\gamma_1=1/(\ba_1^\top \ba_1),\gamma_2=1/(\ba_2^\top \ba_2)$, and one verifies that 
\begin{equation*}
    \rho \bigg( \bM\Big(  \frac{1}{\ba_1^\top \ba_1}, \frac{1}{\ba_2^\top \ba_2} \Big) \bigg)  = (\ba_1^\top \ba_2)^2 / \big( (\ba_1^\top \ba_1) \cdot (\ba_2^\top \ba_2) \big).
\end{equation*}
% the spectral radius of $\bM\Big(  \frac{1}{\ba_1^\top \ba_1}, \frac{1}{\ba_2^\top \ba_2} \Big)$ is $(\ba_1^\top \ba_2)^2 / \big( (\ba_1^\top \ba_1) \cdot (\ba_2^\top \ba_2) \big)$
% and then
% \begin{equation*}
%     \bM\Big(  \frac{1}{\ba_1^\top \ba_1}, \frac{1}{\ba_2^\top \ba_2} \Big) = \begin{bmatrix}
%             0& -\gamma_1 \ba_1^\top \ba_2 / (\ba_1^\top \ba_1) \\ 
%             0 & (\ba_1^\top \ba_2)^2 / \big( (\ba_1^\top \ba_1) \cdot (\ba_2^\top \ba_2) \big)
%         \end{bmatrix}, 
% \end{equation*}
% whose spectral radius is $(\ba_1^\top \ba_2)^2 / \big( (\ba_1^\top \ba_1) \cdot (\ba_2^\top \ba_2) \big)$. 

On the other hand, it follows from \cref{remark:sp=0} that the minimum spectral radius in this case is zero; indeed, with stepsizes $\gamma_1=\frac{1}{\ba_1^\top \ba_1}, \gamma_2=\frac{\ba_1^\top \ba_1}{(\ba_1^\top \ba_1)\cdot (\ba_2^\top \ba_2) - (\ba_1^\top \ba_2)^2}$ we have 
\begin{equation*}%\label{eq:optimal-n=2}
    \bM(  \gamma_1,\gamma_2 ) = \begin{bmatrix}
            0& - \ba_1^\top \ba_2 / (\ba_1^\top \ba_1) \\ 
            0 & 0
        \end{bmatrix},
\end{equation*}
whose eigenvalues are zero. Note that this reasoning only provides a counter-example showing the stepsize sub-optimality of the Gauss-Seidel method, and we are none the wiser: Finding closed-form optimal stepsizes of \ref{eq:BGD} in its full generality has remained non-trivial for $n>2$.

% \subsubsection{Put It Together}
% The hard work has been done in \cref{theorem:S00-S01-S10,theorem:S11-fullrank,theorem:S11-rankdeficient}. Indeed, comparing the minimum spectral radii on the four regions ($S_{00},S_{01},S_{10},S_{11}$) in \cref{theorem:S00-S01-S10,theorem:S11-fullrank,theorem:S11-rankdeficient}, we find  the optimal stepsizes on $S_{11}$, given by either \cref{theorem:S11-fullrank} or \ref{theorem:S11-rankdeficient}. A solution to \cref{eq:minimize-sr} and a closed-form expression for $\rho_{\textnormal{BGD}}^*$ now ensues:

% For example, in the case $r>1$, if $\bC$ is full rank, then the optimal stepsizes and spectral radius to \cref{eq:minimize-sr} are given by \cref{theorem:S11-fullrank}; if $\bC$ is rank-deficient, they are given by \cref{theorem:S11-rankdeficient}. 

\section{Proof for Theorems Related to Generalized Alternating Projection}
\begin{proof}[Proof of \cref{eq:singular-val=angle}]
    It follows directly from Exercise 2.8 of \citet{Vidal-GPCA2016}, and either Property 2.1 of \citet{Zhu-JNM2013} or Theorem 2.7 of \citet{Knyazev-SIAM-J-MAA2007}.
\end{proof}

\begin{proof}[Proof of \cref{lemma:BCD-AP-matrix}]
    With $C:=A_2^\top A_1$, we have 
    \begin{equation*}
        \begin{split}
            A^\top A &= \begin{bmatrix}
            \bI_1 & C^\top \\ 
            C & \bI_2 
        \end{bmatrix}, \\
            A_2^\top A &= \begin{bmatrix}
            % \bI_1 & C^\top \\ 
            C & \bI_2 
        \end{bmatrix}, \\ 
        A_1^\top A &= \begin{bmatrix}
            \bI_1 & C^\top 
            % C & \bI_2 
        \end{bmatrix}. 
        \end{split}
    \end{equation*}
    we need to prove
    \begin{equation*}
        A^\top ( \bI - \gamma_2 A_2 A_2^\top) ( \bI - \gamma_1 A_1 A_1^\top)A = \begin{bmatrix}
        \bI_1 & C^\top \\ 
        C & \bI_2
    \end{bmatrix}  \left(\bI - \begin{bmatrix}
            0 & 0 \\ 
            \gamma_2 C & \gamma_2 \bI_2 
        \end{bmatrix} \Bigg) \Bigg(\bI - \begin{bmatrix}
            \gamma_1 \bI_1 & \gamma_1 C^\top  \\ 
            0 & 0
        \end{bmatrix} \right).
    \end{equation*}
We do so by first simplifying the two terms separately. The left-hand side can be written as
\begin{equation*}
    \begin{split}
        A^\top ( \bI - \gamma_2 A_2 A_2^\top) ( \bI - \gamma_1 A_1 A_1^\top)A &= A^\top A - \gamma_2 A^\top A_2 A_2^\top A - \gamma_1 A^\top A_1A_1^\top A + \gamma_1\gamma_2 A^\top A_2 C A_1^\top A \\
        &=\begin{bmatrix}
            \bI_1 & C^\top \\ 
            C & \bI_2 
        \end{bmatrix}  - \gamma_2 \begin{bmatrix}
            C^\top C & C^\top \\ 
            C & \bI_2
        \end{bmatrix} - \gamma_1 \begin{bmatrix}
            \bI_1 & C^\top \\ 
            C & C C^\top 
        \end{bmatrix} + \gamma_1\gamma_2 \begin{bmatrix}
            C^\top \\
            \bI_2
        \end{bmatrix} C \begin{bmatrix}
            \bI_1 & C^\top 
            % C & \bI_2 
        \end{bmatrix} ,
    \end{split}
\end{equation*}
and the right-hand side can be written as
\begin{equation*}
    \begin{split}
        \begin{bmatrix}
        \bI_1 & C^\top \\ 
        C & \bI_2
    \end{bmatrix}- \gamma_2 \begin{bmatrix}
        \bI_1 & C^\top \\ 
        C & \bI_2
    \end{bmatrix} \begin{bmatrix}
            0 & 0 \\ 
             C &  \bI_2 
        \end{bmatrix} - \gamma_1 \begin{bmatrix}
        \bI_1 & C^\top \\ 
        C & \bI_2
    \end{bmatrix} \begin{bmatrix}
            \bI_1 &  C^\top \\ 
            0 & 0
        \end{bmatrix} + \gamma_1 \gamma_2 \begin{bmatrix}
        \bI_1 & C^\top \\ 
        C & \bI_2
    \end{bmatrix}  \begin{bmatrix}
            0 & 0 \\ 
             C  &   CC^\top 
        \end{bmatrix}.
    \end{split}
\end{equation*}
Then, by inspection, it suffices to prove the following three matrix equations:
\begin{equation*}
    \begin{split}
         \begin{bmatrix}
            C^\top C & C^\top \\ 
            C & \bI_2
        \end{bmatrix}  &= \begin{bmatrix}
        \bI_1 & C^\top \\ 
        C & \bI_2
    \end{bmatrix} \begin{bmatrix}
            0 & 0 \\ 
             C &  \bI_2 
        \end{bmatrix}, \\
          \begin{bmatrix}
            \bI_1 & C^\top \\ 
            C & C C^\top 
        \end{bmatrix}  &= \begin{bmatrix}
        \bI_1 & C^\top \\ 
        C & \bI_2
    \end{bmatrix} \begin{bmatrix}
            \bI_1 &  C^\top \\ 
            0 & 0
        \end{bmatrix}, \\
         \begin{bmatrix}
            C^\top \\
            \bI_2
        \end{bmatrix} C \begin{bmatrix}
            \bI_1 & C^\top 
            % C & \bI_2 
        \end{bmatrix} &= \begin{bmatrix}
        \bI_1 & C^\top \\ 
        C & \bI_2
    \end{bmatrix}  \begin{bmatrix}
            0 & 0 \\ 
             C  &   CC^\top 
        \end{bmatrix}. \\
    \end{split}
\end{equation*}
These can be easily verified, therefore the proof is complete.
\end{proof}

\begin{proof}[Proof of \cref{corollary:GAP}]
    By \cref{corollary:main-result} we have
    \begin{equation*}
        \min_{\gamma_1>0,\gamma_2>0} \rho\big( M(\gamma_1,\gamma_2) \big) = \begin{cases}
           \frac{\sqrt{1-\lambda_r(CC^\top)} - \sqrt{1-\lambda_1(CC^\top)}  }{\sqrt{1-\lambda_r(CC^\top)} + \sqrt{1-\lambda_1(CC^\top)}}  & r = \min\{n_1,n_2\}; \\
           \frac{1 - \sqrt{1-\lambda_1(CC^\top)}  }{1 + \sqrt{1-\lambda_1(CC^\top)}}. & r < \min\{n_1,n_2\}.
        \end{cases}
    \end{equation*}
    By \cref{eq:singular-val=angle} we have $\sqrt{1-\lambda_r(CC^\top)}=\sin(\theta_r)$ and $\sqrt{1-\lambda_1(CC^\top)}=\sin(\theta_1)$. The proof is then complete.
\end{proof}
\begin{remark}\label{remark:opt-stepsizes-AP}
    Here we describe the corresponding stepsizes $\gamma_1^*,\gamma_2^*$ that attain the above minimum. This is derived from \cref{theorem:S11-rankdeficient} and the proof of \cref{theorem:S11-fullrank}.
    \begin{itemize}
        \item If $r=\min\{n_1,n_2\}$, then we have 
        \begin{equation*}
        \begin{split}
            \gamma_1^*= \Bigg( \frac{ \sqrt{(1 + \sin(\theta_1) )( 1 + \sin(\theta_r))} + \sqrt{(1 - \sin(\theta_1) )( 1 - \sin(\theta_r))}   }{ \sin(\theta_r)+\sin(\theta_1) } \Bigg)^2, \\ 
            \gamma_2^*= \Bigg( \frac{ \sqrt{(1 + \sin(\theta_1) )( 1 + \sin(\theta_r))} - \sqrt{(1 - \sin(\theta_1) )( 1 - \sin(\theta_r))}   }{ \sin(\theta_r)+\sin(\theta_1) } \Bigg)^2,
        \end{split} 
    \end{equation*}
    or 
    \begin{equation*}
        \begin{split}
            \gamma_1^*= \Bigg( \frac{ \sqrt{(1 + \sin(\theta_1) )( 1 + \sin(\theta_r))} - \sqrt{(1 - \sin(\theta_1) )( 1 - \sin(\theta_r))}   }{ \sin(\theta_r)+\sin(\theta_1) } \Bigg)^2, \\ 
            \gamma_2^*= \Bigg( \frac{ \sqrt{(1 + \sin(\theta_1) )( 1 + \sin(\theta_r))} + \sqrt{(1 - \sin(\theta_1) )( 1 - \sin(\theta_r))}   }{ \sin(\theta_r)+\sin(\theta_1) } \Bigg)^2.
        \end{split} 
    \end{equation*}
    In other words, in this case there are two sets of optimal parameters and they are symmetric. These stepsizes would coincide with the stepsize rule of \citet{Falt-CDC2017}  if  $\theta_r$ were equal to $\pi/2$. As justified, for generic subspaces we have $\theta_r=\pi/2$ with probability $0$.
        \item if $r<\min\{n_1,n_2\}$, then 
        \begin{equation*}
            \gamma_1^*=\gamma_2^*=\frac{2}{1+\sin(\theta_1)},
        \end{equation*}
        which coincides with the stepsize rule of \citet{Falt-CDC2017}; see \cref{table:compare-methods}.
    \end{itemize}
\end{remark}

\section{Code}\label{section:code}
\subsection{Experiments on Least-Squares With \cref{assumption:BWO} (\cref{fig:sr-convergence})}\label{subsection:code-w.A1}

{ 
\small
\begin{verbatim}
clc; clear all; format longG

cond_num = 1e5;
n_iter = 5000;

num_trials = 1;
m = 1000;
n1 = 300;
n2 = 500;
n = n1 + n2;
noise_level = 0.01;

errorsBGD = zeros(num_trials, n_iter+1);
errorsGD = zeros(num_trials, n_iter+1);
errorsHB = zeros(num_trials, n_iter+1); 

times_BGD = zeros(num_trials, 1);

% stepsizes of methods
global gammaGD alphaHB betaHB gamma1 gamma2;

for t=1:num_trials
    %% generate data
    [A, y, x, lambdamaxCC, lambdaminCC] = gen_data(m, n1, n2, noise_level, cond_num);
    
    %% calculae stepsizes
    minC = sqrt(1-lambdamaxCC);
    maxC = sqrt(1-lambdaminCC);

    lambdamaxAA = 1 + sqrt(lambdamaxCC);
    lambdaminAA = 1 - sqrt(lambdamaxCC);
    
    lambdamaxAA/lambdaminAA
    
    gammaGD = 2 / (lambdamaxAA + lambdaminAA);

    alphaHB = 4 / ( sqrt(lambdamaxAA) + sqrt(lambdaminAA) )^2;
    betaHB = ( sqrt(lambdamaxAA) - sqrt(lambdaminAA) )^2;
    betaHB = betaHB / ( sqrt(lambdamaxAA) + sqrt(lambdaminAA) )^2;

    
    gamma1 = (1+maxC)*(1+minC) / (maxC + minC)^2 ;
    gamma2 = gamma1;

    %% run algorithm
    [BGD_x, BGD_iters, BGD_time] = BGD(A, y, n1, n2, n_iter);
    errorsBGD(t,:) = vecnorm(BGD_iters - x, 2, 1);
    
    [HB_x, HB_iters, HB_time] = HB(A, y, n_iter);
    errorsHB(t,:) = vecnorm(HB_iters - x, 2, 1);

    [GD_x, GD_iterates, GD_time] = GD(A, y, n_iter);
    errorsGD(t,:) = vecnorm(GD_iterates - x, 2, 1);
end

mean_errorsBGD = mean(errorsBGD,1);
mean_errorsGD = mean(errorsGD,1);
mean_errorsHB = mean(errorsHB,1);


[mean_errorsBGD; mean_errorsHB; mean_errorsGD]

function [BGD_x, BGD_iters, BGD_time] = BGD(A, y, n1, n2, n_iter)
    tic;

    global gamma1 gamma2;
    
    c1 = 1-gamma1; c2 = 1-gamma2;
    
    A1 = A(:, 1:n1);
    A2 = A(:, n1+1:end);
    
    C = A2'*A1;
    Ay = A'*y;

    
    BGD_x = zeros(n1+n2,1);
    BGD_iters = BGD_x;
    for i=1:n_iter
        BGD_x(1:n1) = c1*BGD_x(1:n1) - gamma1*(C'*BGD_x(n1+1:end) - Ay(1:n1));
        BGD_x(n1+1:end) = c2*BGD_x(n1+1:end) - gamma2*(C*BGD_x(1:n1) - Ay(n1+1:end));
        
        BGD_iters = [BGD_iters BGD_x];
    end
    BGD_time = toc;
end

function [GD_x, GD_iterates, GD_time] = GD(A, y, n_iter)
    tic;
    
    global gammaGD;
    
    [m,n] = size(A);
    
    AA = A'*A;
    Ay = A'*y;
    
    GD_x = zeros(n,1);
    GD_iterates = GD_x;
    for i=1:n_iter
        GD_x = GD_x - gammaGD * (AA*GD_x - Ay);
        
        GD_iterates = [GD_iterates GD_x];
    end
    
    GD_time = toc;
end

function [HB_x, HB_iters, HB_times] = HB(A, y, n_iter)
    tic;
    
    global alphaHB betaHB;
    
    [m,n] = size(A);

    AA = A'*A;
    Ay = A'*y;
       
    HB_x = zeros(n,1); HB_x_old = HB_x;
    HB_iters = HB_x;
    HB_times = [];
    for i=1:n_iter
        HB_x_new = HB_x - alphaHB * (AA*HB_x - Ay) + betaHB*(HB_x - HB_x_old);
        
        HB_x_old = HB_x; HB_x = HB_x_new;
        
        HB_iters = [HB_iters HB_x];
        HB_times = [HB_times toc];
    end
    
end


function [A, y, x, lambdamaxCC, lambdaminCC]= gen_data(m, n1, n2, noise_level, cond_num)
    % cond_num = ( 1+sqrt(lambdamaxC) ) / ( 1-sqrt(lambdamaxC) )
    
    lambdamaxC = (cond_num - 1) / (cond_num + 1);
    
    %% generate 
    
    [C, lambdamaxCC, lambdaminCC] = gen_mat_random_bounded_svs(n2, n1, lambdamaxC);
   
    [U, ~, ~] = svd(randn(m-n2, n1), 'econ');
    
    scale = vecnorm(C, 2, 1);
    U = U .* sqrt(1- scale.^2);
    
    A1 = [C; U];             % by construction, A1 is orthogonal (A1'*A1 = identity)
    A2 = [eye(n2,n2); zeros(m-n2, n2)];
    
    [U1, ~, ~] = svd(randn(n1,n1)); A1 = A1 * U1;
    [U2, ~, ~] = svd(randn(n2,n2)); A2 = A2 * U2;
    
    A = [A1 A2];
    
    x = randn(n1+n2, 1);
    
    noise = randn(m,1); noise = noise / norm(noise);
    
    y = A*x + noise_level* noise;

    x = A \ y;
end

function [C, lambdamaxCC, lambdaminCC] = gen_mat_random_bounded_svs(n2, n1, L)
    % generate a n2xn1 matrix C such that
    % the maximum and minimum singlar value of C are U and L respectively
    %    the rest singular values are randomly chosen in [0,L2]
    
    num = min(n2,n1) - 1;
    
    svs = L*rand(num,1);
    
    
    svs = sort([svs; L], 'descend');
    
    lambdamaxCC = max(svs)*max(svs);
    lambdaminCC = min(svs)*min(svs);
    
    
    C = diag(svs);
    
    if n1 > n2
        C = [C zeros(n2, n1-n2)];
    elseif n1 < n2
        C = [C; zeros(n2-n1, n1)];
    end
end

\end{verbatim}

}

\subsection{Experiments on Least-Squares (Without \cref{assumption:BWO})}\label{subsection:code-w.o.A1}
{\normalsize Below is the code implementing the heavy ball method in comparison to our proposal, block-wise QR orthogonalization followed by \ref{eq:BGD} with the optimal stepsizes. As analyzed in the main paper, the code makes the point that the proposed method would converge faster than the heavy ball method under the current problem configuration. That being said, it should also be noted that this is still slower than practical least-squares solvers, e.g., QR orthogonalization followed by back substitution to solve an upper triangular system, randomized methods, (preconditioned) conjugate gradient methods, or the MATLAB backslash solver. }

{\small
\begin{verbatim}
clc; clear all; format longG

%% setup
cond_num = 1e5;
n_iter = 5000;

num_trials = 1;
m = 1000;
n1 = 300;
n2 = 500;
n = n1 + n2;
noise_level = 0.01;


%% run experiments
errorsBGD = zeros(num_trials, n_iter+1);
errorsHB = zeros(num_trials, n_iter+1); 

for t=1:num_trials
    %% generate data
    [A, y, x] = gen_data(m, n, cond_num, noise_level);
    %% calculate stepsizes

    %% run algorithm
    [BGD_x, Q1, Q2, R1, R2, BGD_iters, qrtime, BGD_times] = BGD(A, y, n1, n2, n_iter);
    % the LS solution [x1, x2] can be found by backward substitution using R1 and R2.

    errorsBGD(t,:) = Q1 * BGD_iters(1:n1, :) + Q2 * BGD_iters(n1+1:end, :) - y
    errorsBGD(t,:) = vecnorm(errorsBGD(t,:), 2, 1); 
    
    [HB_x, HB_iters, HB_times] = HB(A, y, n_iter);
    % errorsHB(t,:) = vecnorm(HB_iters - x, 2, 1);
    errorsHB(t,:) = vecnorm(A*HB_iters - y, 2, 1);
end

mean_errorsBGD = mean(errorsBGD,1);
mean_errorsHB = mean(errorsHB,1);

BGD_qrtimes = qrtime + BGD_times;

figure(1);
plot(BGD_times, 'LineWidth', 2);
hold on;
plot(BGD_times + qrtime, 'LineWidth', 2);

plot(HB_times, 'LineWidth', 2);
legend("BGD", "BGD+Ortho", "HB")

set(gca, 'FontSize', 18)

figure(2);
plot(mean_errorsBGD, 'LineWidth', 2);
hold on;

plot(mean_errorsHB, 'LineWidth', 2);
legend("BGD+Ortho", "HB")

set(gca, 'YScale', 'log')
set(gca, 'FontSize', 18)



function [BGD_x, Q1, Q2, R1, R2, BGD_iters, qrtime, BGD_times] = BGD(A, y, n1, n2, n_iter)
    
    tic; % orthogonalize
    [Q1, R1] = qr(A(:, 1:n1), 0);
    [Q2, R2] = qr(A(:, n1+1:end), 0);
    qrtime = toc;

    % calculate stepsizes. We don't count the running times here.
    C = Q2'*Q1;
    
    s = svd(C, 0);
    
    gamma1 = 2/(1 + sqrt(1-max(s)^2)); gamma2 = gamma1;
    
    tic;  % run the algorithm
    c1 = 1-gamma1; c2 = 1-gamma2;
    
    Qy = [Q1 Q2]'*y;
    
    BGD_x = zeros(n1+n2,1);
    BGD_iters = BGD_x;
    BGD_times = [];
    for i=1:n_iter

        BGD_x(1:n1) = c1*BGD_x(1:n1) - gamma1*(C'*BGD_x(n1+1:end) - Qy(1:n1));
        BGD_x(n1+1:end) = c2*BGD_x(n1+1:end) - gamma2*(C*BGD_x(1:n1) - Qy(n1+1:end));

        BGD_iters = [BGD_iters BGD_x];
        BGD_times = [BGD_times toc];
    end
end

function [HB_x, HB_iters, HB_times] = HB(A, y, n_iter)
    % calculate stepsizes for HB. We don't count the running times here.
    e = eig(A'*A);
    lambdamaxAA = max(e);
    lambdaminAA = min(e);
    
    alphaHB = 4 / ( sqrt(lambdamaxAA) + sqrt(lambdaminAA) )^2;
    betaHB = ( sqrt(lambdamaxAA) - sqrt(lambdaminAA) )^2;
    betaHB = betaHB / ( sqrt(lambdamaxAA) + sqrt(lambdaminAA) )^2;
    
    tic; % run the algorithm
    [m,n] = size(A);

    AA = A'*A;
    Ay = A'*y;
       
    HB_x = zeros(n,1); HB_x_old = HB_x;
    HB_iters = HB_x;
    HB_times = [];
    for i=1:n_iter
        
        HB_x_new = HB_x - alphaHB * (AA*HB_x - Ay) + betaHB*(HB_x - HB_x_old);
        
        
        HB_x_old = HB_x; HB_x = HB_x_new;
        
        
        HB_iters = [HB_iters HB_x];
        HB_times = [HB_times toc];
    end
    
end

function [A, y, x] = gen_data(m, n, cond_num, noise_level)

    A = randn(m,n);
    
    cond_numA = sqrt(cond_num);
    
    [U, ~, V] = svd(A, 'econ');
    s = 1 + (cond_numA - 1).*rand(n-2,1);
  
    
    s = sort([s; 1; cond_numA], 'descend');
    
    A = U * diag(s) * V';
    
    x = randn(n,1);
    
    noise = randn(m,1); noise = noise / norm(noise);
    
    y = A*x + noise_level* noise; % just to make sure that the minimum loss is not zero

    x = A \ y;
end

\end{verbatim}

}

\end{document}